\documentclass[11pt, a4paper]{article}
\usepackage{amsmath}
\usepackage{wasysym}
\usepackage{latexsym}
\usepackage{amsfonts}
\usepackage{mathrsfs}
\usepackage{amssymb}
\usepackage{ifsym}
\usepackage{dsfont}
\usepackage[all]{xy}
\usepackage{authblk}

\newtheorem{Definitions1}{Definition}[section]

\newtheorem{Theorems1}{Theorem}[section]

\newtheorem{Coroll1}[Theorems1]{Corollary}
\newtheorem{Lemma1}[Theorems1]{Lemma}
\newtheorem{Remark1}[Theorems1]{Remark}

\newtheorem{Quest1}{Question}[section]

\newenvironment{proof}[1][Proof]{\begin{trivlist}
\item[\hskip \labelsep {\bfseries #1}]}{\end{trivlist}}

\newcommand{\qed}{\nobreak \ifvmode \relax \else
      \ifdim\lastskip<1.5em \hskip-\lastskip
      \hskip1.5em plus0em minus0.5em \fi \nobreak
      \vrule height0.75em width0.5em depth0.25em\fi}

\begin{document}
\title{Largest initial segments pointwise fixed by automorphisms of models of set theory}
\author[1]{Ali Enayat}
\author[2]{Matt Kaufmann} 
\author[1]{Zachiri McKenzie}
\affil[1]{Department of Philosophy, Linguistics and Theory of Science\\
Gothenburg University\\
\texttt{\{ali.enayat,zachiri.mckenzie\}@gu.se}}
\affil[2]{Department of Computer Science, University of Texas at Austin\\
\texttt{kaufmann@cs.utexas.edu}}
\maketitle

\begin{abstract}
Given a model $\mathcal{M}$ of set theory, and a nontrivial automorphism $j$
of $\mathcal{M}$, let $\mathcal{I}_{\mathrm{fix}}(j)$ be the submodel of $\mathcal{M}$ whose universe consists of elements $m$ of $\mathcal{M}$ such
that $j(x)=x$ for every $x$ in the transitive closure of $m$ (where the
transitive closure of $m$ is computed within $\mathcal{M}$). Here we study
the class $\mathcal{C}$ of structures of the form $\mathcal{I}_{\mathrm{fix}}(j)$, where the ambient model $\mathcal{M}$ satisfies a frugal yet robust fragment of $\mathrm{ZFC}$ known as $\mathrm{MOST}$, and $j(m)=m$ whenever $m$\ is\ a finite\ ordinal in the sense of $\mathcal{M}.$ Our main achievement is the calculation of the theory of $\mathcal{C}$ as precisely $\mathrm{MOST+\Delta }_{0}^{\mathcal{P}}$-$\mathrm{Collection.}$ The following
theorems encapsulate our principal results:\textrm{\medskip }

\noindent \textbf{Theorem A.} \textit{Every structure in} $\mathcal{C}$
\textit{satisfies} $\mathrm{MOST+\Delta }_{0}^{\mathcal{P}}$-$\mathrm{%
Collection}$.\textrm{\medskip }

\noindent \textbf{Theorem B. }\textit{Each of the following three conditions
is sufficient for a countable structure} $\mathcal{N}$ \textit{to be in }$\mathcal{C}$:$\mathrm{\smallskip }$

\noindent \textbf{(a)} $\mathcal{N}$ \textit{is a transitive model of} $\mathrm{MOST+\Delta }_{0}^{\mathcal{P}}$-$\mathrm{Collection.\smallskip }$

\noindent \textbf{(b)} $\mathcal{N}$ \textit{is a recursively saturated
model of} $\mathrm{MOST+\Delta}_{0}^{\mathcal{P}}$-$\mathrm{Collection.\smallskip }$

\noindent \textbf{(c)} $\mathcal{N}$ \textit{is a model of }$\mathrm{ZFC}$.%
\textrm{\medskip }

\noindent \textbf{Theorem C. }\textit{Suppose }$\mathcal{M}$\textit{\ is a
countable recursively saturated model of }$\mathrm{ZFC}$\textit{\ and }$I$ \textit{\ is a proper initial segment of }$\mathrm{Ord}^{\mathcal{M}}$
\textit{\ that is closed under exponentiation and contains $\omega^\mathcal{M}$. There is a group embedding} $j\longmapsto \check{j}$ \textit{from }$\mathrm{Aut}(\mathbb{Q})$\textit{\
into }$\mathrm{Aut}(\mathcal{M})$\textit{\ such that }$I$ \textit{is the
longest initial segment of} $\mathrm{Ord}^{\mathcal{M}}$ \textit{that is
pointwise fixed by} $\check{j}$ \textit{for every nontrivial }$j\in
\mathrm{Aut}(\mathbb{Q}).$ \textrm{\medskip }

\noindent In Theorem C, $\mathrm{Aut}(X)$ is the group of automorphisms of
the structure $X$, and $\mathbb{Q}$ is the ordered set of rationals.
\end{abstract}

\section[Introduction]{Introduction}

We study automorphisms of models of set theory, focusing on structures of
the form $\mathcal{I}_{\mathrm{fix}}(j)$, where $j$ is a nontrivial
automorphism of a model $\mathcal{M}$ of $\mathrm{MOST}$\ (Mostowski Set
Theory) such that $j(m)=m$ whenever $m$\ is\ a finite\ ordinal in the sense
of $\mathcal{M}$, and $\mathcal{I}_{\mathrm{fix}}(j)$ is the submodel of $%
\mathcal{M}$ whose universe consists of elements $m$ in $M$ such that $j(x)=x
$ for every $x$ in the transitive closure of $m$ (where the transitive
closure is calculated within $\mathcal{M}$). $\mathrm{MOST}$ is a frugal yet
robust fragment of $\mathrm{ZFC}$, introduced and studied in Mathias'
majestic paper \cite{mat01} (see Definition 2.4).\medskip

A principal source of motivation for our work is to be found in the study of
automorphisms of models of $\mathrm{ZF}$'s sister theory $\mathrm{PA}$
(Peano arithmetic). For example: by a theorem of Smory\'{n}ski (\cite{smo82}
, \cite[Theorem 8.4.2]{ks06}), if $\mathcal{I}$ is submodel of a
countable recursively saturated model of $\mathrm{PA}$ whose elements form a
proper initial segment of $\mathcal{M}$, and $\mathcal{I}$ is closed under
the exponential function of $\mathcal{M}$, then there is an automorphism $j$
of $\mathcal{M}$ such that $\mathcal{I}_{\mathrm{fix}}(j)=\mathcal{I}$,
where $\mathcal{I}_{\mathrm{fix}}(j)$ is defined in this context as the
submodel of $\mathcal{M}$ whose universe consists of element $m$ such that $%
j(x)=x$ for all $x<^{\mathcal{M}}m$. Furthermore, it is known
that:\smallskip

\noindent \textbf{(1)} \cite[Lemmas A.0 \& A.2]{ena06} If $\mathcal{M}$%
\ is a model of the fragment $\mathrm{I}\Delta _{0}$ of $\mathrm{PA}$, and $j
$ is a nontrivial automorphisms of $\mathcal{M}$, then $\mathcal{I}_{\mathrm{%
fix}}(j)$ satisfies $\mathrm{I}\Delta _{0}+\mathrm{Exp}+\mathrm{B}\Sigma _{1}
$ (where $\mathrm{Exp}$ expresses the totality of the exponential function,
and $\mathrm{B}\Sigma _{1}$ is the scheme of $\Sigma _{1}$-Collection).
\smallskip

\noindent \textbf{(2)} \cite[Theorem A]{ena06} Every countable model
of $\mathrm{I}\Delta _{0}+\mathrm{Exp}+\mathrm{B}\Sigma _{1}$ arises as $%
\mathcal{I}_{\mathrm{fix}}(j)$ for some nontrivial automorphism $j$ of a
model of $\mathrm{I}\Delta _{0}.$ \medskip

\noindent Our Theorem 3.4 is an analogue of (1); while Theorems 5.5, 5.6, and 5.15 provide analogues of (2). More generally, our results can be viewed as contributing to the project initiated in \cite{ena04} of  investigating the extent to which core results about automorphisms of model of arithmetic can be extended to the set-theoretic realm.

Another source of inspiration for our work is the metamathematics of $%
\mathrm{NFU}$, an urelement-variant of Quine's system \textquotedblleft New
Foundations\textquotedblright , $\mathrm{NF}$. Jensen's pioneering
work \cite{jen69} unveiled a magical link between models of $\mathrm{NFU}$ and automorphisms of models of $\mathrm{ZF}$-style set theories, a link
that has captured the imagination of other researchers, e.g., Holmes \cite{hol01}, Solovay \cite{solXX}, and two of the authors of the
present paper (\cite{ena04}, \cite{mck15}). Our results here have a number of implications for $\mathrm{NFU}$, for example, Theorem 5.15 can be used to show that every countable model of $\mathrm{ZFC}$ can be realized as the strongly cantorian part of a model of $\mathrm{NFU}$.  However, the precise implications of our results to the $\mathrm{NFU}$ setting is yet to be worked out and will be pursued elsewhere.\medskip

The plan of the paper is as follows. After dealing with preliminaries in
Section 2, we introduce the key notion of an \textquotedblleft $H$-cut\textquotedblright\ of a model of set theory in Section 3, where we
establish two important facts, namely: (1) every $H$-cut of a
model of $\mathrm{MOST}$ satisfies $\mathrm{MOST+\Delta }_{0}^{\mathcal{P}}$-
$\mathrm{Collection}$; and (2) $\mathcal{I}_{\mathrm{fix}}(j)$ is an $H$-cut
of $\mathcal{M}$, if $\mathcal{M}$ is a model of $\mathrm{MOST}$, and $j$ is
a nontrivial automorphisms of $\mathcal{M}$ such that $\mathcal{I}_{\mathrm{fix}}(j)$ includes $\omega^{\mathcal{M}}$ (equivalently: $j(m)=m$ whenever
whenever $m$\ is\ a finite\ ordinal in the sense of $\mathcal{M}$). The
central result of Section 4 is Theorem 4.1, whose iterated ultrapower proof
involves a rather intricate set-theoretical adaptation of a machinery that
was invented in the arithmetical context by Paris and Mills (\cite{PM79}, \cite[Theorem 3.5.5]{ks06}) and was further elaborated
in \cite{ena06}. One of the remarkable consequences of Theorem 4.1 is
that if $\mathcal{I}$\ is an $H$-cut of a countable model $\mathcal{M}$\ of $
\mathrm{MOST}+\Delta_{0}^{\mathcal{P}}$-$\mathrm{Collection}$, then $\mathcal{M}$ has a cofinal extension $\mathcal{N}$ that (1) $\mathcal{N}$
does not add any new members to elements of $I$; (2) $\mathcal{N}$ satisfies
$\mathrm{MOST}$; and (3) $\mathcal{N}$ carries an automorphism $j$ such that
$\mathcal{I}=\mathcal{I}_{\mathrm{fix}}(j).$ This consequence of Theorem 4.1
is put to work together with a key construction in Section 5 (Theorem 5.6)
to show that every countable recursively saturated model of $\mathrm{MOST+\Delta }_{0}^{\mathcal{P}}$-$\mathrm{Collection}$ can be realized as $\mathcal{I}_{\mathrm{fix}}(j)$; a result that, together with our work in
Section 3, yields the central theorem of our paper (Theorem 5.8) that
identifies the \textit{theory} of the class of models of the form $\mathcal{I}_{\mathrm{fix}}(j)$ to be precisely $\mathrm{MOST+\Delta }_{0}^{\mathcal{P}}
$-$\mathrm{Collection.}$ In Section 5 we also use our work in Section 4
together with some classical results of Friedman \cite{fri73} and
Hutchinson \cite{hut76} to identify two other sufficient conditions for a
countable model $\mathcal{N}$ to be realizable as $\mathcal{I}_{\mathrm{fix}}(j)$, namely: (1) $\mathcal{N}$ is a transitive model of $\mathrm{MOST}+\Delta_{0}^{\mathcal{P}}$-$\mathrm{Collection}$ and (2) $\mathcal{N}$
is a model of\textit{\ }$\mathrm{ZFC}$. Finally, in Section 6, we fine-tune
a theorem of Togha \cite{tog04} on automorphisms of countable recursively
saturated models of $\mathrm{ZFC}$ in a manner reminiscent of a refinement
of Smory\'{n}ski's aforementioned theorem established in \cite[Theorem B]{ena06}.

\section[Background and definitions]{Background and definitions} \label{Sec:Background}

Throughout this paper $\mathcal{L}$ will denote the language of set theory--- first-order logic endowed with a binary relation symbol $\in$ whose intended interpretation is membership. Structures will usually be denoted using upper-case calligraphy roman letters ($\mathcal{M}, \mathcal{N}, \ldots$) and the corresponding plain font letter ($M, N, \ldots$) will be used to denote the underlying set of that structure. If $\mathcal{M}$ is an $\mathcal{L}^\prime$-structure where $\mathcal{L}^\prime \supseteq \mathcal{L}$ and $a \in M$ then we will use $a^*$ to denote the class $\{x \in M \mid \mathcal{M} \models (x \in a)\}$. As usual $\Delta_0, \Sigma_1, \Pi_1, \ldots$ with be used to denote the L\'{e}vy classes of $\mathcal{L}$-formulae. We will also have cause to consider the class $\Delta_0^\mathcal{P}$ which is the smallest class of $\mathcal{L}$-formulae that contains all atomic formulae, contains all compound formulae formed using the connectives of first-order logic, and is closed under quantification in the form $\mathcal{Q} x \in y$ and $\mathcal{Q} x \subseteq y$ where $x$ and $y$ are distinct variables and $\mathcal{Q}$ is $\exists$ or $\forall$. If $\mathcal{L}^\prime \supseteq \mathcal{L}$ then we use $\Delta_0(\mathcal{L}^\prime)$ ($\Delta_0^\mathcal{P}(\mathcal{L}^\prime)$) to denote the smallest class of formulae that contains all atomic formulae, all compound formulae formed using the connectives of first-order logic, and is closed under quantification in the form $\mathcal{Q} x \in t$ (and $\mathcal{Q} x \subseteq t$) where $t$ is an $\mathcal{L}^\prime$-term and $x$ is a variable that does not appear in $t$, and $\mathcal{Q}$ is $\exists$ or $\forall$. The classes $\Sigma_1^\mathcal{P}, \Pi_1^\mathcal{P}, \ldots$ ($\Sigma_1(\mathcal{L}^\prime), \Pi_1(\mathcal{L}^\prime), \ldots$ and $\Sigma_1^\mathcal{P}(\mathcal{L}^\prime), \Pi_1^\mathcal{P}(\mathcal{L}^\prime), \ldots$) are defined inductively from the class $\Delta_0^\mathcal{P}$ (respectively $\Delta_0(\mathcal{L}^\prime)$ and $\Delta_0^\mathcal{P}(\mathcal{L}^\prime)$) in the same way as the classes $\Sigma_1, \Pi_1, \ldots$ are defined from $\Delta_0$.\\
\\
Let $\mathcal{L}^\prime \supseteq \mathcal{L}$ and let $\mathcal{M}$ and $\mathcal{N}$ be $\mathcal{L}^\prime$-structures. If $\mathcal{M}$ is a substructure of $\mathcal{N}$ then we will write $\mathcal{M} \subseteq \mathcal{N}$. If $\Gamma$ is a class of $\mathcal{L}^\prime$-formulae then we will write $\mathcal{M} \prec_\Gamma \mathcal{N}$ if $\mathcal{M} \subseteq \mathcal{N}$ and for every $\vec{a} \in M$, $\vec{a}$ satisfies the same $\Gamma$-formulae in both $\mathcal{M}$ and $\mathcal{N}$. If $\Gamma$ is $\mathcal{L}^\prime$ or  $\Sigma_n(\mathcal{L}^\prime)$ then we will abbreviate this notation by writing $\mathcal{M} \prec \mathcal{N}$ and $\mathcal{M} \prec_n \mathcal{N}$ respectively. If $\mathcal{M} \subseteq \mathcal{N}$ and for all $x \in M$ and $y \in N$,
$$\textrm{if } \mathcal{N} \models (y \in x) \textrm{ then } y \in M,$$
then we say that $\mathcal{N}$ is an \emph{end-extension} of $\mathcal{M}$ and write $\mathcal{M} \subseteq_e \mathcal{N}$. It is well-known that if $\mathcal{M} \subseteq_e \mathcal{N}$ then $\mathcal{M} \prec_{\Delta_0} \mathcal{N}$. If $\mathcal{N}$ is an end-extension of $\mathcal{M}$ and $\mathcal{M} \prec \mathcal{N}$ then we write $\mathcal{M} \prec_e \mathcal{N}$. In contrast, if $\mathcal{M} \subseteq \mathcal{N}$ and for all $x \in N$, there exists $y \in M$ such that $\mathcal{N} \models (x \in y)$, then we say that $\mathcal{N}$ is a \emph{cofinal extension} of $\mathcal{M}$ and we write $\mathcal{M} \subseteq_{\mathrm{cf}} \mathcal{N}$. And, if $\mathcal{M} \subseteq_{\mathrm{cf}} \mathcal{N}$ and $\mathcal{M} \prec \mathcal{N}$ then we write $\mathcal{M} \prec_{\mathrm{cf}} \mathcal{N}$.\\
\\          
\indent Let $\mathcal{L}^\prime \supseteq \mathcal{L}$ and let $\Gamma$ be a class of $\mathcal{L}^\prime$ formulae. We will use $\Gamma\textrm{-separation}$ and $\Gamma\textrm{-collection}$ to abbreviate the separation and collection schemes, respectively, restricted to formulae in $\Gamma$. We will also use the following axioms and schemes to axiomatize and study weak variants of $\mathrm{ZFC}$:
\begin{itemize}
\item[](transitive containment) $\forall x \exists y(\bigcup y \subseteq y \land x \subseteq y)$,
\item[]($\forall \kappa \exists \kappa^+$) for every initial ordinal $\kappa$, there exists a least initial ordinal $> \kappa$,
\item[]($\forall \kappa(2^\kappa \textrm{ exists})$) for every initial ordinal $\kappa$, there exists an initial $\lambda$ such that $\lambda= |\mathcal{P}(\kappa)|$,
\item[]($\Gamma$-foundation) for all $\Gamma$-formulae $\phi(x, \vec{z})$,
$$\forall \vec{z}(\exists x \phi(x, \vec{z}) \Rightarrow \exists y (\phi(y, \vec{z}) \land (\forall x \in y)\neg \phi(x, \vec{z})))$$
When $\Gamma$ only contains the formula $x \in z$ then we will refer to the single axiom in this scheme as \emph{set foundation}. 
\end{itemize}
The following subsystems of $\mathrm{ZFC}$ are studied in \cite{mat01}:

\begin{Definitions1}
$\mathrm{Mac}$ is the $\mathcal{L}$-theory with axioms: extensionality, pair, emptyset, union, infinity, powerset, transitive containment, $\Delta_0$-separation, set foundation, and the axiom of choice in the form: every set can be well-ordered.
\end{Definitions1}

\begin{Definitions1}
$\mathrm{KP}$ is the $\mathcal{L}$-theory with axioms: extensionality, pair, emptyset, union, $\Delta_0$-separation, $\Delta_0$-collection and $\Pi_1$-foundation.
\end{Definitions1}

\begin{Definitions1}
$\mathrm{KP}^\mathcal{P}$ is the $\mathcal{L}$-theory with axioms: extensionality, pair, emptyset, union, infinity, powerset, $\Delta_0$-separation, $\Delta_0^\mathcal{P}$-collection and $\Pi_1^\mathcal{P}$-foundation.
\end{Definitions1}

Note that $\mathrm{Mac}$ proves that for every set $x$, there exists a unique smallest transitive set, which we will denote $\mathrm{TC}(x)$, that contains $x$. We will use $\mathrm{Ord}$ and $\mathrm{Card}$ to denote the classes (definable in any extension of $\mathrm{Mac}$) of ordinals and cardinals respectively. We record the following well-known facts about fragments of the collection scheme:

\begin{Lemma1}
Let $n \in \omega$.
\begin{itemize}
\item[(i)] $\mathrm{Mac}+ \Pi_n\textrm{-collection} \vdash \Sigma_{n+1}\textrm{-collection}$,
\item[(ii)] $\mathrm{Mac}+ \Pi^\mathcal{P}_n\textrm{-collection} \vdash \Sigma^\mathcal{P}_{n+1}\textrm{-collection}$.
\end{itemize}
\Square
\end{Lemma1} 

\begin{Definitions1}
The theory $\mathrm{MOST}$ is obtained from $\mathrm{Mac}$ by adding $\Sigma_1$-separation and $\Delta_0$-collection.
\end{Definitions1}

Mathias~\cite{mat01} extensively studies the following axiom, originally proposed by Mitchell \cite{mit72}, which asserts that for every set $u$, there exists a universal transitive set which contains every transitive set that is of size $\leq |u|$:
\begin{itemize}
\item[](Axiom $\mathrm{H}$) $\forall u \exists T(\bigcup T \subseteq T \land \forall z(\bigcup z \subseteq z \land |z| \leq |u| \Rightarrow z \subseteq T))$. 
\end{itemize}
Mathias~\cite[Proposition 3.14]{mat01} shows that adding Axiom $\mathrm{H}$ to $\mathrm{Mac}$ yields $\mathrm{MOST}$.

\begin{Lemma1}
$\mathrm{Mac}+\textrm{Axiom }\mathrm{H}= \mathrm{MOST}$.\Square
\end{Lemma1} 

The set theory $\mathrm{MOST}$ is capable of defining the rank function, which we will denote $\rho$. It should be noted, however, that $\mathrm{MOST}$ does not prove enough recursion to ensure the totality of the function $\alpha \mapsto V_\alpha$. The following consequences of $\mathrm{MOST}$ are proved in \cite[Theorem scheme 3.13, Proposition 3.14, Theorem 3.18]{mat01}:

\begin{Lemma1}
$\mathrm{MOST}$ proves the following:
\begin{itemize}
\item[(i)] all instances of $\Delta_0^\mathcal{P}$-separation,
\item[(ii)] all instances of $\Pi_1$-separation,
\item[(iii)] all instances of $\Sigma_1$-collection, 
\item[(iv)] all instances of $\Pi_1$-foundation,
\item[(v)] every well-ordering is isomorphic to an ordinal,
\item[(vi)] $\forall \kappa \exists \kappa^+$,
\item[(vii)] for all cardinals $\kappa$, there is a set $H_\kappa$ of all sets whose transitive closure has cardinality less than $\kappa$.
\end{itemize}
\Square
\end{Lemma1}

\noindent (v) is a special case of what is known as Mostowski Isomorphism Theorem. \cite[Lemma 3.15]{mat01} also proves that a more general special case of this theorem which deals with well-founded extensional relations that are sets is provable in $\mathrm{MOST}$.

\begin{Lemma1} \label{Th:MostowskiIsomorphism}
$\mathrm{MOST}$ proves that if $R \subseteq X \times X$ is well-founded and extensional then the function $\varpi_R: X \longrightarrow V$ defined by
$$\varpi_R(x)= \{\varpi_R(y) \mid y \in X \land (\langle y, x \rangle \in R)\}$$
is defined on all of $X$. Moreover, $\varpi_R``X$ is the unique transitive set that is isomorphic to $R$.
\Square
\end{Lemma1}

In \cite{fk91} Thomas Forster and Richard Kaye introduce the notion of a powerset preserving end-extension. Here we give a slightly more general version of this notion that does not require the powerset axiom to hold in the structures that are being compared.

\begin{Definitions1}
Let $\mathcal{M}$ and $\mathcal{N}$ are $\mathcal{L}^\prime$-structures where $\mathcal{L}^\prime \supseteq \mathcal{L}$. We say that $\mathcal{N}$ is a powerset preserving end-extension of $\mathcal{M}$ and write $\mathcal{M} \subseteq_e^\mathcal{P} \mathcal{N}$ if
\begin{itemize}
\item[(i)] $\mathcal{M} \subseteq_e \mathcal{N}$
\item[(ii)] for all $x \in N$ and for all $y \in M$, if $\mathcal{N} \models (x \subseteq y)$ then $x \in M$.
\end{itemize}
\end{Definitions1}

Just as end-extensions preserve $\Delta_0$ properties, powerset preserving end-extensions preserve $\Delta_0^\mathcal{P}$ properties. The following is a slight modification of a result that appears in \cite{fk91}:

\begin{Lemma1}
Let $\mathcal{M}$ and $\mathcal{N}$ be $\mathcal{L}$-structures that satisfy extensionality. If $\mathcal{M} \subseteq_e^\mathcal{P} \mathcal{N}$ then $\mathcal{M} \prec_{\Delta_0^\mathcal{P}} \mathcal{N}$.
\end{Lemma1}

\begin{proof}
A straightforward induction on the structural complexity of a $\Delta_0^\mathcal{P}$-formula $\phi$.
\Square
\end{proof} 

The next definition captures an important notion that plays a key role in this study.

\begin{Definitions1} \label{Df:ToplessPowersetPreservingEndExtension}
Let $\mathcal{M}$ and $\mathcal{N}$ be $\mathcal{L}$-structures. We say that $\mathcal{N}$ is a topless end-extension of $\mathcal{M}$ and write $\mathcal{M} \subseteq_{\mathrm{topless}} \mathcal{N}$ if
\begin{itemize}
\item[(i)] $\mathcal{M} \subseteq_e \mathcal{N}$,
\item[(ii)] $M \neq N$,
\item[(iii)] if $C \in N$ and $C^* \subseteq \mathrm{Ord}^\mathcal{M}$ then $C \in M$.
\end{itemize}
If $\mathcal{M} \subseteq_{\mathrm{topless}} \mathcal{N}$ and $\mathcal{M} \subseteq_e^\mathcal{P} \mathcal{N}$ then we say that $\mathcal{N}$ is a topless powerset preserving end-extension of $\mathcal{M}$ and write $\mathcal{M} \subseteq_{\mathrm{topless}}^\mathcal{P} \mathcal{N}$. 
\end{Definitions1}

\noindent Note that if the $\mathcal{M}$ and $\mathcal{N}$ in Definition \ref{Df:ToplessPowersetPreservingEndExtension} satisfy $\mathrm{MOST}$ then condition (iii) is a paraphrasing of the assertion that there is no least ordinal in $\mathrm{Ord}^\mathcal{N} \backslash \mathrm{Ord}^\mathcal{M}$.\\ 
\\
\indent This paper studies the largest transitive initial segment of a model of $\mathrm{MOST}$ that is pointwise fixed by a non-trivial automorphism.

\begin{Definitions1}
Let $\mathcal{M}$ be an $\mathcal{L}$-structure with $\mathcal{M} \models \mathrm{MOST}$ and let $j: \mathcal{M} \longrightarrow \mathcal{M}$ be a non-trivial automorphism. We use $\mathrm{fix}(j)$ to denote the class of fixed points of $j$. I.e.
$$\mathrm{fix}(j)= \{x \in M \mid j(x)=x\}.$$
Define $\mathcal{I}_\mathrm{fix}(j)$ to be the substructure of $\mathcal{M}$ with domain
$$I_\mathrm{fix}(j)= \{x \in M \mid \forall y(\mathcal{M} \models (y \in \mathrm{TC}(\{x\})) \Rightarrow j(y)=y)\}.$$ 
\end{Definitions1}

Let $\mathcal{M}$ be an $\mathcal{L}^\prime$-structure. The structure $\mathcal{M}$ is said to be \emph{recursively saturated} if for all $a_1, \ldots, a_n \in M$ and for all recursive finitely realised types $\Gamma(x_1, \ldots, x_m, a_1, \ldots, a_n)$ in the language $\mathcal{L}^\prime$ with parameters $a_1, \ldots, a_n$, $\mathcal{M}$ realises $\Gamma(x_1, \ldots, x_m, a_1, \ldots, a_n)$. We refer the reader to \cite[\S 2.4]{CK90} for a detailed treatment of recursively saturated models. We will make use of the following nice feature of recursive saturation:

\begin{Theorems1} \label{Th:ExistenceOfRecursivelySaturatedModels}
(see \cite{CK90}) Let $\mathcal{L}^\prime$ be a recursive language. If $T$ is a consistent $\mathcal{L}^\prime$-theory with an infinite model then $T$ has a countable recursively saturated model. \Square
\end{Theorems1}

Fix a G\"{o}del coding of $\mathcal{L}$ in the theory $\mathrm{MOST}$ and use $\mathrm{Form}$ to denote the set of G\"{o}del codes of well-formed $\mathcal{L}$-formulae. Let $\mathcal{M}=\langle M, \in^\mathcal{M} \rangle$ be an $\mathcal{L}$-structure with $\mathcal{M} \models \mathrm{MOST}$. A satisfaction class for $\mathcal{M}$ is a class $S \subseteq M$ such that $S$ consists of ordered pairs $\langle a, b \rangle$ where $a \in (\mathrm{Form}^\mathcal{M})^*$ and $b \in M$, and for all $n \in \omega$,
\begin{itemize}
\item[]($S$ is $n$-correct) $S$ satisfies Tarski's inductive conditions for truth for all $\Sigma_n$-formulae.  
\end{itemize}
Let $\mathcal{L}_X$ be the extension of $\mathcal{L}$ obtained by adding a new unary predicate $X$. Use $\mathrm{ZFC}(X)$ to denote the $\mathcal{L}_X$-theory that extends $\mathrm{ZFC}$ with the schemes of $\mathcal{L}_X$-separation and $\mathcal{L}_X$-collection. The following result of Schlipf provides another characterisation of recursive saturation for countable models of $\mathrm{ZFC}$; the theorem below is an immediate consequence of putting Theorem 3.4 and the remark following it in \cite{sch78} with the well-known resplendence property of countable recursively saturated models (see \cite[Theorem 1.3]{sch78}). 

\begin{Theorems1} \label{Th:SatisfactionClassesForRecSatModels}
Let $\mathcal{M}= \langle M, \in^{\mathcal{M}} \rangle$ be a countable $\omega$-nonstandard model of $\mathrm{ZFC}$. The structure $\mathcal{M}$ is recursively saturated if and only if $\mathcal{M}$ can be expanded to an $\mathcal{L}_X$-structure $\mathcal{M}_\mathrm{Sat}= \langle M, \in^{\mathcal{M}}, X^{\mathcal{M}} \rangle$ such that
\begin{itemize}
\item[(I)] $\mathcal{M}_\mathrm{Sat} \models \mathrm{ZFC}(X)$,
\item[(II)] $X^{\mathcal{M}}$ is a satisfaction class for $\mathcal{M}$.
\end{itemize}
\Square 
\end{Theorems1}

\noindent Note that we cannot expect that a satisfaction class obtained from Theorem \ref{Th:SatisfactionClassesForRecSatModels} will be $n$-correct for all internal natural numbers of $\mathcal{M}$, as this would prove the consistency of $\mathrm{ZFC}$. However, if $\mathcal{M}= \langle M, \in^\mathcal{M}, X^\mathcal{M} \rangle$ is an $\mathcal{L}_X$-structure that satisfies (I) and (II) of Theorem \ref{Th:SatisfactionClassesForRecSatModels} and is $\omega$-nonstandard then there is a nonstandard $s \in (\omega^\mathcal{M})^*$ such that $X^\mathcal{M}$ is $s$-correct. This follows from overspill and the fact that there is an $\mathcal{L}_X$-formula with parameter $n$ that expresses that $X$ is $n$-correct.

\section[The structure and first-order theory of $\mathcal{I}_\mathrm{fix}(j)$]{The structure and first-order theory of $\mathcal{I}_\mathrm{fix}(j)$} \label{Sec:StructureOfIFix}

In this section we investigate the first-order theory of $\mathcal{I}_\mathrm{fix}(j)$ and the properties of this structure in relation to the domain of the automorphism $j$. We show that if $j: \mathcal{M} \longrightarrow \mathcal{M}$ is a non-trivial automorphism that hereditarily fixes $\omega^\mathcal{M}$ and $\mathcal{M}$ satisfies $\mathrm{MOST}$ then $\mathcal{I}_\mathrm{fix}(j)$ satisfies $\mathrm{MOST}$ and $\mathcal{M}$ is a topless powerset preserving end-extension of $\mathcal{I}_\mathrm{fix}(j)$. The fact that $\mathcal{M}$ is a topless powerset preserving end-extension of $\mathcal{I}_\mathrm{fix}(j)$ implies that $\mathcal{I}_\mathrm{fix}(j)$ also satisfies all instances of $\Delta_0^\mathcal{P}$-collection. We capture the relationship that we will prove holds between $\mathcal{I}_\mathrm{fix}(j)$ and $\mathcal{M}$ in the following definition:
   
\begin{Definitions1}
Let $\mathcal{M}$ be an $\mathcal{L}$-structure with $\mathcal{M} \models \mathrm{MOST}$. We say that $I \subseteq M$ is an $H$-cut of $\mathcal{M}$ if
\begin{itemize}
\item[(i)] $\langle I, \in^\mathcal{M} \rangle \subseteq_{\mathrm{topless}}^\mathcal{P} \mathcal{M}$,
\item[(ii)] $\langle I, \in^\mathcal{M} \rangle \models \mathrm{MOST}$.
\end{itemize} 
\end{Definitions1}

\noindent Note that if $I$ is an $H$-cut of a model $\mathcal{M}$ of $\mathrm{MOST}$ then $I$ is a union of sets $H_\kappa$ in $\mathcal{M}$. We make this explicit with the following observation:

\begin{Lemma1} \label{Th:UnionOfHs}
Let $\mathcal{M}$ be an $\mathcal{L}$-structure with $\mathcal{M} \models \mathrm{MOST}$. If $I \subseteq M$ is an $H$-cut of $\mathcal{M}$ then for all $\kappa \in I$,
$$\langle I, \in^\mathcal{M} \rangle \models (\kappa \textrm{ is a cardinal}) \textrm{ if and only if } \mathcal{M} \models (\kappa \textrm{ is a cardinal})$$
$$\textrm{and if } \mathcal{M} \models (\kappa \textrm{ is a cardinal}) \textrm{ then } H_\kappa^\mathcal{M}=H_\kappa^{\langle I, \in^\mathcal{M}\rangle}.$$
\end{Lemma1}

\begin{proof}
Let $I \subseteq M$ be an $H$-cut of $\mathcal{M}$ and let $\kappa \in I$. The expression that `$\kappa$ is a cardinal' can be written as a $\Delta_0^\mathcal{P}$-formula which takes the parameter $\kappa \times \kappa$. This shows that
$$\langle I, \in^\mathcal{M} \rangle \models (\kappa \textrm{ is a cardinal}) \textrm{ if and only if } \mathcal{M} \models (\kappa \textrm{ is a cardinal}).$$
The formula $y= \mathrm{TC}(x)$ can be expressed as a $\Delta_0^\mathcal{P}$-formula. The fact that $(H_\kappa^\mathcal{\langle I, \in^\mathcal{M}\rangle})^* \subseteq (H_\kappa^\mathcal{M})^*$ now follows from the fact that $x \in H_\kappa$ can be expressed as a $\Delta_0^\mathcal{P}$-formula with parameters $x \times \kappa$ and $\mathrm{TC}(x)$. To get the reverse inclusion, let $x \in (H_\kappa^\mathcal{M})^*$. All we need to show is that $x \in I$. Work inside $\mathcal{M}$. Since $x \in H_\kappa$, there is an $A \subseteq \kappa$ and a well-founded and extensional $R \subseteq A \times A$ such that $\varpi_R``A= \mathrm{TC}(\{x\})$. Since $\langle I, \in^\mathcal{M} \rangle \subseteq_{e}^\mathcal{P} \mathcal{M}$, it follows that $A, R \in I$ and that $\langle I, \in^\mathcal{M} \rangle$ thinks $R$ is well-founded and extensional. Therefore 
$$\langle I, \in^\mathcal{M} \rangle \models (\varpi_R``A \textrm{ is a transitive set}).$$
Set foundation can then be used to show that $(\varpi_R``A)^\mathcal{M}= (\varpi_R``A)^{\langle I, \in^\mathcal{M}\rangle}$. This shows that $x \in I$ and concludes the proof that $H_\kappa^\mathcal{M}=H_\kappa^{\langle I, \in^\mathcal{M}\rangle}$.         
\Square
\end{proof}

Togha~\cite[Definition 1]{tog04} introduces the notion of a cut of the ordinals of a model of set theory.

\begin{Definitions1} \label{Df:OrdinalCutClosedUnderExponentiation}
Let $\mathcal{M}$ be an $\mathcal{L}$-structure with $\mathcal{M} \models \mathrm{MOST}$. We say that $K \subseteq \mathrm{Ord}^\mathcal{M}$ is an ordinal cut of $\mathcal{M}$ closed under exponentiation if $\langle K, \in^\mathcal{M} \rangle \subseteq_{topless} \langle \mathrm{Ord}^\mathcal{M}, \in^\mathcal{M} \rangle$ and for all $\mathcal{M}$-cardinals $\kappa \in K$, $(2^\kappa)^\mathcal{M} \in K$.
\end{Definitions1}

\noindent If $\mathcal{M} \models \mathrm{MOST}$ and $I \subseteq M$ is an $H$-cut of $\mathcal{M}$ then $\mathrm{Ord}^{\langle I, \in^\mathcal{M} \rangle}$ is an ordinal cut of $\mathcal{M}$ that is closed under exponentiation and, moreover, $\mathrm{Ord}^{\langle I, \in^\mathcal{M} \rangle}$ completely determines $I$.

\begin{Lemma1} \label{Th:HCutDeterminedByOrdinalCut}
Let $\mathcal{M}$ be an $\mathcal{L}$-structure with $\mathcal{M} \models \mathrm{MOST}$. If $I \subseteq M$ is an $H$-cut of $\mathcal{M}$ then
\begin{itemize}
\item[(I)] $\mathrm{Ord}^{\langle I, \in^\mathcal{M} \rangle}$ is an ordinal cut of $\mathcal{M}$ closed under exponentiation,
\item[(II)] $I= \bigcup \{ (H_\kappa^\mathcal{M})^* \mid (\mathcal{M} \models (\kappa \textrm{ is a cardinal})) \land (\kappa \in I)\}$.
\end{itemize}
\end{Lemma1}

\begin{proof}
Let $I \subseteq M$ be an $H$-cut of $\mathcal{M}$. It is clear that $\langle \mathrm{Ord}^{\langle I, \in^\mathcal{M} \rangle}, \in^\mathcal{M} \rangle \subseteq_{topless} \langle \mathrm{Ord}^\mathcal{M}, \in^\mathcal{M}\rangle$. Let $\kappa \in I$ be an $\mathcal{M}$-cardinal. It follows from Lemma \ref{Th:UnionOfHs} that $\kappa$ is an $\langle I, \in^\mathcal{M} \rangle$-cardinal and $(2^\kappa)^\mathcal{M}= (2^\kappa)^{\langle I, \in^\mathcal{M} \rangle}$. Therefore $(2^\kappa)^\mathcal{M} \in I$. This shows that $\mathrm{Ord}^{\langle I, \in^\mathcal{M} \rangle}$ is an ordinal cut of $\mathcal{M}$ closed under exponentiation. Let $J= \bigcup \{ (H_\kappa^\mathcal{M})^* \mid (\mathcal{M} \models (\kappa \textrm{ is a cardinal})) \land (\kappa \in I)\}$. Lemma \ref{Th:UnionOfHs} shows that $J \subseteq I$. Conversely, let $x \in I$. Let $\kappa= |\mathcal{P}^\mathcal{\langle I, \in^\mathcal{M} \rangle}(\mathrm{TC}(x))| \in I$. Therefore, by Lemma \ref{Th:UnionOfHs}, $\kappa$ is an $\mathcal{M}$-cardinal and $x \in (H_\kappa^{\langle I, \in^\mathcal{M} \rangle})^*=(H_\kappa^{\mathcal{M}})^*$. This shows that $I=J$ and completes the proof of the lemma.
\Square  
\end{proof}  

The fact that an $H$-cut sits toplessly inside a powerset preserving end-extension that satisfies $\mathrm{MOST}$ implies that it satisfies all instances of $\Delta_0^\mathcal{P}$-collection.

\begin{Lemma1} \label{Th:CollectionInIfix}
Let $\mathcal{M}$ be an $\mathcal{L}$-structure with $\mathcal{M} \models \mathrm{MOST}$. If $I \subseteq M$ is an $H$-cut of $\mathcal{M}$ then
$$\langle I, \in^\mathcal{M} \rangle \models \Delta_0^\mathcal{P}\textrm{-collection}.$$
\end{Lemma1}

\begin{proof}
Let $\phi(x, y, \vec{z})$ be a $\Delta_0^\mathcal{P}$-formula and let $\vec{a}, b \in I$ such that
$$\langle I, \in^\mathcal{M} \rangle \models (\forall x \in b) \exists y \phi(x, y, \vec{a}).$$
Let $\kappa \in \mathrm{Ord}^\mathcal{M} \backslash \mathrm{Ord}^{\langle I, \in^\mathcal{M} \rangle}$ be a cardinal. Note that $I \subseteq (H_{\kappa^+}^\mathcal{M})^*$. Define $\theta(u, w, w^\prime, \vec{v})$ to be the formula
$$(u \textrm{ is an ordinal})\land$$ 
$$(\exists f, T, y \in w)(\exists x \in w^\prime)\left(\begin{array}{c}
(T= \mathrm{TC}(y)) \land (f: u \longrightarrow T \textrm{ is a bijection}) \land \phi(x, y, \vec{v}) \land\\
(\forall p, g, S, q \in w) \left(\begin{array}{c}
(p \textrm{ is an ordinal}) \land (S= \mathrm{TC}(q)) \land\\ 
(g: p \longrightarrow S \textrm{ is a bijection}) \land \phi(x, q, \vec{v})\\
\Rightarrow u \leq p 
\end{array}\right) 
\end{array}\right).$$
So $\theta(u, w, w^\prime, \vec{v})$ is a $\Delta_0^\mathcal{P}$ formula. Working inside $\mathcal{M}$, $\Delta_0^\mathcal{P}$-separation ensures that
$$C= \{\lambda \in H_{\kappa^+} \mid \theta(\lambda, H_{\kappa^+}, b, \vec{a})\}$$
is a set. Lemma \ref{Th:UnionOfHs} implies that $C^* \subseteq I$. Therefore, since $\langle I, \in^\mathcal{M} \rangle \subseteq_{\mathrm{topless}}^\mathcal{P} \mathcal{M}$, $C \in I$. Working inside $\langle I, \in^\mathcal{M} \rangle$, let $\mu= \sup C$. Now, $\mu$ is a cardinal and 
$$\langle I, \in^\mathcal{M} \rangle \models (\forall x \in b)(\exists y \in H_{\mu^+})\phi(x, y, \vec{a}).$$    
\Square
\end{proof}

The remainder of this section is dedicated to proving that if $j: \mathcal{M} \longrightarrow \mathcal{M}$ is a non-trivial automorphism of a model of $\mathrm{MOST}$ that hereditarily fixes $\omega^\mathcal{M}$ then $I_\mathrm{fix}(j)$ is an $H$-cut of $\mathcal{M}$.

\begin{Theorems1} \label{Th:IfixIsAHcut}
Let $\mathcal{M}$ be an $\mathcal{L}$-structure with $\mathcal{M} \models \mathrm{MOST}$. If $j: \mathcal{M} \longrightarrow \mathcal{M}$ is a non-trivial automorphism with $\omega^\mathcal{M} \in I_\mathrm{fix}(j)$ then $I_\mathrm{fix}(j) \subseteq M$ is an $H$-cut of $\mathcal{M}$. 
\end{Theorems1}

Combined with Lemma \ref{Th:CollectionInIfix} this shows that if $j: \mathcal{M} \longrightarrow \mathcal{M}$ is an automorphism that satisfies the conditions of Theorem \ref{Th:IfixIsAHcut} then $\mathcal{I}_\mathrm{fix}(j)$ satisfies $\mathrm{MOST}$ plus the scheme of $\Delta_0^\mathcal{P}$-collection.

\begin{Coroll1}
Let $\mathcal{M}$ be an $\mathcal{L}$-structure with $\mathcal{M} \models \mathrm{MOST}$. If $j: \mathcal{M} \longrightarrow \mathcal{M}$ is a non-trivial automorphism with $\omega^\mathcal{M} \in I_\mathrm{fix}(j)$ then
$$\mathcal{I}_\mathrm{fix}(j) \models \mathrm{MOST}+\Delta_0^\mathcal{P}\textrm{-collection}.$$
\Square
\end{Coroll1}  

For the remainder of this section we let $j: \mathcal{M} \longrightarrow \mathcal{M}$ be a non-trivial automorphism where $\mathcal{M} \models \mathrm{MOST}$ and $j(x)= x$ for all $x \in (\omega^\mathcal{M})^*$. We first show that $\mathcal{M}$ is a topless powerset preserving end-extension of $\mathcal{I}_\mathrm{fix}(j)$.

\begin{Lemma1} \label{Th:IfixPowersetPreservingSubmodel}
$\mathcal{I}_\mathrm{fix}(j)$ satisfies extensionality and powerset, and $\mathcal{I}_\mathrm{fix}(j) \subseteq_{\mathrm{topless}}^\mathcal{P} \mathcal{M}$.
\end{Lemma1}

\begin{proof}
It is clear from the definition of $\mathcal{I}_\mathrm{fix}(j)$ that $I_\mathrm{fix}(j)$ is a transitive proper subclass of $\mathcal{M}$. It follows that $\mathcal{I}_\mathrm{fix}(j) \models (\textrm{extensionality})$ and $I_\mathrm{fix}(j) \neq M$. Let $y \in I_\mathrm{fix}(j)$. Let $x \in M$ be such that $\mathcal{M} \models (x \subseteq y)$. Since $y \in I_\mathrm{fix}(j)$, $j$ fixes every element of $x^*$. Since $j$ is automorphism it follows that $j(x)= x$ and so $x \in I_\mathrm{fix}(j)$. This shows that for all $x \in \mathcal{P}^\mathcal{M}(y)^*$, $x \in I_\mathrm{fix}(j)$. And so $\mathcal{P}^\mathcal{M}(y) \in I_\mathrm{fix}(j)$. Since $I_\mathrm{fix}(j)$ is transitive, this shows that $\mathcal{I}_\mathrm{fix}(j) \models (\textrm{powerset})$ and $\mathcal{I}_\mathrm{fix}(j) \subseteq_{e}^\mathcal{P} \mathcal{M}$. To show that $\mathcal{I}_\mathrm{fix}(j)$ sits toplessly in $\mathcal{M}$, let $C \in M$ with $C^* \subseteq \mathrm{Ord}^{\mathcal{I}_\mathrm{fix}(j)}$. It immediately follows from the fact that $j$ is an automorphism that $j(C)=C$. Therefore $C \in I_\mathrm{fix}(j)$. This completes the proof that $\mathcal{I}_\mathrm{fix}(j) \subseteq_{\mathrm{topless}}^\mathcal{P} \mathcal{M}$. 
\Square
\end{proof}

\noindent To complete the proof that $\mathcal{I}_\mathrm{fix}(j)$ is an $H$-cut of $\mathcal{M}$ we need to show that $\mathcal{I}_\mathrm{fix}(j)$ satisfies $\mathrm{MOST}$.

\begin{Lemma1}
$\mathcal{I}_\mathrm{fix}(j) \models \mathrm{Mac}$.
\end{Lemma1}

\begin{proof}
Lemma \ref{Th:IfixPowersetPreservingSubmodel} and the fact that $I_\mathrm{fix}(j)$ is a transitive subclass of $\mathcal{M}$ implies that extensionality, union, pair, emptyset, powerset, $\Delta_0$-separation, set foundation and the axiom of choice hold in $\mathcal{I}_\mathrm{fix}(j)$. Since $\omega^\mathcal{M} \in I_\mathrm{fix}(j)$, $\mathcal{I}_\mathrm{fix}(j)$ satisfies the axiom of infinity. To see that transitive containment holds in $\mathcal{I}_\mathrm{fix}(j)$, let $x \in I_\mathrm{fix}(j)$. Since $j$ fixes every element of $\mathrm{TC}^\mathcal{M}(x)^*$, $j$ also fixes $\mathrm{TC}^\mathcal{M}(x)$. Therefore $\mathrm{TC}^\mathcal{M}(x) \in I_\mathrm{fix}(j)$ and $\mathcal{I}_\mathrm{fix}(j)$ believes that $\mathrm{TC}^\mathcal{M}(x)$ is a transitive set which contains $x$.         
\Square
\end{proof}

\noindent It remains to show that $\mathcal{I}_\mathrm{fix}(j)$ satisfies axiom $\mathrm{H}$.

\begin{Lemma1} \label{Th:MostowskiCollapseInIfix}
For all $R, X \in I_\mathrm{fix}(j)$,
$$\mathcal{M} \models (R \subseteq X \times X \textrm{ is well-founded and extensional})$$ 
$$\textrm{if and only if } \mathcal{I}_\mathrm{fix}(j) \models (R \subseteq X \times X \textrm{ is well-founded and extensional})$$
and if 
$$\mathcal{M} \models (R \subseteq X \times X \textrm{ is well-founded and extensional})$$
then $(\varpi_R``X)^\mathcal{M} \in I_\mathrm{fix}(j)$.
\end{Lemma1}

\begin{proof}
Let $R, X \in I_\mathrm{fix}(j)$. The fact that
$$\mathcal{M} \models (R \subseteq X \times X \textrm{ is well-founded and extensional})$$ 
$$\textrm{if and only if } \mathcal{I}_\mathrm{fix}(j) \models (R \subseteq X \times X \textrm{ is well-founded and extensional}),$$
follows from the fact that `$R \subseteq X \times X$ is well-founded and extensional' can be expressed as a $\Delta_0^\mathcal{P}$-formula. Assume that
$$\mathcal{M} \models (R \subseteq X \times X \textrm{ is well-founded and extensional}).$$
The fact that $(\varpi_R``X)^\mathcal{M} \in I_\mathrm{fix}(j)$ follows from the fact that $j$ is an automorphism and $\varpi_R$ is definable.    
\Square
\end{proof}

\begin{Lemma1} \label{Th:CardinalsArePreserved}
For all $\kappa \in I_\mathrm{fix}(j)$,
$$\mathcal{M} \models (\kappa \textrm{ is a cardinal}) \textrm{ if and only if } \mathcal{I}_\mathrm{fix}(j) \models (\kappa \textrm{ is a cardinal}).$$
\end{Lemma1}

\begin{proof}
As in the proof of Lemma \ref{Th:UnionOfHs}, this follows from the fact that `$\kappa$ is a cardinal' can be expressed as a $\Delta_0^\mathcal{P}$-formula with parameters from $I_\mathrm{fix}(j)$. 
\Square
\end{proof}

\begin{Lemma1} \label{Th:EveryWOIsomorphicToAnOrdinal}
$\mathcal{I}_\mathrm{fix}(j) \models (\textrm{every well-ordering is isomorphic to an ordinal})$.
\end{Lemma1}

\begin{proof}
Let $R, X \in I_\mathrm{fix}(j)$ be such that $\mathcal{I}_\mathrm{fix}(j) \models (R \textrm{\ is a well-ordering of } X)$. It follows from Lemma \ref{Th:MostowskiCollapseInIfix} (and the fact that being a linear ordering is absolute) that $\mathcal{M} \models (R \textrm{\ is a well-ordering of } X)$. Now, $(\varpi_R``X)^\mathcal{M}$ is an ordinal and, by Lemma \ref{Th:MostowskiCollapseInIfix}, $(\varpi_R``X)^\mathcal{M} \in I_\mathrm{fix}(j)$. Let $\alpha= (\varpi_R``X)^\mathcal{M}$. In $\mathcal{M}$ there is a bijection $f \subseteq X \times \alpha$ witnessing that $R$ is isomorphic to $\alpha$. Therefore $f \in I_\mathrm{fix}(j)$ and $f$ witnesses that $R$ is isomorphic to $\alpha$ in $\mathcal{I}_\mathrm{fix}(j)$.     
\Square
\end{proof}

\begin{Lemma1} \label{Th:CardinalExponentsExist}
$\mathcal{I}_\mathrm{fix}(j) \models (\textrm{for all cardinals } \kappa, \textrm{ the cardinal } 2^\kappa \textrm{ exists})$.
\end{Lemma1}

\begin{proof}
Let $\kappa \in I_\mathrm{fix}(j)$ be a cardinal. Work inside $\mathcal{I}_\mathrm{fix}(j)$. Let $R \subseteq \mathcal{P}(\kappa) \times \mathcal{P}(\kappa)$ be a well-ordering. By Lemma \ref{Th:EveryWOIsomorphicToAnOrdinal}, $R$ is isomorphic to an ordinal $\alpha \geq 2^\kappa$. 
\Square
\end{proof}

\begin{Lemma1} \label{Th:HFromMinIfix}
If $\kappa \in I_\mathrm{fix}(j)$ is a cardinal then $H_\kappa^\mathcal{M} \in I_\mathrm{fix}(j)$.  
\end{Lemma1}

\begin{proof}
Let $\kappa \in I_\mathrm{fix}(j)$ be a cardinal. Work inside $\mathcal{M}$. Note that $|H_\kappa|= 2^\kappa$ and $2^\kappa \in I_\mathrm{fix}(j)$. Using a bijection $f: H_\kappa \longrightarrow 2^\kappa$ one can find an $R \subseteq 2^\kappa \times 2^\kappa$ such that $\varpi_R``2^\kappa= H_\kappa$. Therefore, it follows from Lemma \ref{Th:MostowskiCollapseInIfix} that $H_\kappa^\mathcal{M} \in I_\mathrm{fix}(j)$.
\Square   
\end{proof}

\begin{Lemma1}
$\mathcal{I}_\mathrm{fix}(j) \models \textrm{Axiom }\mathrm{H}$.
\end{Lemma1}

\begin{proof}
By Lemma \ref{Th:EveryWOIsomorphicToAnOrdinal} it is enough to show that for every cardinal $\kappa$, there exists a transitive set which contains all transitive sets with cardinality $\leq \kappa$. Let $\kappa \in I_\mathrm{fix}(j)$ be a cardinal. Lemmas \ref{Th:CardinalExponentsExist} and \ref{Th:HFromMinIfix} show that $2^\kappa \in I_\mathrm{fix}(j)$ and $H_{2^\kappa}^\mathcal{M} \in I_\mathrm{fix}(j)$. The point $H_{2^\kappa}^\mathcal{M} \in I_\mathrm{fix}(j)$ is transitive and we claim that every transitive set of size $\leq \kappa$ in $I_\mathrm{fix}(j)$ is contained in $H_{2^\kappa}^\mathcal{M}$. To see this, let $z \in I_\mathrm{fix}(j)$ be such that $\mathcal{I}_\mathrm{fix}(j) \models (\bigcup z \subseteq z) \land (|z| \leq \kappa)$. Therefore $\mathcal{M} \models (\bigcup z \subseteq z) \land (|z| \leq \kappa)$ and so $\mathcal{I}_\mathrm{fix}(j) \models (z \subseteq H_{2^\kappa}^\mathcal{M})$. This shows that $\mathcal{I}_\mathrm{fix}(j)$ satisfies axiom $\mathrm{H}$.       
\Square
\end{proof}

\noindent This completes the proof of Theorem \ref{Th:IfixIsAHcut}.

\section[Realizing $H$-cuts as $\mathcal{I}_\mathrm{fix}(j)$]{Realizing $H$-cuts as $\mathcal{I}_\mathrm{fix}(j)$} \label{Sec:RealisingHcutsAsIFix}

In this section we will show if $I$ is an $H$-cut of a countable model $\mathcal{M}$ of $\mathrm{MOST}+\Delta_0^\mathcal{P}\textrm{-collection}$ then there is $\mathcal{M} \subseteq_{\mathrm{cf}} \mathcal{N}$ with $\mathcal{N} \models \mathrm{MOST}$ and a non-trivial automorphism $j: \mathcal{N} \longrightarrow \mathcal{N}$ such that $\mathcal{I}_\mathrm{fix}(j)= \langle I, \in^\mathcal{M} \rangle$ and $\mathrm{fix}(j) = M$. More generally, we will prove the following result:

\begin{Theorems1} \label{Th:MainRealisingHCutTheorem}
Let $\mathcal{M}=\langle M, \in^\mathcal{M} \rangle$ be countable such that $\mathcal{M} \models \mathrm{MOST}+\Delta_0^\mathcal{P}\textrm{-collection}$ and let $I \subseteq M$ be an $H$-cut of $\mathcal{M}$. If $\mathbb{L}$ is a linear order then there is $\mathcal{M} \subseteq_{\mathrm{cf}} \mathcal{N}$ with $\mathcal{N} \models \mathrm{MOST}$ and an embedding $j \mapsto \check{j}$ of $\mathrm{Aut}(\mathbb{L})$ into $\mathrm{Aut}(\mathcal{N})$ such that if $j$ has no fixed points then $\mathcal{I}_\mathrm{fix}(\check{j}) = \langle I, \in^\mathcal{M} \rangle$ and $\mathrm{fix}(\check{j}) = M$. 
\end{Theorems1}

\noindent By setting $\mathbb{L}$ in Theorem \ref{Th:MainRealisingHCutTheorem} to be the linear ordering $\mathbb{Z}$ and by letting $j: \mathbb{Z} \longrightarrow \mathbb{Z}$ be the automorphism that sends $i \mapsto i+1$ for all $i \in \mathbb{Z}$, we obtain:

\begin{Coroll1} \label{Th:EveryHCutCanBeRealisedAsIfix}
Let $\mathcal{M}=\langle M, \in^\mathcal{M} \rangle$ be countable such that $\mathcal{M} \models \mathrm{MOST}+\Delta_0^\mathcal{P}\textrm{-collection}$ and let $I \subseteq M$ be an $H$-cut of $\mathcal{M}$. Then there is $\mathcal{M} \subseteq_{\mathrm{cf}} \mathcal{N}$ with $\mathcal{N} \models \mathrm{MOST}$ and an automorphism $j: \mathcal{N} \longrightarrow \mathcal{N}$ such that $\mathcal{I}_\mathrm{fix}(j)= \langle I, \in^\mathcal{M} \rangle$ and $\mathrm{fix}(j) = M$. \Square
\end{Coroll1}

\noindent Theorem \ref{Th:MainRealisingHCutTheorem} will be proved by adapting a construction, originally due to Jeff Paris and George Mills \cite{PM79} and applied to non-standard models of $\mathrm{PA}$, to non-standard models of $\mathrm{MOST}+\Delta_0^\mathcal{P}\textrm{-collection}$. This construction was also used by the first author in \cite{ena06} to prove that every cut of a countable model of $I\Delta_0+B\Sigma_1$ that is closed under exponentiation can be realised as the largest initial segment of a model of $I\Delta_0$ that is pointwise fixed by some non-trivial automorphism. Theorem \ref{Th:MainRealisingHCutTheorem} and Corollary \ref{Th:EveryHCutCanBeRealisedAsIfix} are the set-theoretic analogues of results proved in Section 3 of \cite{ena06}.\\
\\
\indent Throughout this section we fix an $\mathcal{L}$-structure $\mathcal{M} \models \mathrm{MOST}+\Delta_0^\mathcal{P}\textrm{-collection}$ with $|M|= \aleph_0$, and $I \subseteq M$ an $H$-cut of $\mathcal{M}$. We also fix a linear order $\mathbb{L}$. For simplicity we will use $<$ for the order relation on $\mathbb{L}$. We will write $\mathcal{I}$ for the substructure of $\mathcal{M}$ with domain $I$. We also fix $\bar{\kappa} \in \mathrm{Ord}^\mathcal{M} \backslash \mathrm{Ord}^\mathcal{I}$ such that $\mathcal{M} \models (\bar{\kappa} \textrm{ is a regular cardinal})$. We begin by noting that there is no least cardinal in $\mathrm{Ord}^\mathcal{M} \backslash \mathrm{Ord}^\mathcal{I}$, and so, regardless of our choice of $\bar{\kappa}$, there are infinitely many cardinals in $\mathrm{Ord}^\mathcal{M} \backslash \mathrm{Ord}^\mathcal{I}$ below $\bar{\kappa}$.

\begin{Lemma1}
There is no least cardinal in the class $\mathrm{Ord}^\mathcal{M} \backslash \mathrm{Ord}^\mathcal{I}$.
\end{Lemma1}

\begin{proof}
Suppose that $\kappa \in \mathrm{Ord}^\mathcal{M} \backslash \mathrm{Ord}^\mathcal{I}$ is least such that $\mathcal{M} \models (\kappa \textrm{ is a cardinal})$. Working inside $\mathcal{M}$, consider the set
$$C=\{\alpha \in \kappa \mid |\alpha| < \kappa\}.$$
Note that $C \in M$. Since there is no largest cardinal in $\mathcal{I}$, $C^*= \mathrm{Ord}^\mathcal{I} \subseteq I$. This shows that $C \notin I$, which contradicts the fact that $I$ is an $H$-cut of $\mathcal{M}$.
\Square
\end{proof} 

Using an iterated ultrapower construction we will obtain a cofinal extension $\mathcal{N}$ of $\mathcal{M}$ that satisfies $\mathrm{MOST}$ and such that any fixed point free automorphism of $\mathbb{L}$ generates an automorphism $j: \mathcal{N} \longrightarrow \mathcal{N}$ with $\mathcal{I}_\mathrm{fix}(j)= \mathcal{I}$ and $\mathrm{fix}(j) = M$. The ultrafilter $\mathcal{U}$ used in this iterated ultrapower construction will be an ultrafilter on the subsets of $\bar{\kappa}$ in $\mathcal{M}$. The elements of $\mathcal{U}$ will be large in the following sense:

\begin{Definitions1}
We say that $X \in (\mathcal{P}(\bar{\kappa})^\mathcal{M})^*$ is $I$-large if $|X|^\mathcal{M} \notin I$.
\end{Definitions1}

If $X$ is a set and $\lambda$ is a cardinal then we use $[X]^\lambda$ to denote the set of all subsets of $X$ of size $\lambda$. Recall the canonical partion relation $\kappa \rightarrow *(\lambda)^n$, first considered in \cite{ER50}, that generalises the classical partition relation $\kappa \rightarrow (\lambda)^n_\mu$:

\begin{Definitions1}
Let $\kappa$ and $\lambda$ be cardinals and let $n \in \omega$. If $f$ is a function with domain $[\kappa]^n$, $H \subseteq \kappa$, and $\Delta \subseteq n$ such that for all $s_0 < \cdots < s_{n-1}$ and $t_0 < \cdots < t_{n-1}$ in $H$ we have
$$f(\{s_0, \ldots, s_{n-1}\})= f(\{t_0, \ldots, t_{n-1}\}) \textrm{ if and only if } (\forall i \in \Delta)(s_i=t_i),$$
then we say that $H$ is $f$-canonical. We write $\kappa \rightarrow *(\lambda)^n$ if for all functions $f$ with domain $[\kappa]^n$, there is $H \subseteq \kappa$ with $|H|= \lambda$ and $H$ is $f$-canonical. 
\end{Definitions1}

\noindent This notion allows us to make explicit the combinatorial properties we will require of the ultrafilter used to extend $\mathcal{M}$.

\begin{Definitions1} \label{Df:IComplete}
We say that a non-principal (n.p.) ultrafilter $\mathcal{U} \subseteq (\mathcal{P}(\bar{\kappa})^\mathcal{M})^*$ is $I$-complete if for all $f \in M$, if $\mathcal{M} \models (f: [\bar{\kappa}]^n \longrightarrow \mathrm{TC}(A))$ where $n \in \omega$ and $A \in I$ then there is $X \in \mathcal{U}$ such that $\mathcal{M} \models (f \textrm{ is constant on } [X]^n)$. 
\end{Definitions1}

\begin{Definitions1} \label{Df:CanonicallyRamsey}
We say that an n.p. ultrafilter $\mathcal{U} \subseteq (\mathcal{P}(\bar{\kappa})^\mathcal{M})^*$ is canonically Ramsey if for all $f \in M$, if $\mathcal{M} \models (f \textrm{ is a function with domain }[\bar{\kappa}]^n)$ where $n \in \omega$ then there is $X \in \mathcal{U}$ such that $\mathcal{M} \models (X \textrm{ is } f\textrm{-canonical})$.  
\end{Definitions1}

\begin{Definitions1} \label{Df:ITight}
We say that an n.p. ultrafilter $\mathcal{U} \subseteq (\mathcal{P}(\bar{\kappa})^\mathcal{M})^*$ is $I$-tight if for all $f \in M$, if $\mathcal{M} \models (f \textrm{ is a function with domain } [\bar{\kappa}]^n)$ where $n \in \omega$ then there is $X \in \mathcal{U}$ such that either
\begin{itemize}
\item[(a)] $\mathcal{M} \models (f \textrm{ is constant on } [X]^n)$, or
\item[(b)] there is an $I$-large cardinal $\theta$ such that
$$\mathcal{M} \models \forall x_1 \cdots \forall x_n\left(\bigwedge_{1 \leq i \leq n} (x_i \in X) \Rightarrow f(\{x_1, \ldots, x_n\}) \notin H_\theta\right).$$
\end{itemize}  
\end{Definitions1}

\noindent It is important to note that Definitions \ref{Df:IComplete}, \ref{Df:CanonicallyRamsey} and \ref{Df:ITight} only refer to functions in $\mathcal{M}$ whose domain is $[\bar{\kappa}]^n$ where $n$ is a natural number of the meta-theory. We now turn to showing that $\mathcal{P}(\bar{\kappa})$ in $\mathcal{M}$ is rich enough to ensure the existence of an (external) ultrafilter on $(\mathcal{P}(\bar{\kappa})^\mathcal{M})^*$ that is simultaneously $I$-complete, canonically Ramsey, $I$-tight and contains arbitrarily small $I$-large sets.

\begin{Lemma1} \label{Th:CompleteExistence}
If $X \in (\mathcal{P}(\bar{\kappa})^\mathcal{M})^*$ is $I$-large, $\lambda \in I$ is such that $\mathcal{M} \models (\lambda \textrm{ is a cardinal})$ and $f \in M$ is such that $\mathcal{M} \models f: \bar{\kappa} \longrightarrow H_\lambda$, then there is $Y \in (\mathcal{P}(X)^\mathcal{M})^*$ such that $Y$ is $I$-large and $f$ is constant on $Y$. 
\end{Lemma1}

\begin{proof}
Work inside $\mathcal{M}$. Let $X \subseteq \bar{\kappa}$ be $I$-large and let $\lambda$ be a cardinal in $I$. Let $f: \bar{\kappa} \longrightarrow H_\lambda$. The fact that $I$ is an $H$-cut of $\mathcal{M}$ ensures that $H_\lambda$ is in $I$ and $|H_\lambda|$ is not $I$-large. Let $\mu$ be an $I$-large regular cardinal such that $\mu < |X|$. So, $|H_\lambda| < \mu$. Suppose that for all $x \in H_\lambda$, $f^{-1}[x] \cap X$ is not $I$-large, in particular $|f^{-1}[x] \cap X| < \mu$. So,
$$|X|= \left| \bigcup_{x \in H_\lambda} f^{-1}[x] \cap X \right| \leq \mu,$$
which is a contradiction.
\Square
\end{proof}

In \cite[Corollary 2]{bau75} James Baumgartner completely determines the canonical partition relation on infinite cardinals and reveals the following relationship between the canonical and classical partition relations:

\begin{Theorems1} \label{Th:BaumgartnerTheorem}
(Baumgartner) If $\kappa$ and $\lambda$ are infinite cardinals and $n \in \omega$ then
\begin{equation} \label{eq:ClassicalVsCanonicallPartitionArrow}
\kappa \rightarrow *(\lambda)^{n+1} \textrm{ if and only if for all } \mu < \lambda, \kappa \rightarrow (\lambda)^{n+1}_\mu.
\end{equation}
\Square
\end{Theorems1}

\noindent The setting for Baumgartner's \cite{bau75} is $\mathrm{ZFC}$, however an examination of the proof of Theorem \ref{Th:BaumgartnerTheorem} reveals that for fixed $\kappa$ and $\lambda$ all the instances of separation and replacement appealed to in the proof of the equivalence described by (\ref{eq:ClassicalVsCanonicallPartitionArrow}) can be restricted to sets that provably exist in $\mathrm{MOST}$. This means that for fixed $\kappa$ and $\lambda$, the equivalence described by (\ref{eq:ClassicalVsCanonicallPartitionArrow}) is provable in $\mathrm{MOST}$. In $\mathrm{ZFC}$, Theorem \ref{Th:BaumgartnerTheorem} coupled with the Erd\H{o}s-Rado Theorem \cite{ER56} shows that for every infinite successor cardinal $\lambda$ and $n \in \omega$, there exists a cardinal $\kappa$ such that $\kappa \rightarrow *(\lambda)^{n+1}$. Specifically:

\begin{Theorems1} \label{Th:ErdosRadoTheorem}
(Erd\H{o}s-Rado) If $\kappa$ is an infinite cardinal and $n \in \omega$ then
$$\beth_n(\kappa)^+ \rightarrow (\kappa^+)_\kappa^{n+1}.$$
\Square
\end{Theorems1}

If both $\kappa^+$ and $\beth_n(\kappa)^+$ exist then all instances of separation and replacement appealed to in the proof that $\beth_n(\kappa)^+ \rightarrow (\kappa^+)_\kappa^{n+1}$ can be bounded by sets that provably exist in $\mathrm{MOST}$. Thus, as long as both $\kappa^+$ and $\beth_n(\kappa)^+$ exist, $\mathrm{MOST}$ proves that $\beth_n(\kappa)^+ \rightarrow (\kappa^+)_\kappa^{n+1}$. It should be noted, however, that $\mathrm{MOST}$ is incapable of proving that $\beth_n(\aleph_0)$ exists for all natural numbers $n$. Despite this, for any infinite cardinal $\kappa$, $\mathrm{MOST}$ proves that $2^\kappa$ exists. Therefore, if $n$ is (externally) a standard natural number and $\kappa$ is an infinite cardinal then $\mathrm{MOST}$ proves that $\beth_n(\kappa)^+$ exists. In the context of $\mathrm{MOST}$ Theorem \ref{Th:ErdosRadoTheorem} becomes:

\begin{Lemma1} \label{Th:ErdosRadoInMOST}
$\mathrm{MOST}$ proves the theorem scheme: for all $n \in \omega$,
$$\forall \kappa((\kappa \textrm{ is an infinite cardinal}) \Rightarrow (\beth_n(\kappa)^+\rightarrow (\kappa^+)_\kappa^{n+1})).$$
\end{Lemma1}

Therefore combining Theorems \ref{Th:BaumgartnerTheorem} and \ref{Th:ErdosRadoTheorem} in the context of the theory $\mathrm{MOST}$ we get:

\begin{Lemma1} \label{Th:CanonicalPartitionInMOST}
$\mathrm{MOST}$ proves the theorem scheme: for all $n \in \omega$,
$$\forall \kappa((\kappa \textrm{ is an infinite cardinal}) \Rightarrow (\beth_n(\kappa)^+\rightarrow *(\kappa^+)^{n+1})).$$
\Square
\end{Lemma1}     

\begin{Lemma1} \label{Th:BethIterate}
Let $n \in \omega$. If $\lambda \in \bar{\kappa}^*\cup \{\bar{\kappa}\}$ is $I$-large then there exists an $I$-large $\mu \in \bar{\kappa}^*$ with $\mathcal{M} \models (\mu \textrm{ is a cardinal})$ such that $\beth_n(\mu) < \lambda$.
\end{Lemma1}

\begin{proof}
Work inside $\mathcal{M}$. Assume that $\lambda \leq \bar{\kappa}$ is $I$-large and for all $I$-large cardinals $\gamma$, $\beth_n(\gamma) \geq \lambda$.  
Consider
$$C= \{ \gamma \in \bar{\kappa} \mid (\gamma \textrm{ is a cardinal}) \land (\beth_n(\gamma) < \lambda)\}.$$
By bounding all the quantifiers in the defining formula of $C$ by $H_{\bar{\kappa}^+}$ we can see that $\Delta_0$-separation implies that $C$ is a set. Our assumption about $\lambda$ implies that $C^* \subseteq I$. Therefore, since $I$ is a proper $H$-cut, $C \in I$. But this means that $C$ is bounded in $I$. But if $\gamma \in I$ is a cardinal with $\gamma \notin C^*$ then $\beth_n(\gamma) \notin I$, which contradicts the fact that $\mathcal{I} \models \mathrm{MOST}$. 
\Square
\end{proof}

\begin{Lemma1} \label{Th:fCanonicalExistence}
Let $n \in \omega$. If $X \in (\mathcal{P}(\bar{\kappa})^\mathcal{M})^*$ is $I$-large and $f \in M$ is such that $\mathcal{M} \models(f \textrm{ is a function with domain } [X]^{n+1})$ then there is an $I$-large $Y \in (\mathcal{P}(X)^\mathcal{M})^*$ that is $f$-canonical.
\end{Lemma1}

\begin{proof}
Work inside $\mathcal{M}$. Let $X \subseteq \bar{\kappa}$ be $I$-large and let $f$ be a function with domain $[X]^{n+1}$. Let $\lambda= |X|$. So, $\lambda$ is $I$-large. Using Lemma \ref{Th:BethIterate} we can find an $I$-large cardinal $\mu$ such that $\beth_n(\mu) < \lambda$. Therefore, by Lemma \ref{Th:CanonicalPartitionInMOST}, $\lambda \rightarrow *(\mu^+)^{n+1}$. And so there is an $I$-large $Y \subseteq X$ that is $f$-canonical.
\Square
\end{proof}

\begin{Lemma1} \label{Th:TightExistence}
Let $n \in \omega$. If $X \in (\mathcal{P}(\bar{\kappa})^\mathcal{M})^*$ is $I$-large and $f \in M$ is such that $\mathcal{M} \models(f \textrm{ is a function with domain } [X]^{n+1})$ then there is an $I$-large $Y \in (\mathcal{P}(X)^\mathcal{M})^*$ such that either
\begin{itemize}
\item[(a)] $\mathcal{M} \models (f \textrm{ is constant on } [Y]^{n+1})$, or 
\item[(b)] there is an $I$-large $\theta \in \bar{\kappa}^*$ with $\mathcal{M} \models (\theta \textrm{ is a cardinal})$ such that
$$\mathcal{M} \models (\forall A \in [Y]^{n+1})(f(A) \notin H_\theta).$$
\end{itemize}
\end{Lemma1}

\begin{proof}
Work inside $\mathcal{M}$. Let $X \subseteq \bar{\kappa}$ be $I$-large and let $f$ be a function with domain $[X]^{n+1}$. Let $\kappa= |X|$. Using Lemma \ref{Th:BethIterate} we can find an $I$-large cardinal $\lambda$ such that $\kappa \rightarrow (\lambda)_2^{n+1}$. Applying Lemma \ref{Th:BethIterate} we can then obtain an $I$-large cardinal $\mu^+$ such that $\lambda \rightarrow (\mu^+)_\mu^{n+1}$. Then, using Lemma \ref{Th:BethIterate} again, let $\theta$ be an $I$-large cardinal such that $2^\theta < \mu^+$. Therefore $|H_\theta| \leq \mu$. Now, define $g:[X]^{n+1} \longrightarrow 2$ such that for all $\{ x_0, \ldots, x_n \} \in [X]^{n+1}$,
$$g(\{x_0, \ldots, x_n\})= \left\{\begin{array}{ll}
0 & \textrm{if } f(\{x_0, \ldots, x_n\}) \in H_\theta\\
1 & \textrm{otherwise}
\end{array}\right.$$
Let $Z \subseteq X$ be such that $|Z| \geq \lambda$ and $g$ is constant on $[Z]^{n+1}$. If $g``[Z]^{n+1}=\{1\}$ then we are done. If $g``[Z]^{n+1}=\{0\}$ then let $f^\prime$ be the restriction of $f$ to $[Z]^{n+1}$.  Since $|\mathrm{rng}(f^\prime)| \leq |H_\theta| \leq \mu$, there is a $Y \subseteq Z$ with $|Y| \geq \mu^+$ such that $f^\prime$ is constant on $[Y]^{n+1}$.
\Square
\end{proof}

Equipped with Lemmas \ref{Th:CompleteExistence}, \ref{Th:fCanonicalExistence} and \ref{Th:TightExistence} we now show that there exists an external non-principle ultrafilter on the subsets of $\bar{\kappa}$ in $\mathcal{M}$ which is simultaneously $I$-complete, canonically Ramsey, $I$-tight and contains arbitrarily small $I$-large sets. We use $\mathcal{F}$ to denote the class of all points in $\mathcal{M}$ which correspond to a function whose domain is $[\bar{\kappa}]^n$ for some external natural number $n$. I.e.
$$\mathcal{F}= \{f \in M \mid (\exists n \in \omega)(\mathcal{M} \models (f \textrm{ is a function with domain }[\bar{\kappa}]^n))\}.$$

\begin{Theorems1} \label{Th:UltrafilterExistence}
There exists an n.p. ultrafilter $\mathcal{U} \subseteq (\mathcal{P}(\bar{\kappa})^\mathcal{M})^*$ which is $I$-complete, canonically Ramsey, $I$-tight and is such that $\{|X|^\mathcal{M} \mid X \in \mathcal{U}\}$ is downward cofinal in $\mathrm{Ord}^\mathcal{M} \backslash \mathrm{Ord}^\mathcal{I}$.
\end{Theorems1}

\begin{proof}
Let $\langle f_n \mid n \in \omega \rangle$ be an enumeration of $\mathcal{F}$. Let $\langle k_n \mid n \in \omega \rangle$ be a sequence of natural numbers such that for all $n \in \omega$, $\mathcal{M} \models (f_n \textrm{ is a function with domain } [\bar{\kappa}]^{k_n})$. Let $\langle \lambda_n \mid n \in \omega \rangle$ be a decreasing sequence of $\mathcal{M}$-cardinals that is downward cofinal in $\mathrm{Ord}^\mathcal{M} \backslash \mathrm{Ord}^\mathcal{I}$ with $\lambda_0 \leq \bar{\kappa}$. Using Lemmas \ref{Th:CompleteExistence}, \ref{Th:fCanonicalExistence} and \ref{Th:TightExistence} inductively build sequences $\langle W_n \mid n \in \omega \rangle$, $\langle X_n \mid n \in \omega \rangle$, $\langle Y_n \mid n \in \omega \rangle$ and $\langle Z_n \mid n \in \omega \rangle$ of $I$-large elements of $(\mathcal{P}(\bar{\kappa})^\mathcal{M})^*$ such that for all $n \in \omega$,
\begin{enumerate}
\item $\mathcal{M} \models (W_n \supseteq X_n \supseteq Y_n \supseteq Z_n \supseteq W_{n+1})$,
\item $\mathcal{M} \models (W_n \textrm{ is } f_n\textrm{-canonical})$,
\item if $k_n= 1$ and there is an $\mathcal{M}$-cardinal $\mu \in I$ such that $\mathcal{M} \models (f_n: \bar{\kappa} \longrightarrow H_\mu)$ then $f_n$ is constant on $X_n$, otherwise $X_n= W_n$, 
\item $\mathcal{M} \models (f_n \textrm{ is constant on } [Y_n]^{k_n})$ or there is an $I$-large $\mathcal{M}$-cardinal $\mu \in \bar{\kappa}^*$ such that 
$$\mathcal{M} \models (\forall A \in [Y_n]^{k_n})(f_n(A) \notin H_\mu),$$
\item $\mathcal{M} \models (|Z_n| < \lambda_n)$.  
\end{enumerate}
Define $\mathcal{U}= \{X \in (\mathcal{P}(\bar{\kappa})^\mathcal{M})^* \mid (\exists n \in \omega)(\mathcal{M} \models W_n \subseteq X)\}$. It is clear from the construction that $\mathcal{U}$ is an n.p. ultrafilter that is $I$-complete, canonically Ramsey, $I$-tight and is such that $\{|X|^\mathcal{M} \mid X \in \mathcal{U}\}$ is downward cofinal in $\mathrm{Ord}^\mathcal{M} \backslash \mathrm{Ord}^\mathcal{I}$.    
\Square
\end{proof}

Let $\mathcal{U} \subseteq (\mathcal{P}(\bar{\kappa})^\mathcal{M})^*$ be an n.p. ultrafilter obtained from Theorem \ref{Th:UltrafilterExistence}, so $\mathcal{U}$ is $I$-complete, canonically Ramsey, $I$-tight and $\{|X|^\mathcal{M} \mid X \in \mathcal{U}\}$ is downward cofinal in $\mathrm{Ord}^\mathcal{M} \backslash \mathrm{Ord}^\mathcal{I}$. Let $\mathcal{L}_\mathcal{F}$ be the language extending $\mathcal{L}$ such that for all $f \in \mathcal{F}$ and $n \in \omega$, if $\mathcal{M} \models (f \textrm{ has domain }[\bar{\kappa}]^n)$ then $\mathcal{L}_\mathcal{F}$ has a new $n$-ary function symbol $\hat{f}$. Let $\mathcal{M}_\mathcal{F}$ be the $\mathcal{L}_\mathcal{F}$-structure obtained from $\mathcal{M}$ by defining, for all $f \in \mathcal{F}$ with $\mathcal{M} \models (f \textrm{ has domain }[\bar{\kappa}]^n)$,
$$\hat{f}^{\mathcal{M}_\mathcal{F}}(x_1, \ldots, x_n)= \left\{\begin{array}{ll}
f(\{x_1, \ldots, x_n\})^\mathcal{M} & \textrm{if } x_1 < \cdots < x_n \in \bar{\kappa}^*\\
\emptyset & \textrm{otherwise}
\end{array}\right.$$
For each $n \in \omega$, define the $n+1$-partial type $\Gamma_n(x_0, \ldots, x_n) \subseteq \mathcal{L}_\mathcal{F}$ by
$$\begin{array}{lcl}
\phi(x_0, \ldots, x_n) \in \Gamma_n(x_0, \ldots, x_n) & \textrm{if and only if} & \exists X \in \mathcal{U} \textrm{ such that } \mathcal{M}_\mathcal{F} \models \phi(a_0, \ldots, a_n)\\
& & \textrm{for all } a_0 < \cdots < a_n \in X^*
\end{array}$$            
Define 
$$T_\mathcal{U}= \bigcup_{n \in \omega} \Gamma_n(x_0, \ldots, x_n).$$
Let $\mathcal{L}_{\mathcal{F}, \mathbb{L}}$ be the extension of $\mathcal{L}_\mathcal{F}$ obtained by adding constant symbols $c_i$ for each $i \in \mathbb{L}$. Define
$$T_{\mathcal{U}, \mathbb{L}}= \{ \phi(c_{i_0}, \ldots, c_{i_n}) \mid \phi \in T_{\mathcal{U}} \textrm{ and } i_0 < \cdots < i_n \in \mathbb{L}\}.$$
The fact that for each $f \in \mathcal{F}$, the interpretation of the function symbol $\hat{f}$ in $\mathcal{M}_\mathcal{F}$ is coded by the point $f \in M$ yields the following:

\begin{Lemma1} \label{Th:SeparationAndCollectionInF}
$\mathcal{M}_\mathcal{F} \models \Delta_0^\mathcal{P}(\mathcal{L}_\mathcal{F}) \textrm{-separation}+\Sigma_1^\mathcal{P}(\mathcal{L}_\mathcal{F})\textrm{-collection}$. \Square
\end{Lemma1}

\begin{Lemma1} \label{Th:ConsistentAndComplete}
$T_{\mathcal{U}, \mathbb{L}}$ is consistent and is $\Delta_0^{\mathcal{P}}(\mathcal{L}_{\mathcal{F}, \mathbb{L}})$-complete.
\end{Lemma1}

\begin{proof}
The fact that $T_{\mathcal{U}, \mathbb{L}}$ is finitely realisable implies that it is consistent. Let $\phi(x_0, \ldots, x_n)$ be a $\Delta_0^{\mathcal{P}}(\mathcal{L}_{\mathcal{F}})$-formula. We work inside $\mathcal{M}_{\mathcal{F}}$. $\Delta_0^{\mathcal{P}}(\mathcal{L}_{\mathcal{F}})$-separation ensures that the function $f: [\bar{\kappa}]^{n+1} \longrightarrow 2$ defined such that for all $\alpha_0 < \cdots < \alpha_n \in \bar{\kappa}$,
$$f(\alpha_0, \ldots, \alpha_n)= \left\{ \begin{array}{ll}
1 & \textrm{if } \phi(\alpha_0, \ldots, \alpha_n)\\
0 & \textrm{otherwise} 
\end{array}\right.$$
is a set. Now, there is an $X \in \mathcal{U}$ such that
$$\mathcal{M} \models (f \textrm{ is constant on } [X]^{n+1}).$$ 
Therefore, for all $i_0 < \cdots < i_n \in \mathbb{L}$, $T_{\mathcal{U}, \mathbb{L}}$ decides $\phi(c_{i_0}, \ldots, c_{i_n})$.        
\Square
\end{proof}

\begin{Remark1}
Lemma \ref{Th:ConsistentAndComplete} generalises to show that for all $n \in \omega$, if $\mathcal{M} \models \Sigma_n \textrm{-separation}$ then $T_{\mathcal{U}, \mathbb{L}}$ is $\Sigma_n(\mathcal{L}_{\mathcal{F}, \mathbb{L}})$-complete. In particular, if $\mathcal{M} \models \mathrm{ZFC}$ then $T_{\mathcal{U}, \mathbb{L}}$ is $\mathcal{L}_{\mathcal{F}, \mathbb{L}}$-complete. 
\end{Remark1}

Let $\mathrm{TERM} \subseteq \mathcal{L}_{\mathcal{F}, \mathbb{L}}$ be the class of terms of the form $\hat{f}(c_{i_0}, \ldots, c_{i_n})$ where $i_0 < \cdots < i_n \in \mathbb{L}$. Define
$$\hat{f}(c_{i_0}, \ldots, c_{i_n}) \sim \hat{g}(c_{j_0}, \ldots, c_{j_m}) \textrm{ if and only if } (\hat{f}(c_{i_0}, \ldots, c_{i_n}) = \hat{g}(c_{j_0}, \ldots, c_{j_m})) \in T_{\mathcal{U}, \mathbb{L}}.$$
Now, $\sim$ is an equivalence class on $\mathrm{TERM}$. If $\tau \in \mathrm{TERM}$ then we use $[\tau]$ to denote the equivalence class to which $\tau$ belongs. We now turn to defining an $\mathcal{L}$-structure $\mathcal{N}_{\mathcal{U}, \mathbb{L}}= \langle N_{\mathcal{U}, \mathbb{L}}, \in^{\mathcal{N}_{\mathcal{U}, \mathbb{L}}} \rangle$. Let 
$$N_{\mathcal{U}, \mathbb{L}}= \{ [\tau] \mid \tau \in \mathrm{TERM} \}.$$
Define
$$[\hat{f}(c_{i_0}, \ldots, c_{i_n})] \in^{\mathcal{N}_{\mathcal{U}, \mathbb{L}}} [\hat{g}(c_{j_0}, \ldots, c_{j_m})] \textrm{ if and only if}$$ 
$$(\hat{f}(c_{i_0}, \ldots, c_{i_n}) \in \hat{g}(c_{j_0}, \ldots, c_{j_m})) \in T_{\mathcal{U}, \mathbb{L}}.$$

\begin{Lemma1} \label{Th:SkolemLemma}
Let $\phi(x_0, \ldots, x_n)$ be a $\Delta_0(\mathcal{L}_{\mathcal{F}})$-formula. If 
$$\mathcal{M}_{\mathcal{F}} \models (\forall x_1, \ldots x_n \in \bar{\kappa})((x_1 < \cdots < x_n) \Rightarrow \exists y \phi(y, x_1, \ldots, x_n))$$
then there exists $f \in \mathcal{F}$ such that 
$$\mathcal{M}_{\mathcal{F}} \models (\forall x_1, \ldots, x_n \in \bar{\kappa})((x_1 < \cdots < x_n) \Rightarrow \phi(\hat{f}(x_1, \ldots, x_n), x_1, \ldots, x_n)).$$ 
\end{Lemma1}

\begin{proof}
Assume that $\mathcal{M}_{\mathcal{F}} \models (\forall x_1, \ldots x_n \in \bar{\kappa})((x_1 < \cdots < x_n) \Rightarrow \exists y \phi(y, x_1, \ldots, x_n))$. We work inside $\mathcal{M}_{\mathcal{F}}$. Consider $[\bar{\kappa}]^n$ and define 
$$\phi^\prime(y, x) \textrm{ if and only if } (\exists x_1, \ldots, x_n \in x)((x_1 < \cdots < x_n) \land \phi(y, x_1, \ldots, x_n)).$$
Note that $(\forall x \in [\bar{\kappa}]^n) \exists y \phi^\prime(y, x)$. Therefore $\Delta_0^\mathcal{P}(\mathcal{L}_{\mathcal{F}})$-collection ensures that there is a set $A$ such that $(\forall x \in [\bar{\kappa}]^n)(\exists y \in A) \phi^\prime(y, x)$. Let $\LHD \subseteq A \times A$ be a well-ordering of $A$. Let $\psi(y, x)$ be the $\Delta_0^{\mathcal{P}}(\mathcal{L}_{\mathcal{F}})$-formula
$$\phi^\prime(y, x) \land (\forall z \in A)(\phi^\prime(z, x) \Rightarrow \langle z, y \rangle \notin \LHD).$$
Now, $\Delta_0^\mathcal{P}(\mathcal{L}_{\mathcal{F}})$-separation ensures that $f: [\bar{\kappa}]^n \longrightarrow A$, defined such that $f(x)$ is the unique $y$ such that $\psi(y, x)$, is a set. 
\Square
\end{proof}

\begin{Lemma1} \label{Th:LosLemma}
If $\phi(x_0, \ldots, x_n)$ a $\Sigma_1^{\mathcal{P}}$-formula then for all $[\tau_0], \ldots, [\tau_n] \in N_{\mathcal{U}, \mathbb{L}}$,
$$\mathcal{N}_{\mathcal{U}, \mathbb{L}} \models \phi([\tau_0], \ldots, [\tau_n]) \textrm{ if and only if } \phi(\tau_0, \ldots, \tau_n) \in T_{\mathcal{U}, \mathbb{L}}.$$
\end{Lemma1}

\begin{proof}
Let $\phi(x_0, \ldots, x_n)$ be a $\Sigma_1^{\mathcal{P}}$-formula and let $[\tau_0], \ldots, [\tau_n] \in N_{\mathcal{U}, \mathbb{L}}$. Without loss of generality we can assume that $\phi$ only uses the quantifier $\exists$ and the logical connectives $\land$ and $\neg$. We prove the result by induction on the structure of $\phi$. It is clear that the result holds for all atomic formulae and conjunctions of atomic formulae. Assume that the result holds for a $\Delta_0^{\mathcal{P}}$-formula $\theta(x_0, \ldots, x_n)$. Lemma \ref{Th:ConsistentAndComplete} implies that the result also holds for $\neg \theta$.\\ 
We now turn to showing that the class of formula satisfying the result is closed under existential quantification. Assume that the result holds for a $\Delta_0^{\mathcal{P}}$-formula $\theta(y, x_0, \ldots, x_n)$. Let $[\tau_0], \ldots, [\tau_n] \in N_{\mathcal{U}, \mathbb{L}}$. Suppose that $\mathcal{N}_{\mathcal{U}, \mathbb{L}} \models \exists y \theta(y, [\tau_0], \ldots, [\tau_n])$. Let $[\tau] \in N_{\mathcal{U}, \mathbb{L}}$ be such that $\mathcal{N}_{\mathcal{U}, \mathbb{L}} \models \theta([\tau], [\tau_0], \ldots, [\tau_n])$. Therefore $\theta(\tau, \tau_0, \ldots, \tau_n) \in T_{\mathcal{U}, \mathbb{L}}$. This implies that $\exists y \theta(y, \tau_0, \ldots, \tau_n) \in T_{\mathcal{U}, \mathbb{L}}$.\\ 
Conversely, suppose that $\exists y \theta(y, \tau_0, \ldots, \tau_n) \in T_{\mathcal{U}, \mathbb{L}}$. Assume, without loss of generality, that for all $0 \leq j \leq n$, $\tau_j= \hat{f}_j(c_{i_0}, \ldots, c_{i_m})$ for $i_0 < \cdots < i_m$ in $\mathbb{L}$. Let $X \in \mathcal{U}$ be such that for all $a_0 < \cdots < a_m \in X^*$,
$$\mathcal{M}_{\mathcal{F}} \models \exists y \theta(y, \hat{f}_0(a_0, \ldots, a_m), \ldots, \hat{f}_n(a_0, \ldots, a_m)).$$
Let $\theta^\prime(y, x_0, \ldots, x_m)$ be the $\Delta_0^{\mathcal{P}}(\mathcal{L}_{\mathcal{F}})$-formula
$$\left(\bigwedge_{0 \leq i \leq m} (x_i \in X)\right)\land (x_0 < \cdots < x_m) \Rightarrow \theta(y, \hat{f}_0(x_0, \ldots, x_m), \ldots, \hat{f}_n(x_0, \ldots, x_m)).$$
Using Lemma \ref{Th:SkolemLemma} we can find an $f \in \mathcal{F}$ such that for all $x_0 < \cdots < x_m \in X^*$,
$$\mathcal{M}_{\mathcal{F}} \models \theta^\prime(\hat{f}(x_0, \ldots, x_m), x_0, \ldots, x_m).$$
Let $\tau= \hat{f}(c_{i_0}, \ldots, c_{i_m})$. Now, $\theta(\tau, \tau_0, \ldots, \tau_n) \in T_{\mathcal{U}, \mathbb{L}}$. Therefore $\mathcal{N}_{\mathcal{U}, \mathbb{L}} \models \theta([\tau], [\tau_0], \ldots, [\tau_n])$ and so $\mathcal{N}_{\mathcal{U}, \mathbb{L}} \models \exists y \theta(y, [\tau_0], \ldots, [\tau_n])$. 
\Square
\end{proof}

\begin{Remark1}
The proof of Lemma \ref{Th:LosLemma} generalises to show that for all $n \in \omega$, if $\mathcal{M} \models \Sigma_n\textrm{-collection}$ then the conclusion of Lemma \ref{Th:LosLemma} holds for all $\Sigma_n$-formulae. In particular, if $\mathcal{M} \models \mathrm{ZFC}$ then $\mathcal{M} \prec \mathcal{N}_{\mathcal{U}, \mathbb{L}}$.
\end{Remark1}

The structure $\mathcal{M}$ embeds into $\mathcal{N}_{\mathcal{U}, \mathbb{L}}$. To see this observe that for all $m \in M$ there exists an $h_m \in \mathcal{F}$ such that 
$$\mathcal{M}_{\mathcal{F}} \models (\forall x \in \bar{\kappa})(\hat{h}_m(x)= m).$$
It is clear that for all $i, j \in \mathbb{L}$ and for all $m \in M$, $[\hat{h}_m(c_i)]$ is equal to $[\hat{h}_m(c_j)]$. In order to make it easier to refer to these elements of $N_{\mathcal{U}, \mathbb{L}}$ we will fix an element $0 \in \mathbb{L}$ so that we can write $[\hat{h}_m(c_0)]$. We can also see that the linear order $\mathbb{L}$ embeds in the ordinals of $\mathcal{N}_{\mathcal{U}, \mathbb{L}}$. To see this consider the map $\mathrm{id} \in \mathcal{F}$ such that
$$\mathcal{M}_{\mathcal{F}} \models (\hat{\mathrm{id}}: \bar{\kappa} \longrightarrow \bar{\kappa}) \land (\forall x \in \bar{\kappa})(\hat{\mathrm{id}}(x)=x).$$
For each $i \in \mathbb{L}$, the term $[\hat{\mathrm{id}}(c_i)]$ is an ordinal in $\mathcal{N}_{\mathcal{U}, \mathbb{L}}$.

\begin{Lemma1} \label{Th:Indiscernibles}
The class $\{ [\hat{\mathrm{id}}(c_i)] \mid i \in \mathbb{L} \} \subseteq N_{\mathcal{U}, \mathbb{L}}$ is a class of $\Sigma_1^{\mathcal{P}}$-indiscernibles.
\end{Lemma1}

\begin{proof}
This follows immediately from Lemma \ref{Th:LosLemma}.
\Square
\end{proof}

The fact that the image of the embedding of $\mathbb{L}$ into $\mathcal{N}_{\mathcal{U}, \mathbb{L}}$ is a class of $\Sigma_1^{\mathcal{P}}$-indiscernibles means that any automorphism of $\mathbb{L}$ raises to an automorphism of $\mathcal{N}_{\mathcal{U}, \mathbb{L}}$. For all $j \in \mathrm{Aut}(\mathbb{L})$, define $\check{j}: \mathcal{N}_{\mathcal{U}, \mathbb{L}} \longrightarrow \mathcal{N}_{\mathcal{U}, \mathbb{L}}$ by
\begin{equation} \label{eq:AutomorphismOfN}
\check{j}([\hat{f}(c_{i_1}, \ldots c_{i_n})])= [\hat{f}(c_{j(i_1)}, \ldots c_{j(i_n)})] \textrm{ for all } f \in \mathcal{F} \textrm{ and } i_1 < \cdots < i_n \in \mathbb{L}. 
\end{equation}
Lemma \ref{Th:Indiscernibles} implies that for all $j \in \mathrm{Aut}(\mathbb{L})$, $\check{j}: \mathcal{N}_{\mathcal{U}, \mathbb{L}} \longrightarrow \mathcal{N}_{\mathcal{U}, \mathbb{L}}$ is an automorphism. Moreover, the map from $\mathrm{Aut}(\mathbb{L})$ into $\mathrm{Aut}(\mathcal{N}_{\mathcal{U}, \mathbb{L}})$ defined by $j \mapsto \check{j}$ for all $j \in \mathrm{Aut}(\mathbb{L})$, is an injective group homomorphism. It is also immediate from (\ref{eq:AutomorphismOfN}) that for all $j \in \mathrm{Aut}(\mathbb{L})$ and for all $x \in M$, $\check{j}$ fixes $[\hat{h}_x(c_0)]$, so $\mathrm{fix}(\check{j}) \supseteq M$.\\ 
\\
\indent Lemma \ref{Th:LosLemma} also implies that embedding of $\mathcal{M}$ into $\mathcal{N}_{\mathcal{U}, \mathbb{L}}$ defined by $x \mapsto [\hat{h}_x(c_0)]$ preserves $\Pi_2^\mathcal{P}$ properties of tuples from $\mathcal{M}$.

\begin{Lemma1} \label{Th:Pi2PPropertiesPreserved}
Let $\phi(x_0, \ldots, x_n)$ be a $\Pi_2^\mathcal{P}$-formula. For all $a_0, \ldots, a_n \in M$, if $\mathcal{M} \models \phi(a_0, \ldots, a_n)$ then
$$\mathcal{N}_{\mathcal{U}, \mathbb{L}} \models \phi([\hat{h}_{a_0}(c_0)], \ldots, [\hat{h}_{a_n}(c_0)]).$$
\end{Lemma1}

\begin{proof}
The formula $\phi(x_0, \ldots, x_n)$ is in the form $\forall y_0 \cdots \forall y_m \theta(y_0, \ldots, y_m, x_0, \ldots, x_n)$ where $\theta(y_0, \ldots, y_m, x_0, \ldots, x_n)$ is $\Sigma_1^\mathcal{P}$. Let $a_0, \ldots, a_n \in M$. Assume that
$$\mathcal{M} \models \forall y_0 \cdots \forall y_m \theta(y_0, \ldots, y_m, a_0, \ldots, a_n).$$
Let $f_0, \ldots, f_m \in \mathcal{F}$. Without loss of generality assume that each $\hat{f}_0, \ldots, \hat{f}_m$ are $k+1$-ary function symbols in $\mathcal{L}_\mathcal{F}$. For all $x_0 < \cdots < x_k \in \bar{\kappa}^*$,
$$\mathcal{M}_\mathcal{F} \models \theta(\hat{f}_0(x_0, \ldots, x_k), \ldots, \hat{f}_m(x_0, \ldots, x_k), a_0, \ldots, a_n).$$
Therefore, for all $i_0 < \cdots < i_k \in \mathbb{L}$,
$$\theta(\hat{f}_0(c_{i_0}, \ldots, c_{i_k}), \ldots, \hat{f}_m(c_{i_0}, \ldots, c_{i_k}), \hat{h}_{a_0}(c_0), \ldots, \hat{h}_{a_n}(c_0)) \in T_{\mathcal{U}, \mathbb{L}}.$$
And so, by Lemma \ref{Th:LosLemma},
$$\mathcal{N}_{\mathcal{U}, \mathbb{L}} \models \theta([\hat{f}_0(c_{i_0}, \ldots, c_{i_k})], \ldots, [\hat{f}_m(c_{i_0}, \ldots, c_{i_k})], [\hat{h}_{a_0}(c_0)], \ldots, [\hat{h}_{a_n}(c_0)]).$$
Since $f_0, \ldots, f_m \in \mathcal{F}$ and $i_0 < \cdots < i_k \in \mathbb{L}$ were arbitrary, the Lemma follows.     
\Square
\end{proof}

Using this result we can show that $\mathcal{N}_{\mathcal{U}, \mathbb{L}}$ satisfies all of the axioms of $\mathrm{Mac}$.

\begin{Lemma1} \label{Th:MacInN}
$\mathcal{N}_{\mathcal{U}, \mathbb{L}} \models \mathrm{Mac}$.
\end{Lemma1}

\begin{proof}
By Lemma \ref{Th:Pi2PPropertiesPreserved} it is enough to show that every axiom of $\mathrm{Mac}$ can be written as a $\Pi_2^\mathcal{P}$-sentence.
\begin{itemize}
\item[] Extensionality: $\forall x \forall y(x = y \iff (\forall z \in x)(z \in y) \land (\forall z \in y)(z \in x))$
\item[] Emptyset: $\exists x (\forall y \in x)(y \neq y)$
\item[] Union: $\forall x \exists y (\forall z \in x)(\forall w \in z)(w \in y)$ 
\item[] Pairing: $\forall x \forall y \exists z ((x \in z) \land (y \in z) \land (\forall w \in z)(w= x \lor w = y))$ 
\item[] $\Delta_0$-separation: for all $\Delta_0$-formulae $\phi(x, \vec{z})$,
$$\forall a \forall \vec{z} (\exists y \subseteq a)(\forall x \in a)(x \in y \iff \phi(x, \vec{z})).$$
\item[] Set Foundation: $\forall x ((\forall w \in x)(w \neq w) \lor (\exists y \in x)(\forall z \in y)(z \notin x))$
\item[] Infinity: 
$$\exists S( \emptyset \in S \land (\forall x \in S)(\exists y \in S)(x \in y \land (\forall z \in x)(z \in y) \land (\forall z \in y)(z= x \lor z \in x)))$$
\item[] Powerset: $\forall x \exists y ((\forall z \subseteq x)(z \in y) \land (\forall z \in y)(\forall w \in z)(w \in x))$
\item[] Transitive Containment: $\forall x \exists y (\bigcup y \subseteq y \land (\forall z \in x)(z \in y))$
\item[] Axiom of Choice: $\forall x \exists R(R \textrm{ is a well-ordering of } x)$
\end{itemize}
\Square
\end{proof}

Before showing that $\mathcal{N}_{\mathcal{U}, \mathbb{L}}$ satisfies Axiom $H$ we first show that $\mathcal{N}_{\mathcal{U}, \mathbb{L}}$ is an end-extension of $\mathcal{I}$ and a cofinal extension of $\mathcal{M}$.

\begin{Lemma1} \label{Th:NEndExtendsI}
For all $x \in I$ and for all $[\tau] \in N_{\mathcal{U}, \mathbb{L}}$, if $\mathcal{N}_{\mathcal{U}, \mathbb{L}} \models ([\tau] \in [\hat{h}_x(c_0)])$ then there exists $y \in I$ such that
$$\mathcal{N}_{\mathcal{U}, \mathbb{L}} \models ([\tau]= [\hat{h}_y(c_0)]) \textrm{ and } \mathcal{M} \models (y \in x).$$
\end{Lemma1}

\begin{proof}
Let $x \in I$ and let $[\tau] \in N_{\mathcal{U}, \mathbb{L}}$. Suppose that $\tau= \hat{f}(c_{j_1}, \ldots, c_{j_n})$ where $f \in \mathcal{F}$ and $j_1 < \cdots < j_n \in \mathbb{L}$, and $\mathcal{N}_{\mathcal{U}, \mathbb{L}} \models ([\tau] \in [\hat{h}_x(c_0)])$. Therefore, by Lemma \ref{Th:LosLemma}, there exists $X \in \mathcal{U}$ such that for all $x_1 < \cdots < x_n \in X^*$,
$$\mathcal{M}_\mathcal{F} \models \hat{f}(x_1, \ldots, x_n) \in x.$$
Since $\mathcal{U}$ is $I$-complete, there is $Y \in \mathcal{U}$ and $y \in x^*$ such that for all $x_1 < \cdots < x_n \in Y^*$,
$$\mathcal{M}_\mathcal{F} \models \hat{f}(x_1, \ldots, x_n)=y.$$ 
Therefore $(\hat{f}(c_{j_1}, \ldots, c_{j_n})= \hat{h}_y(c_0)) \in T_{\mathcal{U}, \mathbb{L}}$. So by Lemma \ref{Th:LosLemma}, $\mathcal{N}_{\mathcal{U}, \mathbb{L}} \models ([\tau]= [\hat{h}_y(c_0)])$;
and since $y \in x^*$, $\mathcal{M} \models (y \in x)$. 
\Square
\end{proof}

Lemma \ref{Th:NEndExtendsI} shows that for all $j \in \mathrm{Aut}(\mathbb{L})$, $I_{\mathrm{fix}(\check{j})} \supseteq I$.

\begin{Lemma1} \label{Th:MCofinalInN}
For all $[\tau] \in N_{\mathcal{U}, \mathbb{L}}$, there exists $x \in M$ such that,
$$\mathcal{N}_{\mathcal{U}, \mathbb{L}} \models [\tau] \in [\hat{h}_x(c_0)].$$
\end{Lemma1}

\begin{proof}
Let $[\tau] \in N_{\mathcal{U}, \mathbb{L}}$. Suppose that $\tau= \hat{f}(c_{j_1}, \ldots, c_{j_n})$ where $f \in M$ and $j_1 < \cdots < j_n \in \mathbb{L}$. Let $x \in M$ be such that $\mathcal{M} \models (\mathrm{rng}(f) \subseteq x)$. Now, for all $x_1 < \cdots < x_n \in \bar{\kappa}^*$,
$$\mathcal{M}_\mathcal{F} \models \hat{f}(x_1, \ldots, x_n) \in x.$$
Therefore, $(\hat{f}(c_{j_1}, \ldots, c_{j_n}) \in \hat{h}_x(c_0)) \in T_{\mathcal{U}, \mathbb{L}}$. So by Lemma \ref{Th:LosLemma},
$$\mathcal{N}_{\mathcal{U}, \mathbb{L}} \models [\tau] \in [\hat{h}_x(c_0)].$$ 
\Square
\end{proof}

The fact that $\mathcal{U}$ is $I$-tight also ensures that $\mathrm{Card}^\mathcal{M} \backslash \mathrm{Card}^\mathcal{I}$ is downward cofinal in $\mathrm{Card}^{\mathcal{N}_{\mathcal{U}, \mathbb{L}}} \backslash \mathrm{Card}^\mathcal{I}$.

\begin{Lemma1} \label{Th:MCardinalsDownwardCofinalInN}
If $\lambda \in \mathrm{Card}^{\mathcal{N}_{\mathcal{U}, \mathbb{L}}} \backslash \mathrm{Card}^\mathcal{I}$ then there is $\mu \in \mathrm{Card}^\mathcal{M} \backslash \mathrm{Card}^\mathcal{I}$ with $\mathcal{N}_{\mathcal{U}, \mathbb{L}} \models ([\hat{h}_\mu(c_0)] \leq \lambda)$.
\end{Lemma1}

\begin{proof}
Let $\lambda \in \mathrm{Card}^{\mathcal{N}_{\mathcal{U}, \mathbb{L}}} \backslash \mathrm{Card}^\mathcal{I}$. Therefore $\lambda= [\hat{f}(c_{i_1}, \ldots, c_{i_n})]$ where $f \in \mathcal{F}$ and $i_1 < \cdots < i_n \in \mathbb{L}$, and $[\hat{f}(c_{i_1}, \ldots, c_{i_n})] \notin I$ and 
$$\mathcal{N}_{\mathcal{U}, \mathbb{L}} \models ([\hat{f}(c_{i_1}, \ldots, c_{i_n})] \textrm{ is a cardinal}).$$
If there is $\mu \in M$ such that 
$$\mathcal{N}_{\mathcal{U}, \mathbb{L}} \models ([\hat{f}(c_{i_1}, \ldots, c_{i_n})]= [\hat{h}_\mu(c_0)])$$
then, by Lemma \ref{Th:Pi2PPropertiesPreserved}, $[\hat{f}(c_{i_1}, \ldots, c_{i_n})] \in \mathrm{Card}^\mathcal{M} \backslash \mathrm{Card}^\mathcal{I}$ and we are done. Therefore, assume that $[\hat{f}(c_{i_1}, \ldots, c_{i_n})] \in N_{\mathcal{U}, \mathbb{L}} \backslash M$. Since $\mathcal{U}$ is $I$-tight and $[\hat{f}(c_{i_1}, \ldots, c_{i_n})] \notin I$, there is an $I$-large $\mathcal{M}$-cardinal $\mu \in M$ and $X \in \mathcal{U}$ such that for all $x_1 < \cdots < x_n \in X^*$,
$$\mathcal{M}_\mathcal{F} \models (\hat{f}(x_1, \ldots, x_n) \notin H_\mu).$$
Since $[\hat{f}(c_{i_1}, \ldots, c_{i_n})]$ is an $\mathcal{N}_{\mathcal{U}, \mathbb{L}}$-cardinal, this implies that for all $x < x_1 < \cdots < x_n \in X^*$,
$$\mathcal{M}_\mathcal{F} \models (\hat{h}_\mu(x) \leq \hat{f}(x_1, \ldots, x_n)).$$
And so, by Lemma \ref{Th:LosLemma}, $\mathcal{N}_{\mathcal{U}, \mathbb{L}} \models (\hat{h}_\mu(c_0) \leq \hat{f}(c_{i_1}, \ldots, c_{i_n}))$.   
\Square
\end{proof}  

We now turn to showing that $\mathcal{N}_{\mathcal{U}, \mathbb{L}}$ satisfies Axiom $H$.

\begin{Lemma1} \label{Th:InitialOrdinalsExistInN}
$\mathcal{N}_{\mathcal{U}, \mathbb{L}} \models \forall x(|x| \textrm{ exists})$.
\end{Lemma1}

\begin{proof}
Let $\phi(x)$ be the formula
$$\exists f \exists \alpha((\alpha \textrm{ is an ordinal}) \land (f: x \longrightarrow \alpha \textrm{ is an injection})).$$
Now, $\phi(x)$ is $\Sigma_1$ and $\mathcal{M} \models \forall x \phi(x)$. Therefore, by Lemma \ref{Th:Pi2PPropertiesPreserved}, $\mathcal{N}_{\mathcal{U}, \mathbb{L}} \models \forall x \phi(x)$ and the Lemma follows.
\Square
\end{proof}

\begin{Lemma1} \label{Th:AxiomHinN}
$\mathcal{N}_{\mathcal{U}, \mathbb{L}} \models \textrm{Axiom } \mathrm{H}$.
\end{Lemma1}

\begin{proof}
Let $u \in N_{\mathcal{U}, \mathbb{L}}$. Using Lemma \ref{Th:InitialOrdinalsExistInN}, let $\lambda \in N_{\mathcal{U}, \mathbb{L}}$ be such that 
$$\mathcal{N}_{\mathcal{U}, \mathbb{L}} \models (\lambda \textrm{ is a cardinal}) \land (|u| \leq \lambda).$$
By Lemma \ref{Th:MCofinalInN} there is $\mu \in M$ such that $\mu$ is an $\mathcal{M}$-cardinal and $\mathcal{N}_{\mathcal{U}, \mathbb{L}} \models (\lambda \leq \hat{h}_\mu(c_0))$. Let $T \in M$ be such that $\mathcal{M} \models (T= H_{\mu^+})$. Therefore
$$\mathcal{M} \models \forall z\left(\left(\bigcup z \subseteq z \land |z| \leq \mu\right) \Rightarrow z \subseteq T\right).$$
So, by Lemma \ref{Th:Pi2PPropertiesPreserved},
$$\mathcal{N}_{\mathcal{U}, \mathbb{L}} \models \forall z \left(\left( \bigcup z \subseteq z \land |z| \leq [\hat{h}_\mu(c_0)]\right) \Rightarrow z \subseteq [\hat{h}_T(c_0)]\right)$$
$$\textrm{and } \mathcal{N}_{\mathcal{U}, \mathbb{L}} \models \left(\bigcup [\hat{h}_T(c_0)] \subseteq [\hat{h}_T(c_0)]\right).$$
In particular
$$\mathcal{N}_{\mathcal{U}, \mathbb{L}} \models \forall z \left(\left( \bigcup z \subseteq z \land |z| \leq |u|\right) \Rightarrow z \subseteq [\hat{h}_T(c_0)]\right).$$
And this shows that
$$\mathcal{N}_{\mathcal{U}, \mathbb{L}} \models \forall u \exists T \left(\bigcup T \subseteq T \land \forall z \left(\left( \bigcup z \subseteq z \land |z| \leq |u|\right) \Rightarrow z \subseteq T\right)\right).$$
\Square
\end{proof}

Combining Lemma \ref{Th:MacInN} and Lemma \ref{Th:AxiomHinN} we get: 

\begin{Theorems1}
$\mathcal{N}_{\mathcal{U}, \mathbb{L}} \models \mathrm{MOST}$. \Square
\end{Theorems1}

We now turn to showing that if $j \in \mathrm{Aut}(\mathbb{L})$ has no fixed points then $\mathcal{I}_\mathrm{fix}(\check{j}) = \mathcal{I}$ and $\mathrm{fix}(\check{j})= M$.

\begin{Theorems1} \label{Th:FixjIsM}
If $j \in \mathrm{Aut}(\mathbb{L})$ has no fixed points then $\mathrm{fix}(\check{j})=M$.
\end{Theorems1}

\begin{proof}
Let $j \in \mathrm{Aut}(\mathbb{L})$ have no fixed points. It follows immediately from (\ref{eq:AutomorphismOfN}) that for all $x \in M$, $\check{j}$ fixes $[\hat{h}_x(c_0)]$. Therefore, we need to show that if $[\tau] \in N_{\mathcal{U}, \mathbb{L}}$ and $\check{j}([\tau])=[\tau]$ then there exists $x \in M$ such that $[\tau]= [\hat{h}_x(c_0)]$. Let $\tau= \hat{f}(c_{i_0}, \ldots, c_{i_n})$ where $f \in \mathcal{F}$ and $i_0 < \cdots < i_n \in \mathbb{L}$. Assume that $\check{j}([\tau])= [\tau]$. Therefore $[\hat{f}(c_{i_0}, \ldots, c_{i_n})]= [\hat{f}(c_{j(i_0)}, \ldots, c_{j(i_n)})]$. Since $\mathcal{U}$ is canonically Ramsey, there is $X \in \mathcal{U}$ and $\Delta \subseteq n+1$ such that for all $x_0 < \cdots < x_n$ and $y_0 < \cdots < y_n$ in $X^*$,
$$(\mathcal{M}_\mathcal{F} \models \hat{f}(x_0, \ldots, x_n)= \hat{f}(y_0, \ldots, y_n)) \textrm{ if and only if } (\forall m \in \Delta)(x_m=y_m).$$
Therefore, for all $j_0 < \cdots < j_n$ and $k_0 < \cdots <k_n$ in $\mathbb{L}$,
$$[\hat{f}(c_{j_0}, \ldots, c_{j_n})]=[\hat{f}(c_{k_0}, \ldots, c_{k_n})] \textrm{ if and only if } (\forall m \in \Delta)(j_m=k_m).$$
Therefore, since $j$ has no fixed points, $\Delta= \emptyset$. And so there exists $y \in M$ such that for all $x_0 < \cdots < x_n \in X^*$,
$$\mathcal{M}_\mathcal{F} \models \hat{f}(x_0, \ldots, x_n)=y.$$
Therefore $\mathcal{N}_{\mathcal{U}, \mathbb{L}} \models ([\tau]= [\hat{h}_y(c_0)])$. 
\Square
\end{proof}

\begin{Lemma1} \label{Th:TCsArePreserved}
If $x \in M$ then $\mathcal{N}_{\mathcal{U}, \mathbb{L}} \models (\mathrm{TC}([\hat{h}_x(c_0)])=[\hat{h}_{\mathrm{TC}(x)}(c_0)])$.
\end{Lemma1}

\begin{proof}
Let $x \in M$. Now, for all $y \in \bar{\kappa}^*$,
$$\mathcal{M}_{\mathcal{F}} \models \mathrm{TC}(\hat{h}_x(y))= \hat{h}_{\mathrm{TC}(x)}(y).$$
Therefore $(\mathrm{TC}(\hat{h}_x(c_0))= \hat{h}_{\mathrm{TC}(x)}(c_0)) \in T_{\mathcal{U}, \mathbb{L}}$. Now, ``$Y$ is the transitive closure of $X$" can be expressed by a $\Delta_0^{\mathcal{P}}$-formula:
$$X \subseteq Y \land \bigcup Y \subseteq Y \land (\forall x \subseteq Y)\left(X \subseteq x \land \bigcup x \subseteq x \Rightarrow Y \subseteq x\right).$$
Therefore, by Lemma \ref{Th:LosLemma}, 
$$\mathcal{N}_{\mathcal{U}, \mathbb{L}} \models \mathrm{TC}([\hat{h}_x(c_0)])= [\hat{h}_{\mathrm{TC}(x)}(c_0)].$$ 
\Square
\end{proof}

\begin{Theorems1} \label{Th:IFixIsI}
If $j \in \mathrm{Aut}(\mathbb{L})$ has no fixed points then $\mathcal{I}_{\mathrm{fix}}(\check{j}) = \mathcal{I}$.
\end{Theorems1}

\begin{proof}
Let $j \in \mathrm{Aut}(\mathbb{L})$ have no fixed points. Lemma \ref{Th:NEndExtendsI} implies that the embedding of $\mathcal{M}$ into $\mathcal{N}_{\mathcal{U}, \mathbb{L}}$ embeds $\mathcal{I}$ into an initial segment of $\mathcal{N}_{\mathcal{U}, \mathbb{L}}$. Therefore if $x \in I$ and $y \in N_{\mathcal{U}, \mathbb{L}}$ is such that $\mathcal{N}_{\mathcal{U}, \mathbb{L}} \models (y \in \mathrm{TC}(\{x\}))$ then $\check{j}(y)= y$. Therefore, we need to show that if $x \in N_{\mathcal{U}, \mathbb{L}} \backslash I$ then there is some $y \in (\mathrm{TC}(\{x\})^{\mathcal{N}_{\mathcal{U}, \mathbb{L}}})^*$ that is moved by $\check{j}$. In light of Theorem \ref{Th:FixjIsM} it is enough to show that for all $x \in N_{\mathcal{U}, \mathbb{L}} \backslash I$, there exists $y \in N_{\mathcal{U}, \mathbb{L}} \backslash M$ with $\mathcal{N}_{\mathcal{U}, \mathbb{L}} \models (y \in \mathrm{TC}(\{x\}))$. Let $x \in N_{\mathcal{U}, \mathbb{L}} \backslash I$. If $x \in N_{\mathcal{U}, \mathbb{L}} \backslash M$ then we are done, so assume that $x \in M \backslash I$. It follows  that $|\mathrm{TC}(\{x\})|^\mathcal{M} \notin I$; hence $|\mathrm{TC}(\{x\})|^\mathcal{M}$ is an $I$-large $\mathcal{M}$-cardinal.  Therefore there is $X \in \mathcal{U}$ such that $\mathcal{M} \models (|X| \leq |\mathrm{TC}(\{x\})|)$. Let $g \in M$ be such that
$$\mathcal{M} \models (g: X \longrightarrow \mathrm{TC}(\{x\})) \land (g \textrm{ is an injection}).$$
Therefore there exists $f \in M$ such that $\mathcal{M} \models (f: \bar{\kappa} \longrightarrow \mathrm{TC}(\{x\}))$ and for all $z_1 < z_2 \in X^*$,
$$\mathcal{M}_\mathcal{F} \models \hat{f}(z_1) \neq \hat{f}(z_2).$$
Therefore, for all $w \in M$, $\mathcal{N}_{\mathcal{U}, \mathbb{L}} \models ([\hat{f}(c_0)] \neq [\hat{h}_w(c_0)])$. We need to show that $[\hat{f}(c_0)]$ is in the transitive closure of $[\hat{h}_{\{x\}}(c_0)]$ in $\mathcal{N}_{\mathcal{U}, \mathbb{L}}$. For all $z \in X^*$,
$$\mathcal{M}_\mathcal{F} \models (\hat{f}(z) \in \hat{h}_{\mathrm{TC}(\{x\})}(z)).$$
Therefore, for all $i \in \mathbb{L}$, $(\hat{f}(c_i) \in \hat{h}_{\mathrm{TC}(\{x\})}(c_i)) \in T_{\mathcal{U}, \mathbb{L}}$. So, by Lemma \ref{Th:LosLemma},
$$\mathcal{N}_{\mathcal{U}, \mathbb{L}} \models [\hat{f}(c_0)] \in [\hat{h}_{\mathrm{TC}(\{x\})}(c_0)].$$
Therefore, by Lemma \ref{Th:TCsArePreserved},
$$\mathcal{N}_{\mathcal{U}, \mathbb{L}} \models [\hat{f}(c_0)] \in \mathrm{TC}([\hat{h}_{x}(c_0)]),$$
which completes the proof of the theorem.      
\Square
\end{proof}

This completes the proof of Theorem \ref{Th:MainRealisingHCutTheorem}.

\section[Realizing countable models as $H$-cuts]{Realizing countable models as $H$-cuts}

This section tackles the question of which countable $\mathcal{L}$-structures can be realised as an $H$-cut of a model of set theory. We show:
\begin{itemize}
\item every countable transitive model of $\mathrm{MOST}+\Delta_0^\mathcal{P}\textrm{-collection}$ can be realized as an $H$-cut of a model of $\mathrm{MOST}+\Delta_0^\mathcal{P}\textrm{-collection}$,   
\item every countable recursively saturated model of $\mathrm{MOST}+\Delta_0^\mathcal{P}\textrm{-collection}$ can be realised as an $H$-cut of a model of $\mathrm{MOST}+\Delta_0^\mathcal{P}\textrm{-collection}$,
\item every countable model of $\mathrm{ZFC}$ can be realised as an $H$-cut of a model of $\mathrm{ZFC}$. 
\end{itemize}
Combined with the results of section \ref{Sec:RealisingHcutsAsIFix} we then have:
\begin{enumerate}
\item every countable transitive model of $\mathrm{MOST}+\Delta_0^\mathcal{P}\textrm{-collection}$ can be realised as $\mathcal{I}_{\mathrm{fix}}(j)$ for some non-trivial automorphism $j$ of a model of $\mathrm{MOST}$,
\item every countable recursively saturated model of $\mathrm{MOST}+\Delta_0^\mathcal{P}\textrm{-collection}$ can be realised as $\mathcal{I}_{\mathrm{fix}}(j)$ for some non-trivial automorphism $j$ of a model of $\mathrm{MOST}$,
\item every countable model of $\mathrm{ZFC}$ can be realised as $\mathcal{I}_{\mathrm{fix}}(j)$ for some non-trivial automorphism $j$ of a model of $\mathrm{ZFC}$.
\end{enumerate}
Together with the results of Section \ref{Sec:StructureOfIFix}, (1) and (2) yield a complete characterisation of the countable transitive and countable recursively saturated $\mathcal{L}$-structures satisfying infinity which can be realised as $\mathcal{I}_{\mathrm{fix}}(j)$ for some non-trivial automorphism $j$ of a model of $\mathrm{MOST}$. Since every consistent theory with an infinite model has a countable recursively saturated model (Theorem \ref{Th:ExistenceOfRecursivelySaturatedModels}), (2) implies that \emph{the theory of the class of $\mathcal{L}$-structures satisfying infinity that appear as $\mathcal{I}_{\mathrm{fix}}(j)$ for some non-trivial automorphism $j$ of a model of $\mathrm{MOST}$ is exactly $\mathrm{MOST}+\Delta_0^\mathcal{P}\textrm{-collection}$.}

\subsection[Countable transitive models of $\mathrm{MOST}+\Delta_0^\mathcal{P}\textrm{-collection}$]{Countable transitive models of $\mathrm{MOST}+\Delta_0^\mathcal{P}\textrm{-collection}$} \label{Sec:CountableTransModels}

The realisation of a countable transitive model of $\mathrm{MOST}+\Delta_0^\mathcal{P}\textrm{-collection}$ as an $H$-cut of a model of $\mathrm{MOST}+\Delta_0^\mathcal{P}\textrm{-collection}$ will be achieved by applying a theorem, due to Harvey Friedman \cite{fri73}, which combined with a result of Mathias~\cite{mat01} shows that every countable transitive model of $\mathrm{MOST}+\Delta_0^\mathcal{P}\textrm{-collection}$ can be toplessly end-extended to a model of $\mathrm{MOST}+\Delta_0^\mathcal{P}\textrm{-collection}$.\\
\\
In \cite{fri73} Friedman studies countable transitive models of a theory that he calls Power Admissible Set Theory ($\mathrm{PAdm}^s$). This theory is axiomatised using classes of $\mathcal{L}$-formulae that Friedman calls $\Delta_0^s(\mathcal{P})$. Let $\mathcal{L}_\mathcal{P}$ be the extension of the $\mathcal{L}$ obtained by adding a new unary function symbol $\mathcal{P}$. Define $\mathrm{pseudo}\textrm{-}\Delta_0^s(\mathcal{P})$ to be the class of $\Delta_0(\mathcal{L}_\mathcal{P})$ formulae that contain quantification in the form $\exists x \in y$ or $\forall x \in y$ where $x$ and $y$ are distinct \emph{variables}. The class $\Delta_0^s(\mathcal{P})$ is obtained by translating the formulae in $\mathrm{pseudo}\textrm{-}\Delta_0^s(\mathcal{P})$ into $\mathcal{L}$-formulae using translations generated by the defining axiom $x \in \mathcal{P}(y) \iff \forall z(z \in x \Rightarrow z \in y)$ \footnote{an explicit list of translations can be found in \cite{mat01}}.          

\begin{Definitions1}
Power Admissible Set Theory ($\mathrm{PAdm}^s$) is obtained from $\mathrm{KP}$ by adding powerset, $\Delta_0^s(\mathcal{P})$-collection and $\mathcal{L}$-foundation.  
\end{Definitions1}

Friedman \cite[Section 1]{fri73} also introduces a class theory that corresponds to Power Admissible Set Theory. We use $\mathcal{L}_{\mathrm{Cl}}$ to denote the two-sorted extension of $\mathcal{L}$ with set variables $x, y, z, \ldots$ and class variables $X, Y, Z, \ldots$. The well-formed formulae of $\mathcal{L}_{\mathrm{Cl}}$ are built inductively from atomic formulae in the form $x \in y$, $x=y$, $X=Y$, $x=Y$ and $y \in X$ using the connectives and quantifiers of first-order logic. The class $\Delta_0^c$ is the smallest class of $\mathcal{L}_{\mathrm{Cl}}$-formulae that contains all atomic formulae in the form $x=y$, $x=Y$, $x \in y$ and $y \in X$, contains all compound formulae formed using the connectives of first-order logic, and is closed under quantification in the form $\exists x \in y$ and $\forall x \in y$ where $x$ and $y$ are distinct variables. The class $\Sigma^c$ is the smallest class of $\mathcal{L}_{\mathrm{Cl}}$-formulae that contains all $\Delta_0^c$-formulae, contains all compound formulae formed using the connectives $\land$ and $\lor$
of first-order logic, and is closed under quantification in the form $\exists x$ and $\forall x \in y$ where $x$ and $y$ are distinct variables.     

\begin{Definitions1}
Power Admissible Class Theory ($\mathrm{PAdm}^c$) is the $\mathcal{L}_{\mathrm{Cl}}$-theory with axioms: extensionality for both sets and classes, $\forall x \exists Y(x=Y)$, pairing and union for sets, powerset for sets, $\mathcal{L}$-foundation, and the following:
\begin{itemize}
\item[]($\Delta_0^c$-separation) for all $\Delta_0^c$-formulae $\phi(x, \vec{Z})$,
$$\forall \vec{Z} \forall w \exists y \forall x(x \in y \iff x \in w \land \phi(x, \vec{Z}))$$ 
\item[]($\Delta_0^c$-collection) for all $\Delta_0^c$-formulae $\phi(x, y, \vec{Z})$,
$$\forall \vec{Z} \forall w((\forall x \in w) \exists y \phi(x, y, \vec{Z}) \Rightarrow \exists C (\forall x \in w)(\exists y \in C) \phi(x, y, \vec{Z}))$$ 
\item[]($\Delta^c$-CA) for all $\Sigma^c$-formulae $\phi(x, \vec{Z})$ and $\psi(x, \vec{Z})$,
$$\forall \vec{Z} (\forall x (\phi(x, \vec{Z}) \iff \neg \psi(x, \vec{Z})) \Rightarrow \exists Y \forall x(x \in Y \iff \phi(x, \vec{Z})))$$
\item[](Class Powerset)
$$\exists X \forall y(y \in X \iff \exists w \exists z(y=\langle w, z \rangle \land \forall v(v \in z \iff \forall u(u \in v \Rightarrow u \in w))))$$  
\end{itemize} 
\end{Definitions1}

Friedman \cite[Theorem 1.6]{fri73} notes that $\mathrm{PAdm}^c$ is a conservative extension of $\mathrm{PAdm}^s$:

\begin{Theorems1} \label{Th:PAdmClassConservativeOverPAdmSet}
(Friedman) If $\mathcal{M}= \langle M, \in^\mathcal{M} \rangle$ is an $\mathcal{L}$-structure with $\mathcal{M} \models \mathrm{PAdm}^s$ then there exists $C$ such that $\langle M, C, \in^\mathcal{M} \rangle \models \mathrm{PAdm}^c$. \Square 
\end{Theorems1}    

On the other hand, as shown by Mathias~\cite[Metatheorem 6.20]{mat01} we have:

\begin{Theorems1} \label{Th:KPPandPAdm}
(Mathias) The theories $\mathrm{PAdm}^s+\textrm{Infinity}$ and $\mathrm{KP}^\mathcal{P}$ have the same transitive models. \Square
\end{Theorems1}

Note that apart from instances of the $\Pi_1^\mathcal{P}$-foundation scheme, every axiom of $\mathrm{KP}^\mathcal{P}$ is also an axiom of $\mathrm{MOST}+\Delta_0^\mathcal{P}\textrm{-collection}$. Therefore, every transitive model of $\mathrm{MOST}+\Delta_0^\mathcal{P}\textrm{-collection}$ is a model of $\mathrm{KP}^\mathcal{P}$.\\  
\\
\indent Using a version of the Barwise Compactness Theorem, Friedman \cite[Theorem 2.3]{fri73} shows that any countable transitive model of $\mathrm{PAdm}^s$ has a topless powerset-preserving end extension that is a model of $\mathrm{PAdm}^s$.

\begin{Theorems1} \label{Th:FriedmanStandardPart}
(Friedman) Let $\mathcal{M}= \langle M, C, \in \rangle$ be countable and transitive with $\mathcal{M} \models \mathrm{PAdm}^c$. If $T \in C$ is an $\mathcal{L}$-theory with $\langle M, \in \rangle \models T$ then there is an $\mathcal{L}$-structure $\mathcal{N}= \langle N, \in^\mathcal{N} \rangle$ such that 
\begin{itemize}
\item[(i)] $N \neq M$,
\item[(ii)] $\langle M, \in \rangle \subseteq_{\mathrm{topless}}^{\mathcal{P}} \mathcal{N}$,
\item[(iii)] $\mathcal{N} \models T$.
\end{itemize}
\Square   
\end{Theorems1}

Combining Theorems \ref{Th:PAdmClassConservativeOverPAdmSet}, \ref{Th:KPPandPAdm} and \ref{Th:FriedmanStandardPart} yields:

\begin{Coroll1}
If $\langle I, \in\rangle$ is a countable transitive model of $\mathrm{MOST}+\Delta_0^\mathcal{P}\textrm{-collection}$ then there is $\mathcal{M}=\langle M, \in^\mathcal{M}\rangle$ with $\mathcal{M} \models \mathrm{MOST}+\Delta_0^\mathcal{P}\textrm{-collection}+ \mathcal{L}\textrm{-foundation}$ and $I \subseteq M$ is an $H$-cut of $\mathcal{M}$. \Square
\end{Coroll1}

Combining this with the results of sections \ref{Sec:StructureOfIFix} and \ref{Sec:RealisingHcutsAsIFix} gives a characterisation of the countable transitive structures $\langle I, \in\rangle$ with $\omega \in I$ that can be realised as $\mathcal{I}_{\mathrm{fix}}(j)$ for some non-trivial automorphism $j$ of a model of $\mathrm{MOST}$. 

\begin{Theorems1}
Let $\langle I, \in \rangle$ be a countable transitive structure with $\omega \in I$. The following are equivalent:
\begin{itemize}
\item[(I)] $\langle I, \in \rangle \models \mathrm{MOST}+\Delta_0^\mathcal{P}\textrm{-collection}$  
\item[(II)] there is an $\mathcal{L}$-structure $\mathcal{M} \models \mathrm{MOST}$ and a non-trivial automorphism $j: \mathcal{M} \longrightarrow \mathcal{M}$ such that
$$\mathcal{I}_{\mathrm{fix}}(j) = \langle I, \in \rangle.$$
\end{itemize}
\Square
\end{Theorems1}

\subsection[Countable recursively saturated models of $\mathrm{MOST}+\Delta_0^\mathcal{P}\textrm{-collection}$]{Countable recursively saturated models of $\mathrm{MOST}+\Delta_0^\mathcal{P}\textrm{-collection}$} \label{Sec:CountableRecSatModels}

In this section we will show that every countable recursively saturated model of $\mathrm{MOST}+\Delta_0^\mathcal{P}\textrm{-collection}$ can be realised as an $H$-cut of a model of $\mathrm{MOST}+\Delta_0^\mathcal{P}\textrm{-collection}$. This will be achieved by proving the following refined version of Friedman's Self-Embedding Theorem \cite[Section 4]{fri73} for non-standard models of set theory:

\begin{Theorems1} \label{Th:SelfEmbeddingTheorem}
If $\mathcal{M}=\langle M, \in^\mathcal{M} \rangle$ is a countable recursively saturated model of $\mathrm{MOST}+\Delta_0^\mathcal{P}\textrm{-collection}$ then there exists an embedding $h: \mathcal{M} \longrightarrow \mathcal{M}$ such that $\mathrm{rng}(h) \subseteq M$ is an $H$-cut of $\mathcal{M}$.  
\end{Theorems1} 

Combined with sections \ref{Sec:StructureOfIFix} and \ref{Sec:RealisingHcutsAsIFix} Theorem \ref{Th:SelfEmbeddingTheorem} yields a characterisation of the countable recursively saturated structures satisfying infinity that can be realised as $\mathcal{I}_{\mathrm{fix}}(j)$ for some non-trivial automorphism $j$ of a model of $\mathrm{MOST}$. 

\begin{Theorems1}
Let $\mathcal{I}= \langle I, \in^\mathcal{I} \rangle$ be a countable recursively saturated structure with $\mathcal{I} \models \textrm{Infinity}$. The following are equivalent:
\begin{itemize}
\item[(I)] $\mathcal{I} \models \mathrm{MOST}+\Delta_0^\mathcal{P}\textrm{-collection}.$
\item[(II)] there is an $\mathcal{L}$-structure $\mathcal{M} \models \mathrm{MOST}$ and a non-trivial automorphism $j: \mathcal{M} \longrightarrow \mathcal{M}$ such that
$$\mathcal{I}_{\mathrm{fix}}(j) = \mathcal{I}.$$
\end{itemize}
\Square
\end{Theorems1}

Combined with the observation that every consistent theory $T$ extending $\mathrm{MOST}+\Delta_0^\mathcal{P}\textrm{-collection}$ has a countable recursively saturated model (Theorem \ref{Th:ExistenceOfRecursivelySaturatedModels}), Theorem \ref{Th:SelfEmbeddingTheorem} also yields the first-order theory of the class of $\mathcal{L}$-structures that can appear as $\mathcal{I}_{\mathrm{fix}}(j)$ for some non-trivial automorphism $j$ of a model of $\mathrm{MOST}$.

\begin{Theorems1}
Let $T$ be a complete, consistent $\mathcal{L}$-theory such that $T \vdash \textrm{Infinity}$. The following are equivalent:
\begin{itemize}
\item[(I)] $T \vdash \mathrm{MOST}+\Delta_0^\mathcal{P}\textrm{-collection}.$
\item[(II)] there is an $\mathcal{L}$-structure $\mathcal{M} \models \mathrm{MOST}$ and a non-trivial automorphism $j: \mathcal{M} \longrightarrow \mathcal{M}$ such that
$$\mathcal{I}_{\mathrm{fix}}(j) \models T.$$ 
\end{itemize}
\Square
\end{Theorems1}
   
We now turn to proving Theorem \ref{Th:SelfEmbeddingTheorem}. For the remainder of this section let $\mathcal{M}=\langle M, \in^\mathcal{M} \rangle$ be a countable recursively saturated structure with $\mathcal{M} \models \mathrm{MOST}+\Delta_0^\mathcal{P}\textrm{-collection}$. We need to build an embedding $h: \mathcal{M} \longrightarrow \mathcal{M}$ such that $\mathrm{rng}(h) \subseteq M$ is an $H$-cut of $\mathcal{M}$. This will be achieved by a three-stage back-and-forth construction in which the `back' and `forth' steps are similar to the proofs of the classical versions of Friedman's Self-Embedding Theorem for set theory and arithmetic (e.g., as in \cite[Section 4]{fri73} and \cite[Chapter 12]{kay91}). The `third' stage of the back-and-forth construction will be used to ensure that the range of the embedding sits toplessly inside $\mathcal{M}$. Let $\langle m_i \mid i \in \omega \rangle$ be an enumeration of $M$ in which each element of $M$ appears infinitely often. Let $\langle \lambda_i \mid i \in \omega \rangle$ be an enumeration of the class
$$\{\lambda \in M \mid \mathcal{M} \models (\lambda \textrm{ is a limit cardinal})\}.$$

\begin{Lemma1} \label{Th:BackAndForthStartup}
There exists $X \in M$ such that
$$\mathcal{M} \models \exists \kappa((\kappa \textrm{ is a cardinal})\land \forall y(|\mathrm{TC}(y)| < \kappa \iff y \in X))$$
and for all $\Sigma_1^\mathcal{P}$-sentences $\phi$, if $\mathcal{M} \models \phi$ then $\mathcal{M} \models \phi^X$.
\end{Lemma1}

\begin{proof}
We use the fact that $\mathcal{M}$ is recursively saturated. Let $\Gamma(x)$ be the one-type that consists of the following formulae:
\begin{itemize}
\item[(i)] $\exists \kappa ((\kappa \textrm{ is a cardinal}) \land \forall y(|\mathrm{TC}(y)| < \kappa \iff y \in x))$ 
\item[(ii)] for all $\Sigma_1^\mathcal{P}$-sentences $\phi$,
$$\phi \Rightarrow \phi^x.$$
\end{itemize}
Note that $\Gamma(x)$ is a recursive type. We need to show that $\Gamma(x)$ is finitely realised. Suppose that $\Delta(x) \subseteq \Gamma(x)$ is finite and that the instances of (ii) mentioned in $\Delta(x)$ are exactly 
$$\psi_i \Rightarrow \psi_i^x \textrm{ where } \psi_i \textrm{ is a } \Sigma_1^\mathcal{P} \textrm{-sentence for } 0 \leq i < k.$$
Without loss of generality we may assume that for all $0 \leq i < k$, $\mathcal{M} \models \psi_i$. Suppose that for all $0 \leq i < k$, $\psi_i$ is the sentence $\exists z \theta_i(z)$ where $\theta_i(z)$ is a $\Delta_0^\mathcal{P}$-formula. Let $a_0, \ldots, a_{k-1} \in M$ be such that for all $0 \leq i < k$,
$$\mathcal{M} \models \theta_i(a_i).$$
Now, work inside $\mathcal{M}$. Let $\kappa= \sup\{|\mathrm{TC}(\mathcal{P}(a_i))| \mid 0 \leq i < k\}$. The fact that we are working in a model of $\mathrm{MOST}$ ensures that $H_{\kappa^+}$ exists. Since $a_0, \ldots, a_{k-1}, \mathcal{P}(a_0), \ldots, \mathcal{P}(a_{k-1}) \in H_{\kappa^+}$ and $\theta_0(z), \ldots, \theta_{k-1}(z)$ are $\Delta_0^\mathcal{P}$, it follows that
$$\langle H_{\kappa^+}, \in \rangle \models \theta_i(a_i) \textrm{ for all } 0 \leq i < k.$$
Since $H_{\kappa^+}^\mathcal{M}$ also satisfies (i), $H_{\kappa^+}^\mathcal{M}$ realises $\Delta(x)$ in $\mathcal{M}$. This shows that $\Gamma(x)$ is finitely realised. Since $\mathcal{M}$ is recursively saturated, there is $X \in M$ that realises $\Gamma(x)$. This proves the lemma.  
\Square
\end{proof}

Lemma \ref{Th:BackAndForthStartup} allows us to initiate the back-and-forth proof of Theorem \ref{Th:SelfEmbeddingTheorem}. Let $X_0 \in M$ be such that
$$\exists \kappa ((\kappa \textrm{ is a cardinal}) \land \forall y(|\mathrm{TC}(y)| < \kappa \iff y \in X_0)) \textrm{ and}$$
$$\textrm{for all } \Sigma_1^\mathcal{P} \textrm{-sentences } \phi, \textrm{ if } \mathcal{M} \models \phi \textrm{ then } \mathcal{M} \models \phi^{X_0}.$$
We will construct an embedding $h: \mathcal{M} \longrightarrow \mathcal{M}$ by constructing sequences $\langle u_i \mid i \in \omega \rangle$, $\langle v_i \mid i \in \omega \rangle$ and $\langle X_i \mid i \in \omega \rangle$ of elements of $M$,
such that for all $i, j \in \omega$: if $i < 2j$ then $\mathcal{M} \models v_i \in X_j$,
and if $i < j$ then $\mathcal{M} \models X_j \subseteq X_i$.  Then we define:
$$h(u_i)= v_i \textrm{ for all } i \in \omega.$$
At stage $j$, after having defined $u_0, \ldots, u_{2j-1}, v_0, \ldots, v_{2j-1}$ and $X_j$, we will ensure that the following condition is maintained:
\begin{itemize}
\item[]($\dagger_j$) for all $\Sigma_1^\mathcal{P}$-formulae $\phi(x_0, \ldots, x_{2j-1})$,
$$\textrm{if } \mathcal{M} \models \phi(u_0, \ldots, u_{2j-1}) \textrm{ then } \mathcal{M} \models \phi^{X_j}(v_0, \ldots, v_{2j-1}).$$
\end{itemize}
Suppose that we have chosen $u_0, \ldots, u_{2k-1}, v_0, \ldots, v_{2k-1}$ and $X_k$ and maintained $(\dagger_k)$. Stage $k$ of the construction comprises three steps:\\
\textbf{Step 1:} This step will ensure that the image of $h$ sits toplessly inside $\mathcal{M}$. Our aim at stage $k$ is to prevent $(H_{\lambda_k}^\mathcal{M})^*$ from being the image of $h$ (where $\lambda_k$ is as in the definition preceding Lemma~\ref{Th:BackAndForthStartup}). Consider the following conditions:
\begin{itemize}
\item[(a)] $\mathcal{M} \models (H_{\lambda_k} \subseteq X_k)$
\item[(b)] $\langle (H_{\lambda_k}^\mathcal{M})^*, \in^\mathcal{M} \rangle \models \mathrm{MOST}$
\item[(c)] $v_0, \ldots, v_{2k-1} \in (H_{\lambda_k}^\mathcal{M})^*$
\item[(d)] for all $\Sigma_1^\mathcal{P}$-formulae $\phi(x_0, \ldots, x_{2k-1})$,
$$\textrm{if } \mathcal{M} \models \phi(u_0, \ldots, u_{2k-1}) \textrm{ then } \mathcal{M} \models \phi^{H_{\lambda_k}}(v_0, \ldots, v_{2k-1}).$$
\end{itemize}
If any of the conditions (a), (b), (c) or (d) fail then it is already impossible for $(H_{\lambda_k}^\mathcal{M})^*$ to be the image of $h$. Therefore, if any of (a), (b), (c) or (d) fail then let $X_{k+1}=X_k$. Since $\dagger_k$ holds, for all $\Sigma_1^\mathcal{P}$-formulae $\phi(x_0, \ldots, x_{2k-1})$,
\begin{equation} \label{eq:ConditionAfterStep1}
\textrm{if } \mathcal{M} \models \phi(u_0, \ldots, u_{2k-1}) \textrm{ then } \mathcal{M} \models \phi^{X_{k+1}}(v_0, \ldots, v_{2k-1}).
\end{equation}     
If (a), (b), (c) and (d) all hold then we will choose $X_{k+1}$ so that condition (a) fails for $X_k$. This will prevent $(H_{\lambda_k}^\mathcal{M})^*$ from being the image of $h$.

\begin{Lemma1} \label{Th:RefineXLemma}
If conditions (a), (b), (c) and (d) all hold then there exists $X \in (H_{\lambda_k}^\mathcal{M})^*$ such that
$$\mathcal{M} \models \exists \kappa((\kappa \textrm{ is a cardinal}) \land \forall y (|\mathrm{TC}(y)| < \kappa \iff y \in X))$$
and for all $\Sigma_1^\mathcal{P}$-formulae $\phi(x_0, \ldots, x_{2k-1})$,
$$\textrm{if } \mathcal{M} \models \phi(u_0, \ldots, u_{2k-1}) \textrm{ then } \mathcal{M} \models \phi^X(v_0, \ldots, v_{2k-1}).$$ 
\end{Lemma1}

\begin{proof}
Assume that conditions (a), (b), (c) and (d) all hold. Let $\Gamma(x)$ be the one-type that contains the following formulae:
\begin{itemize}
\item[(i)] $x \in H_{\lambda_k}$,
\item[(ii)] $\exists \kappa((\kappa \textrm{ is a cardinal}) \land \forall y(|\mathrm{TC}(y)| < \kappa \iff y \in x))$,
\item[(iii)] for all $\Sigma_1^\mathcal{P}$-formulae $\phi(x_0, \ldots, x_{2k-1})$,
$$\phi(u_0, \ldots, u_{2k-1}) \Rightarrow \phi^x(v_0, \ldots, v_{2k-1}).$$
\end{itemize}
$\Gamma(x)$ is a recursive type. We need to show that $\Gamma(x)$ is finitely realised. Suppose that $\Delta(x) \subseteq \Gamma(x)$ is finite and that the instances of (iii) mentioned in $\Delta(x)$ are exactly 
$$\psi_i(u_0, \ldots, u_{2k-1}) \Rightarrow \psi_i^x(v_0, \ldots, v_{2k-1}) \textrm{ where } \psi_i \textrm{ is a } \Sigma_1^\mathcal{P} \textrm{-formula for } 0 \leq i < m.$$
Without loss of generality we may assume that for all $0 \leq i < m$,
$$\mathcal{M} \models \psi_i(u_0, \ldots, u_{2k-1}).$$
Suppose that for each $0 \leq i < m$, the formula $\psi_i(x_0, \ldots, x_{2k-1})$ is $\exists z \theta_i(z, x_0, \ldots, x_{2k-1})$ where $\theta_i(z, x_0, \ldots, x_{2k-1})$ is a $\Delta_0^\mathcal{P}$-formula. It follows from (d) that for all $0 \leq i < m$, 
$$\mathcal{M} \models \psi_i^{H_{\lambda_k}}(v_0, \ldots, v_{2k-1}).$$
Let $a_0, \ldots, a_{m-1} \in (H_{\lambda_k}^\mathcal{M})^*$ be such that for all $0 \leq i < m$,
$$\mathcal{M} \models \theta_i^{H_{\lambda_k}}(a_i, v_0, \ldots, v_{2k-1}).$$
Work inside $\mathcal{M}$. Let $\mu_1 = \sup\{|\mathrm{TC}(\mathcal{P}(a_i))| \mid 0 \leq i < m\}$ and let $\mu_2= \sup\{|\mathrm{TC}(\mathcal{P}(v_i))| \mid 0 \leq i < 2k\}$. Let $\kappa= \max\{\mu_1, \mu_2\}$. The fact that we are working in a model of $\mathrm{MOST}$ ensures that $H_{\kappa^+}$ exists. It follows from condition (b) that $H_{\kappa^+} \in H_{\lambda_k}$. Now,
$$a_0, \ldots, a_{m-1}, \mathcal{P}(a_0), \ldots, \mathcal{P}(a_{m-1}), v_0, \ldots, v_{2k-1}, \mathcal{P}(v_0), \ldots, \mathcal{P}(v_{2k-1}) \in H_{\kappa^+}.$$
Therefore, since each $\theta_i(z, x_0, \ldots, x_{2k-1})$ is $\Delta_0^\mathcal{P}$,
$$\langle H_{\kappa^+}, \in \rangle \models \theta_i(a_i, v_0, \ldots, v_{2k-1}) \textrm{ for all } 0 \leq i < m.$$
Therefore $H_{\kappa^+}^\mathcal{M}$ realizes $\Delta(x)$ in $\mathcal{M}$. This shows that $\Gamma(x)$ is finitely realized. Since $\mathcal{M}$ is recursively saturated, it follows that there is an $X \in M$ that realizes $\Gamma(x)$. This proves the lemma. 
\Square
\end{proof}

Let $X_{k+1}$ be the point in $\mathcal{M}$ guaranteed by Lemma \ref{Th:RefineXLemma}. It follows that
$$\mathcal{M} \models \neg(H_{\lambda_k} \subseteq X_{k+1}).$$
This prevents $(H_{\lambda_k}^\mathcal{M})^*$ from being the image of $h$. We also have that for all $\Sigma_1^\mathcal{P}$-formulae $\phi(x_0, \ldots, x_{2k-1})$, (\ref{eq:ConditionAfterStep1}) holds. This completes Step 1.\\
\\
\textbf{Step 2:} This is the usual `forth' step in the proof of Friedman's Embedding Theorem (see \cite[Theorem 12.3]{kay91}). Let $u_{2k}= m_k$. This choice will eventually ensure that the domain of $h$ is all of $M$. We need to choose $v_{2k} \in X_{k+1}^*$ such that for all $\Sigma_1^\mathcal{P}$-formulae $\phi(x_0, \ldots, x_{2k})$,
$$\textrm{if } \mathcal{M} \models \phi(u_0, \ldots, u_{2k}) \textrm{ then } \mathcal{M} \models \phi(v_0, \ldots, v_{2k}).$$
The following Lemma shows that we can successfully make this choice:

\begin{Lemma1} \label{Th:ForthLemma}
There exists $v \in X_{k+1}^*$ such that for all $\Sigma_1^\mathcal{P}$-formulae $\phi(x_0, \ldots, x_{2k})$,
$$\textrm{if } \mathcal{M} \models \phi(u_0, \ldots, u_{2k}) \textrm{ then } \mathcal{M} \models \phi(v_0, \ldots, v_{2k-1}, v).$$
\end{Lemma1}

\begin{proof}
We use the fact that $\mathcal{M}$ is recursively saturated. Let $\Gamma(x)$ be the one-type that contains the following formulae:
\begin{itemize}
\item[(i)] $x \in X_{k+1}$, 
\item[(ii)] for all $\Sigma_1^\mathcal{P}$-formulae $\phi(x_0, \ldots, x_{2k})$,
$$\phi(u_0, \ldots, u_{2k}) \Rightarrow \phi^{X_{k+1}}(v_0, \ldots, v_{2k-1}, x).$$ 
\end{itemize}
$\Gamma(x)$ is a recursive type. We need to show that $\Gamma(x)$ is finitely realised. Suppose that $\Delta(x) \subseteq \Gamma(x)$ is finite and that the instances of (ii) mentioned in $\Delta(x)$ are exactly 
$$\psi_i(u_0, \ldots, u_{2k}) \Rightarrow \psi_i^{X_{k+1}}(v_0, \ldots, v_{2k-1}, x) \textrm{, where } \psi_i \textrm{ is a } \Sigma_1^\mathcal{P} \textrm{-formula for } 0 \leq i < m.$$
Without loss of generality we may assume that for all $0 \leq i < m$,
$$\mathcal{M} \models \psi_i(u_0, \ldots, u_{2k}).$$
Suppose that for each $0 \leq i < m$, the formulae $\psi_i(x_0, \ldots, x_{2k})$ is $\exists z_i \theta_i(z_i, x_0, \ldots, x_{2k})$ where $\theta_i(z_i, x_0, \ldots, x_{2k})$ is a $\Delta_0^\mathcal{P}$-formula and without loss of generality the $z_i$ are distinct.  We have
$$\mathcal{M} \models \exists y \exists z_0 \cdots \exists z_{m-1} \bigwedge_{0 \leq i < m} \theta_i(z_i, u_0, \ldots, u_{2k-1}, y).$$
Since for all $\Sigma_1^\mathcal{P}$-formulae $\phi(x_0, \ldots, x_{2k-1})$, (\ref{eq:ConditionAfterStep1}) holds, it follows that
$$\mathcal{M} \models (\exists y \in X_{k+1})(\exists z_0 \in X_{k+1}) \cdots (\exists z_{m-1} \in X_{k+1}) \bigwedge_{0 \leq i < m} \theta_i^{X_{k+1}}(z_i, v_0, \ldots, v_{2k-1}, y).$$
Let $v \in X_{k+1}^*$ be such that for all $0 \leq i < m$,
$$\mathcal{M} \models \psi_i^{X_{k+1}}(v_0, \ldots, v_{2k-1}, v).$$
Therefore $v \in M$ realizes $\Delta(x)$. This shows that $\Gamma(x)$ is finitely realised. Since $\mathcal{M}$ is recursively saturated, there is a $v \in M$ which realizes $\Gamma(x)$. This proves the lemma.  
\Square
\end{proof}

Let $v \in X_{k+1}^*$ be the point in $\mathcal{M}$ guaranteed by Lemma \ref{Th:ForthLemma}. Let
$$v_{2k}=\left\{\begin{array}{ll}
v_i & \textrm{if } u_{2k}= u_i \textrm{ for some } 0 \leq i < 2k,\\
v & \textrm{otherwise}
\end{array}\right.$$
This ensures that for all $\Sigma_1^\mathcal{P}$-formulae $\phi(x_0, \ldots, x_{2k})$,
\begin{equation} \label{eq:ConditionAfterStep2}
\textrm{if } \mathcal{M} \models \phi(u_0, \ldots, u_{2k}) \textrm{ then } \mathcal{M} \models \phi^{X_{k+1}}(v_0, \ldots, v_{2k}).
\end{equation}
This completes Step 2.\\
\\
\textbf{Step 3:} This is the usual `back' step in the proof of Friedman's Embedding Theorem (see \cite[Theorem 12.3]{kay91}). This step will eventually ensure that $\mathcal{M}$ is a powerset-preserving end-extension of the image of $h$. In this step we have two cases to consider:\\
\textbf{Case 1:} For all $0 \leq i < 2k+1$, $m_k \nsubseteq v_i$. In this case let $v_{2k+1}=v_0$ and let $u_{2k+1}=u_0$. This choice clearly satisfies ($\dagger_{k+1}$).\\
\textbf{Case 2:} There exists $0 \leq i < 2k+1$ such that $m_k \subseteq v_i$. Note that in this case it immediately follows that $m_k \in X_{k+1}^*$. Let $v_{2k+1} = m_k$. We need to choose $u_{2k+1}$ such that for all $\Sigma_1^\mathcal{P}$-formulae $\phi(x_0, \ldots, x_{2k+1})$,
$$\textrm{if } \mathcal{M} \models \phi(u_0, \ldots, u_{2k+1}) \textrm{ then } \mathcal{M} \models \phi^{X_{k+1}}(v_0, \ldots, v_{2k+1}).$$

\begin{Lemma1} \label{Th:PowersetForParameters}
For all $0 \leq j < 2k+1$, $\mathcal{P}^\mathcal{M}(v_j) \in X_{k+1}^*$.
\end{Lemma1}

\begin{proof}
Let $0 \leq j < 2k+1$. Since $\mathcal{M} \models \mathrm{MOST}$, it follows that
$$\mathcal{M} \models \exists y (\forall x \subseteq u_j(x \in y) \land \forall x \in y(x \subseteq u_j)).$$
By (\ref{eq:ConditionAfterStep2}):
$$\mathcal{M} \models  (\exists y \in X_{k+1})(\forall x \subseteq v_j(x \in y) \land \forall x \in y(x \subseteq v_j))^{X_{k+1}}.$$
Since every subset of $v_j$ is a member of $X_{k+1}$, it follows that $\mathcal{P}^\mathcal{M}(v_j) \in X_{k+1}^*$. 
\Square
\end{proof}

The following lemma ensures that we can choose $u_{2k+1}$ to satisfy ($\dagger_{k+1}$):

\begin{Lemma1} \label{Th:BackLemma}
There exists $u \in M$ such that for all $\Sigma_1^\mathcal{P}$-formulae $\phi(x_0, \ldots, x_{2k+1})$,
$$\textrm{if } \mathcal{M} \models \phi(u_0, \ldots, u_{2k}, u) \textrm{ then } \mathcal{M} \models \phi^{X_{k+1}}(v_0, \ldots, v_{2k+1}).$$
\end{Lemma1}

\begin{proof}
We use the fact that $\mathcal{M}$ is recursively saturated. Let $\Gamma(x)$ be the one-type that contains the following formulae:
\begin{itemize}
\item[(i)] $x \subseteq u_i$
\item[(ii)] for all $\Delta_0^\mathcal{P}$-formulae $\phi(z, x_0, \ldots, x_{2k+1})$,
$$(\forall z \in X_{k+1}) \phi^{X_{k+1}}(z, v_0, \ldots, v_{2k+1}) \Rightarrow \forall z \phi(z, u_0, \ldots, u_{2k}, x).$$
\end{itemize}
$\Gamma(x)$ is a recursive type. We need to show that $\Gamma(x)$ is finitely realised. Suppose that $\Delta(x) \subseteq \Gamma(x)$ is finite and that the instances of (ii) mentioned in $\Delta(x)$ are exactly 
$$(\forall z \in X_{k+1}) \psi_j^{X_{k+1}}(z, v_0, \ldots, v_{2k+1}) \Rightarrow \forall z \psi_j(z, u_0, \ldots, u_{2k}, x)$$ 
$$\textrm{where } \psi_j \textrm{ is a } \Delta_0^\mathcal{P} \textrm{-formula for } 0 \leq j < m.$$
Without loss of generality we may assume that for all $0 \leq j < m$,
$$\mathcal{M} \models (\forall z \in X_{k+1}) \psi_j^{X_{k+1}}(z, v_0, \ldots v_{2k+1}).$$
Suppose, for a contradiction, that
$$\mathcal{M} \models (\forall x \subseteq u_i) \exists z \bigvee_{0 \leq j < m} \neg \psi_j(z, u_0, \ldots, u_{2k}, x).$$
Therefore
$$\mathcal{M} \models (\forall x \in \mathcal{P}(u_i)) \exists z \bigvee_{0 \leq j < m} \neg \psi_j(z, u_0, \ldots, u_{2k}, x).$$
By applying $\Delta_0^\mathcal{P}$-collection we can conclude that
$$\mathcal{M} \models \exists t (\forall x \in \mathcal{P}(u_i))(\exists z \in t) \bigvee_{0 \leq j < m} \neg \psi_j(z, u_0, \ldots, u_{2k}, x).$$
Therefore
$$\mathcal{M} \models \exists t (\forall x \subseteq u_i)(\exists z \in t) \bigvee_{0 \leq j < m} \neg \psi_j(z, u_0, \ldots, u_{2k}, x).$$
By (\ref{eq:ConditionAfterStep2}) and Lemma \ref{Th:PowersetForParameters}, we have
\begin{equation} \label{eq:BackEq1}
\mathcal{M} \models (\exists t \in X_{k+1})(\forall x \subseteq v_i)(\exists z \in t) \bigvee_{0 \leq j < m} \neg \psi_j^{X_{k+1}}(z, v_0, \ldots, v_{2k}, x).
\end{equation}
But $v_{2k+1} \subseteq v_i$ and
$$\mathcal{M} \models (\forall z \in X_{k+1}) \bigwedge_{0\leq j < m} \neg \psi_j^{X_{k+1}}(z, v_0, \ldots, v_{2k+1})$$
which contradicts (\ref{eq:BackEq1}). Therefore $\Gamma(x)$ is finitely satisfied. Since $\mathcal{M}$ is recursively saturated, there is a $u \in M$ that realises $\Gamma(x)$. We claim that this $u \in M$ is the point that is required by the lemma. Suppose that this is not the case and that $\phi(x_0, \ldots, x_{2k+1})$ is a $\Sigma_1^\mathcal{P}$-formula such that
$$\mathcal{M} \models \phi(u_0, \ldots, u_{2k}, u) \textrm{ and } \mathcal{M} \models \neg \phi^{X_{k+1}}(v_0, \ldots, v_{2k+1}).$$
The formula $\phi(x_0, \ldots, x_{2k+1})$ is equivalent to a formula $\exists z \theta(z, x_0, \ldots, x_{2k+1})$ where $\theta(z, x_0, \ldots, x_{2k+1})$ is a $\Delta_0^\mathcal{P}$-formula. Therefore
$$\mathcal{M} \models (\forall z \in X_{k+1}) \neg \theta^{X_{k+1}}(z, v_0, \ldots, v_{2k+1}) \textrm{ and } \mathcal{M} \models \neg \forall z \neg \theta(z, u_0, \ldots, u_{2k}, u).$$
But this contradicts the fact that $u$ realizes $\Gamma(x)$. This proves the lemma.   
\Square
\end{proof}

Let $u \in M$ be the point guaranteed by Lemma \ref{Th:BackLemma}. Let
$$u_{2k+1}=\left\{\begin{array}{ll}
u_j & \textrm{if } v_{2k+1}= v_j \textrm{ for some } 0 \leq j < 2k+1,\\
u & \textrm{otherwise}
\end{array}\right.$$
It follows from Lemma \ref{Th:BackLemma} that this choice of $u_{2k+1}$ satisfies ($\dagger_{k+1}$).\\
\\
This completes the $k^{th}$ stage of the back-and-forth construction. Continuing this process yields a map $h: \mathcal{M} \longrightarrow \mathcal{M}$. The fact that ($\dagger_n$) holds at the beginning of each stage $n \in \omega$ ensures that $h: \mathcal{M} \longrightarrow \mathcal{M}$ is an embedding. Step 2 of each stage ensures that the domain of $h: \mathcal{M} \longrightarrow \mathcal{M}$ is all of $M$. Step 3 of each stage ensures that 
$$\mathrm{rng}(h) \subseteq_e^\mathcal{P} \mathcal{M}.$$
Step 1 of each stage ensures both that $\mathrm{rng}(h) \neq M$ and that $\mathrm{rng}(h)$ sits toplessly in $\mathcal{M}$. This completes the proof of Theorem \ref{Th:SelfEmbeddingTheorem}.

\subsection[Countable models of $\mathrm{ZFC}$]{Countable models of $\mathrm{ZFC}$} 

Sections \ref{Sec:CountableTransModels} and \ref{Sec:CountableRecSatModels} show that if $\mathcal{I}$ is countable with $\mathcal{I} \models \mathrm{MOST}+\Delta_0^\mathcal{P}\textrm{-collection}$ and $\mathcal{I}$ is either transitive or recursively saturated then $\mathcal{I}$ can be realised as $\mathcal{I}_\mathrm{fix}(j)$ for some non-trivial automorphism $j$ of a model of $\mathrm{MOST}$. This raises the following:

\begin{Quest1} \label{Q:CanWeClassifyCountableIfix}
Can every countable model of $\mathrm{MOST}+\Delta_0^\mathcal{P}\textrm{-collection}$ be realised as $\mathcal{I}_\mathrm{fix}(j)$ for some non-trivial automorphism $j$ of a model of $\mathrm{MOST}$?
\end{Quest1}

\noindent A positive answer to Question \ref{Q:CanWeClassifyCountableIfix} would yield a complete classification of the countable $\mathcal{L}$-structures satisfying infinity that can be realised as $\mathcal{I}_\mathrm{fix}(j)$ for some non-trivial automorphism $j$ of a model of $\mathrm{MOST}$. A result proved by John Hutchinson in \cite[Theorem 3.1]{hut76} shows that Question \ref{Q:CanWeClassifyCountableIfix} has a positive answer if both $\mathrm{MOST}+\Delta_0^\mathcal{P}\textrm{-collection}$ and $\mathrm{MOST}$ are replaced by $\mathrm{ZFC}$. \footnote{Hutchinson's result was generalized to models of countable cofinality by the second author \cite[Theorem 5.1]{kau83}.}

\begin{Theorems1} \label{Th:ToplessEndExtensionsOfModelsOfZFC}
(Hutchinson) If $\mathcal{M}= \langle M, \in^\mathcal{M} \rangle$ is countable with $\mathcal{M} \models \mathrm{ZFC}$ then there is a countable $\mathcal{L}$-structure $\mathcal{N}$ such that $\mathcal{M} \prec_\mathrm{topless} \mathcal{N}$. \Square   
\end{Theorems1}

\noindent Note that if $\mathcal{M} \models \mathrm{ZFC}$ and $\mathcal{M} \prec_e \mathcal{N}$ then for all $\alpha \in \mathrm{Ord}^\mathcal{M}$, $V_\alpha^\mathcal{M}=V_\alpha^\mathcal{N}$, and so $\mathcal{M} \prec_e^\mathcal{P} \mathcal{N}$. Therefore Theorem \ref{Th:ToplessEndExtensionsOfModelsOfZFC} shows that every countable model of $\mathrm{ZFC}$ can be realised as an $H$-cut of a model of $\mathrm{ZFC}$. Combined with the construction in Section \ref{Sec:RealisingHcutsAsIFix} this shows:

\begin{Theorems1}
Let $\mathcal{M}=\langle M, \in^\mathcal{M} \rangle$ be countable with $\mathcal{M} \models \mathrm{ZFC}$. There exists a countable $\mathcal{L}$-structure $\mathcal{M} \prec \mathcal{N}$ and an automorphism $j: \mathcal{N} \longrightarrow \mathcal{N}$ such that $\mathcal{I}_\mathrm{fix}(j)= \mathcal{M}$. \Square
\end{Theorems1}

\section[An extension of Togha's Theorem]{An extension of Togha's Theorem}

Togha~\cite[Theorem 3]{tog04} proves the following set-theoretic analogue of a result, due to Smory\'{n}ski \cite[Theorem A]{smo82}, about automorphisms of countable recursively saturated models of $\mathrm{PA}$.

\begin{Theorems1} \label{Th:ToghaTheorem}
(Togha) Let $\mathcal{M}$ be a countable recursively saturated model of $\mathrm{ZFC}$. If $I$ is an ordinal cut of $\mathcal{M}$ closed under cardinal exponentiation (Definition \ref{Df:OrdinalCutClosedUnderExponentiation}) then there is an automorphism $j: \mathcal{M} \longrightarrow \mathcal{M}$ such that $I$ is the largest initial segment of $\mathrm{Ord}^\mathcal{M}$ that is pointwise fixed by $j$. \Square  
\end{Theorems1}

\noindent In light of the correspondence revealed by Lemmas \ref{Th:HCutDeterminedByOrdinalCut} and \ref{Th:UnionOfHs}, Togha's Theorem shows that if $\mathcal{M}$ is a countable recursively saturated model of $\mathrm{ZFC}$ and $I \subseteq M$ is an $H$-cut of $\mathcal{M}$ then there is $j \in \mathrm{Aut}(\mathcal{M})$ such that $\mathcal{I}_\mathrm{fix}(j)=\langle I, \in^\mathcal{M}\rangle$. In this section we will generalise Togha's Theorem by showing that if $I$ is an $H$-cut of a countable recursively saturated model $\mathcal{M}$ of $\mathrm{ZFC}$ then $\mathcal{M}$ is endowed with
continuum-many automorphisms $j$ with the property that $\mathcal{I}_\mathrm{fix}(j)$ is exactly $\langle I, \in^\mathcal{M} \rangle$. This generalisation of Togha's Theorem is analogous to the generalisation of Smory\'{n}ski's result proved by the first author in \cite[Theorem B]{ena06}.    

\begin{Theorems1} \label{Th:ExtensionOfToghasTheorem}
Let $\mathcal{M}=\langle M, \in^\mathcal{M} \rangle$ be a countable recursively saturated model of $\mathrm{ZFC}$. Let $I \subseteq M$ be an $H$-cut of $\mathcal{M}$. There is an embedding $j \mapsto \check{j}$ of $\mathrm{Aut}(\mathbb{Q})$ into $\mathrm{Aut}(\mathcal{M})$ such that if $j \in \mathrm{Aut}(\mathbb{Q})$ is non-trivial then $\mathcal{I}_\mathrm{fix}(\check{j})= \langle I, \in^\mathcal{M} \rangle$. 
\end{Theorems1}

Theorem \ref{Th:ExtensionOfToghasTheorem} can also be viewed as a refinement of a result due to Schmerl \cite{sch85} which shows that if $\mathcal{M}$ is a countable recursively saturated model of $\mathrm{ZFC}$ (or indeed any theory equipped with enough coding) then there is a group embedding of $\mathrm{Aut}(\mathbb{Q})$ into $\mathrm{Aut}(\mathcal{M})$.\\
\\ 
\indent It should be noted that Lemma \ref{Th:HCutDeterminedByOrdinalCut} shows that any $H$-cut is uniquely determined by its ordinal spine, which forms an ordinal cut that is closed under exponentiation and contains $\omega^\mathcal{M}$. Therefore Theorem \ref{Th:ExtensionOfToghasTheorem} is equivalent to Theorem C mentioned in the abstract.\\
\\  
\indent Before proving Theorem \ref{Th:ExtensionOfToghasTheorem} we first need to prove a model-theoretic result that allows us to expand cofinal elementary extensions of models of $\mathrm{ZFC}$. We show that if $\mathcal{M}= \langle M, \in^\mathcal{M}, X^\mathcal{M} \rangle$ is a structure with $X^\mathcal{M} \subseteq M$ such that $\mathcal{M}$ satisfies $\mathrm{ZFC}$ plus the full collection scheme in the language of $\mathcal{M}$, and $\mathcal{N}$ is an $\mathcal{L}$-structure with $\langle M, \in^\mathcal{M} \rangle \prec_{\mathrm{cf}} \mathcal{N}$, then $\mathcal{N}$ can be expanded to a structure $\mathcal{N}^\prime$ in the language of $\mathcal{M}$ that is a cofinal elementary extension of $\mathcal{M}$. This model-theoretic tool is the set-theoretic version of a result that was independently proved by Schmerl \cite[Theorem 1.2]{sch81} and Kotlarski \cite[Theorem 8]{kot83} for cofinal elementary extensions of models of $\mathrm{PA}$. We use $\mathcal{L}_X$ to denote the language obtained by extending $\mathcal{L}$ with a new unary predicate $X$. Recall that $\mathrm{ZFC}(X)$ is obtained from $\mathrm{ZFC}$ by adding the schemes of $\mathcal{L}_X$-separation and $\mathcal{L}_X$-collection. 

\begin{Theorems1} \label{Th:CofinalExpansionLemma}
Let $\mathcal{M}= \langle M, \in^\mathcal{M}, X^\mathcal{M} \rangle$ be an $\mathcal{L}_X$-structure with $\mathcal{M} \models \mathrm{ZFC}(X)$. If $\mathcal{N}= \langle N, \in^\mathcal{N} \rangle$ is such that $\langle M, \in^\mathcal{M} \rangle \prec_{\mathrm{cf}} \mathcal{N}$ then there exists $X^\mathcal{M} \subseteq X^\mathcal{N} \subseteq N$ such that $\mathcal{M} \prec \langle N, \in^\mathcal{N}, X^\mathcal{N} \rangle$.    
\end{Theorems1}

\begin{proof}
Let $\mathcal{N}= \langle N, \in^\mathcal{N} \rangle$ be such that $\langle M, \in^\mathcal{M} \rangle \prec_{\mathrm{cf}} \mathcal{N}$. Note that this immediately implies that $\mathcal{N} \models \mathrm{ZFC}$. We begin by defining $X^\mathcal{N}$. Define $F: \mathrm{Ord}^\mathcal{M} \longrightarrow M$ such that for all $\alpha \in \mathrm{Ord}^\mathcal{M}$, 
$$\mathcal{M} \models \forall y(y \in  F(\alpha) \iff y \in V_\alpha \land X(y)).$$
The fact that $\mathcal{M} \models \mathrm{ZFC}(X)$ ensures that $F(\alpha)$ exists for all $\alpha \in \mathrm{Ord}^\mathcal{M}$. Also note that since $M \subseteq N$, for all $\alpha \in \mathrm{Ord}^\mathcal{M}$, $F(\alpha) \in N$.  Define $X^\mathcal{N} \subseteq N$ by: for all $x \in N$,
$$x \in X^\mathcal{N} \textrm{ if and only if there exists } \alpha \in \mathrm{Ord}^\mathcal{M} \textrm{ such that } \mathcal{N} \models (x \in F(\alpha)).$$
The fact that $\langle M, \in^\mathcal{M} \rangle \prec_{\mathrm{cf}} \mathcal{N}$, implies that $X^\mathcal{M} \subseteq X^\mathcal{N}$.\\
We will prove by induction on $n \in \omega$ that $\mathcal{M} \prec_n \langle N, \in^\mathcal{N}, X^\mathcal{N} \rangle$.\\
Firstly, note that since $\langle M, \in^\mathcal{M} \rangle \prec \mathcal{N}$, for all $\alpha \in \mathrm{Ord}^\mathcal{M}$, $V_\alpha^\mathcal{M}= V_\alpha^\mathcal{N}$. Therefore, for all $\alpha \in \mathrm{Ord}^\mathcal{M}$,
$$\mathcal{M} \models \forall x(x \in F(\alpha) \iff x \in V_\alpha \land X(x)) \textrm{ and } \mathcal{N} \models \forall x(x \in F(\alpha) \iff x \in V_\alpha \land X(x)).$$
And, for all $\Delta_0(\mathcal{L}_X)$-formulae $\phi(x_1, \ldots, x_m)$,
\begin{itemize}
\item[(I)] for all $a_1, \ldots, a_m \in N$, if $\alpha \in \mathrm{Ord}^\mathcal{M}$ is such that $a_1, \ldots, a_m \in (V_\alpha^\mathcal{N})^*$, then
$$\langle N, \in^\mathcal{N}, X^\mathcal{N} \rangle \models \phi(a_1, \ldots, a_m) \textrm{ if and only if } \mathcal{N} \models (\langle V_\alpha, \in, F(\alpha) \rangle \models \phi(a_1, \ldots, a_m)).$$ 
\item[(II)] for all $a_1, \ldots, a_m \in M$, if $\alpha \in \mathrm{Ord}^\mathcal{M}$ is such that $a_1, \ldots, a_m \in (V_\alpha^\mathcal{M})^*$, then
$$\mathcal{M} \models \phi(a_1, \ldots, a_m) \textrm{ if and only if } \langle M, \in^\mathcal{M} \rangle \models (\langle V_\alpha, \in, F(\alpha) \rangle \models \phi(a_1, \ldots, a_m)).$$
\end{itemize}
Let $\phi(x_1, \ldots, x_m)$ be a $\Delta_0(\mathcal{L}_X)$-formula and let $a_1, \ldots, a_m \in M$. Let $\alpha \in \mathrm{Ord}^\mathcal{M}$ be such that $a_1, \ldots, a_m \in (V_\alpha^\mathcal{M})^*$. Now,
$$\begin{array}{lcl}
\mathcal{M} \models \phi(a_1, \ldots, a_m) & \textrm{ if and only if } & \langle M, \in^\mathcal{M} \rangle \models (\langle V_\alpha, \in, F(\alpha) \rangle \models \phi(a_1, \ldots, a_m))\\
& \textrm{ if and only if } & \mathcal{N} \models (\langle V_\alpha, \in, F(\alpha) \rangle \models \phi(a_1, \ldots, a_m))\\
& \textrm{ if and only if } & \langle N, \in^\mathcal{N}, X^\mathcal{N} \rangle \models \phi(a_1, \ldots, a_m)
\end{array}.$$
This shows that $\mathcal{M} \prec_0 \langle N, \in^\mathcal{N}, X^\mathcal{N} \rangle$ and completes the base case of the induction. Now, let $n \in \omega$ and suppose that $\mathcal{M} \prec_n \langle N, \in^\mathcal{N}, X^\mathcal{N} \rangle$;
we prove this relation for $n+1$ in place of $n$.  We leave the case $n=0$ as an exercise, as it requires a simpler version of the argument below. So assume $n>0$.
Let $\phi(x_1, \ldots, x_m)$ be an $\mathcal{L}_X$-formula $\exists y \forall z \psi(y, z, x_1, \ldots, x_m)$ where $\psi(y, z, x_1, \ldots, x_m)$ is $\Sigma_{n-1}(\mathcal{L}_X)$. Let $a_1, \ldots, a_m \in M$. It follows immediately from the fact that $\mathcal{M} \prec_n \langle N, \in^\mathcal{N}, X^\mathcal{N} \rangle$ that if $\mathcal{M} \models \phi(a_1, \ldots, a_m)$, then $\langle N, \in^\mathcal{N}, X^\mathcal{N} \rangle \models \phi(a_1, \ldots, a_m)$. Conversely, suppose that
$$\langle N, \in^\mathcal{N}, X^\mathcal{N} \rangle \models \phi(a_1, \ldots, a_m) \textrm{ and } \mathcal{M} \models \neg \phi(a_1, \ldots, a_m).$$
Let $\alpha \in \mathrm{Ord}^\mathcal{M}$ be such that
$$\langle N, \in^\mathcal{N}, X^\mathcal{N} \rangle \models (\exists y \in V_\alpha) \forall z \psi(y, z, a_1, \ldots, a_m).$$
Note that
$$\mathcal{M} \models (\forall y \in V_\alpha) \exists z \neg \psi(y, z, a_1, \ldots, a_m).$$
Work inisde $\mathcal{M}$. Since $\mathcal{M} \models \mathrm{ZFC}(X)$, we can find a set $C$ such that for all $y \in V_\alpha$, there exists $z \in C$ such that $\neg\psi(y, z, a_1, \ldots, a_m)$ holds. Let $f: V_\alpha \longrightarrow C$ be such that for all $y \in V_\alpha$, $\neg \psi(y, f(y), a_1, \ldots, a_m)$ holds. Working in the meta-theory again, note that
$$\mathcal{M} \models (\forall y \in V_\alpha) \neg \psi(y, f(y), a_1, \ldots, a_m).$$
The expression $(\forall y \in V_\alpha)\neg \psi(y, f(y), a_1, \ldots, a_m)$ is $\Pi_n(\mathcal{L}_X)$ (indeed, $\Pi_{n-1}(\mathcal{L}_X)$) with parameters $a_1, \ldots, a_m, f \in M$. Therefore, by the induction hypothesis
$$\langle N, \in^\mathcal{N}, X^\mathcal{N} \rangle \models (\forall y \in V_\alpha) \neg \psi(y, f(y), a_1, \ldots, a_m).$$
Now, let $y_0 \in N$ be such that
$$\langle N, \in^\mathcal{N}, X^\mathcal{N} \rangle \models \forall z \psi(y_0, z, a_1, \ldots, a_m) \land (y_0 \in V_\alpha).$$
But,
$$\langle N, \in^\mathcal{N}, X^\mathcal{N} \rangle \models \neg \psi(y_0, f(y_0), a_1, \ldots, a_m),$$
which is a contradiction. This completes the induction step and proves the theorem.        
\Square
\end{proof}

We now turn to proving Theorem \ref{Th:ExtensionOfToghasTheorem}. For the remainder of this section fix a countable recursively saturated structure $\mathcal{M}= \langle M, \in^\mathcal{M} \rangle$ with $\mathcal{M} \models \mathrm{ZFC}$, and fix an $H$-cut $I \subseteq M$. Again, we will write $\mathcal{I}$ for the $\mathcal{L}$-structure $\langle I, \in^\mathcal{M} \rangle$. Let $\bar{\kappa} \in \mathrm{Ord}^\mathcal{M} \backslash \mathrm{Ord}^\mathcal{I}$ be such that $\mathcal{M} \models (\bar{\kappa} \textrm{ is a regular cardinal})$. Using the construction presented in Section \ref{Sec:RealisingHcutsAsIFix} we will construct $\mathcal{N}_{\mathcal{U}, \mathbb{Q}}$ such that $\mathcal{M} \prec_\mathrm{cf} \mathcal{N}_{\mathcal{U}, \mathbb{Q}}$ and there is an embedding $j \mapsto \check{j}$ of $\mathrm{Aut}(\mathbb{Q})$ into $\mathrm{Aut}(\mathcal{N}_{\mathcal{U}, \mathbb{Q}})$ such that for all fixed point free $j \in \mathbb{Q}$, $\mathcal{I}_{\mathrm{fix}}(\check{j})= \mathcal{I}$. Carefully choosing the ultrafilter $\mathcal{U}$ will allow us to build an isomorphism between $\mathcal{M}$ and $\mathcal{N}_{\mathcal{U}, \mathbb{Q}}$ that fixes $I$. In order to ensure that such an isomorphism exists we need $\mathcal{M}$ and $\mathcal{N}_{\mathcal{U}, \mathbb{Q}}$ to code the same subsets of $I$.

\begin{Definitions1}
Let $\mathcal{N}=\langle N, \in^\mathcal{N} \rangle$ be an $\mathcal{L}$-structure and let $K \subseteq N$. Define
$$\mathrm{SSy}_K(\mathcal{N})= \{c^* \cap K \mid c \in N\}.$$  
\end{Definitions1}

\begin{Definitions1}
A filter $\mathcal{U} \subseteq (\mathcal{P}(\bar{\kappa}))^*$ is $I$-conservative if for every $n \in \omega$, and for every $f \in M$, if $\mathcal{M} \models (f: H_{\bar{\kappa}} \longrightarrow \mathcal{P}([\bar{\kappa}]^{n+1}))$ then there exists $X \in \mathcal{U}$ and an $I$-large $\mathcal{M}$-cardinal $\lambda \in \bar{\kappa}^*$ such that for all $x \in (H_{\lambda^+}^\mathcal{M})^*$,
$$\mathcal{M} \models ([X]^{n+1} \subseteq f(x)) \textrm{ or } \mathcal{M} \models ([X]^{n+1} \subseteq [\bar{\kappa}]^{n+1} \backslash f(x)).$$ 
\end{Definitions1}

As we did in Section \ref{Sec:RealisingHcutsAsIFix}, let
$$\mathcal{F}= \{f \in M \mid (\exists n \in \omega)(\mathcal{M} \models f \textrm{ is a function with domain } [\bar{\kappa}]^{n+1})\}$$
Let $\mathcal{L}_\mathcal{F}$ be the extension of $\mathcal{L}$ defined in Section \ref{Sec:RealisingHcutsAsIFix} that adds new function symbols $\hat{f}$ for every $f \in \mathcal{F}$. Let $\mathcal{M}_\mathcal{F}$ be the expansion of $\mathcal{M}$ to an $\mathcal{L}_\mathcal{F}$-structure defined in Section \ref{Sec:RealisingHcutsAsIFix}. Since each new function symbol $\hat{f}$ in $\mathcal{L}_\mathcal{F}$ is coded by a point in $\mathcal{M}$ and $\mathcal{M} \models \mathrm{ZFC}$, we immediately get the following extension of Lemma \ref{Th:SeparationAndCollectionInF}:

\begin{Lemma1}
$\mathcal{M}_\mathcal{F} \models \mathcal{L}_\mathcal{F}\textrm{-separation}+\mathcal{L}_\mathcal{F}\textrm{-collection}$. \Square
\end{Lemma1}

We will extend Theorem \ref{Th:UltrafilterExistence} to show that we find an external non-principle ultrafilter on the subsets of $\bar{\kappa}$ in $\mathcal{M}$ which is simultaneously $I$-complete, canonically Ramsey, $I$-tight, $I$-conservative, and contains arbitrarily small $I$-large sets. The fact that we can extend Theorem \ref{Th:UltrafilterExistence} will follow from the following lemma:

\begin{Lemma1} \label{Th:ConservativeExistence}
Let $n \in \omega$. Let $X \in (\mathcal{P}(\bar{\kappa})^\mathcal{M})^*$ be $I$-large and let $f \in M$ be such that $\mathcal{M} \models (f: H_{\bar{\kappa}} \longrightarrow \mathcal{P}([X]^{n+1}))$. There exists an $I$-large $\mathcal{M}$-cardinal $\lambda \in \bar{\kappa}^*$ and an $I$-large $Y \subseteq X$ such that for all $x \in (H_{\lambda^+}^\mathcal{M})^*$,
$$\mathcal{M} \models ([Y]^{n+1} \subseteq f(x)) \textrm{ or } \mathcal{M} \models ([Y]^{n+1} \subseteq [\bar{\kappa}]^{n+1} \backslash f(x)).$$  
\end{Lemma1}

\begin{proof}
Work inside $\mathcal{M}$. Let $\mu= |X|$. Using Lemma \ref{Th:BethIterate} and Theorem \ref{Th:ErdosRadoTheorem} we can find $I$-large cardinals $\gamma, \lambda < \mu$ such that $\gamma=|H_{\lambda^+}|$ and $\mu \rightarrow ((2^\gamma)^+)_{2^\gamma}^{n+1}$. Define $F: [X]^{n+1} \longrightarrow 2^{H_{\lambda^+}}$ such that for all $A \in [X]^{n+1}$,
$$F(A)=g_A \textrm{ where for all } x \in H_{\lambda^+},$$
$$g_A(x)=1 \textrm{ if and only if } A \in f(x).$$
Since $\mu \rightarrow ((2^\gamma)^+)_{2^\gamma}^{n+1}$, we can find an $I$-large $Y \subseteq X$ such that $F$ is monochromatic on $[Y]^{n+1}$. Therefore, for all $A, B \in [Y]^{n+1}$ and for all $x \in H_{\lambda^+}$, either $A, B \in f(x)$ or $A, B \notin f(x)$. Therefore, for all $x \in (H_{\lambda^+}^\mathcal{M})^*$,
$$\mathcal{M} \models ([Y]^{n+1} \subseteq f(x)) \textrm{ or } \mathcal{M} \models ([Y]^{n+1} \subseteq [\bar{\kappa}]^{n+1} \backslash f(x)).$$
\Square
\end{proof}

\begin{Theorems1} \label{Th:ConservativeUltrafilterExistence}
There exists an n.p. ultrafilter $\mathcal{U} \subseteq (\mathcal{P}(\bar{\kappa})^\mathcal{M})^*$ which is $I$-complete, canonically Ramsey, $I$-tight, $I$-conservative and such that $\{|X|^\mathcal{M} \mid X \in \mathcal{U}\}$ is downward cofinal in $\mathrm{Ord}^\mathcal{M} \backslash \mathrm{Ord}^\mathcal{I}$.
\end{Theorems1}

\begin{proof}
We use exactly the same method as we used to prove Theorem \ref{Th:UltrafilterExistence}. Let $\langle f_n \mid n \in \omega \rangle$ be an enumeration of $\mathcal{F}$ and let $\langle k_n \mid n \in \omega \rangle$ be a sequence of natural numbers such that for all $n \in \omega$, $\mathcal{M} \models (f_n \textrm{ is a function with domain } [\bar{\kappa}]^{k_n})$. Let $\langle \lambda_n \mid n \in \omega \rangle$ be a decreasing sequence of $\mathcal{M}$-cardinals that is downward cofinal in $\mathrm{Ord}^\mathcal{M} \backslash \mathrm{Ord}^\mathcal{I}$ with $\lambda_0 \in \bar{\kappa}^*$. Let 
$$\mathcal{G}=\{g \in M \mid (\exists n \in \omega)(\mathcal{M} \models g: H_{\bar{\kappa}} \longrightarrow \mathcal{P}([\bar{\kappa}]^{n+1}))\}.$$
Let $\langle g_n \mid n \in \omega \rangle$ be an enumeration of $\mathcal{G}$ and let $\langle l_n \mid n \in \omega \rangle$ be a sequence of natural numbers such that for all $n \in \omega$, $\mathcal{M} \models (g_n:H_{\bar{\kappa}} \longrightarrow \mathcal{P}([\bar{\kappa}]^{l_n}))$. Using Lemmas \ref{Th:CompleteExistence}, \ref{Th:fCanonicalExistence}, \ref{Th:TightExistence} and \ref{Th:ConservativeExistence} inductively build sequences $\langle W_n \mid n \in \omega \rangle$, $\langle Q_n \mid n \in \omega \rangle$, $\langle X_n \mid n \in \omega \rangle$, $\langle Y_n \mid n \in \omega \rangle$ and $\langle Z_n \mid n \in \omega \rangle$ of $I$-large elements of $(\mathcal{P}(\bar{\kappa})^\mathcal{M})^*$ such that for all $n \in \omega$,
\begin{enumerate}
\item $\mathcal{M} \models (W_n \supseteq Q_n \supseteq X_n \supseteq Y_n \supseteq Z_n \supseteq W_{n+1})$,
\item $\mathcal{M} \models (W_n \textrm{ is } f_n\textrm{-canonical})$,
\item there exists an $I$-large $\mathcal{M}$-cardinal $\lambda \in \bar{\kappa}^*$ such that for all $x \in (H_{\lambda^+}^\mathcal{M})^*$,
$$\mathcal{M} \models ([Q_n]^{l_n} \subseteq g_n(x)) \textrm{ or } \mathcal{M} \models ([Q_n]^{l_n} \subseteq [\bar{\kappa}]^{l_n} \backslash g_n(x)),$$ 
\item if $k_n= 1$ and there is an $\mathcal{M}$-cardinal $\mu \in I$ such that $\mathcal{M} \models (f_n: \bar{\kappa} \longrightarrow H_\mu)$ then $f_n$ is constant on $X_n$, otherwise $X_n= Q_n$, 
\item $\mathcal{M} \models (f_n \textrm{ is constant on } [Y_n]^{k_n})$ or there is an $I$-large $\mathcal{M}$-cardinal $\mu \in \bar{\kappa}^*$ such that 
$$\mathcal{M} \models (\forall A \in [Y_n]^{k_n})(f_n(A) \notin H_\mu),$$
\item $\mathcal{M} \models (|Z_n| < \lambda_n)$.  
\end{enumerate}
Define $\mathcal{U}= \{X \in (\mathcal{P}(\bar{\kappa})^\mathcal{M})^* \mid (\exists n \in \omega)(\mathcal{M} \models W_n \subseteq X)\}$. It is clear from the construction that $\mathcal{U}$ is an n.p. ultrafilter that is $I$-complete, canonically Ramsey, $I$-tight, $I$-conservative and is such that $\{|X|^\mathcal{M} \mid X \in \mathcal{U}\}$ is downward cofinal in $\mathrm{Ord}^\mathcal{M} \backslash \mathrm{Ord}^\mathcal{I}$.
\Square
\end{proof}

Let $\mathcal{U} \subseteq (\mathcal{P}(\bar{\kappa})^\mathcal{M})^*$ be an n.p. ultrafilter obtained from Theorem \ref{Th:ConservativeUltrafilterExistence}, so $\mathcal{U}$ is $I$-complete, canonically Ramsey, $I$-tight, $I$-conservative and $\{|X|^\mathcal{M} \mid X \in \mathcal{U}\}$ is downward cofinal in $\mathrm{Ord}^\mathcal{M} \backslash \mathrm{Ord}^\mathcal{I}$. Let $\mathcal{N}_{\mathcal{U}, \mathbb{Q}}= \langle N_{\mathcal{U}, \mathbb{Q}}, \in^\mathcal{N} \rangle$ be the iterated ultrapower, constructed in Section \ref{Sec:RealisingHcutsAsIFix}, of $\mathcal{M}$ and $\mathcal{U}$ endowed with a class of indiscernibles of order-type $\mathbb{Q}$. The results proved in Section \ref{Sec:RealisingHcutsAsIFix} and the fact $\mathcal{M} \models \mathrm{ZFC}$ imply that $\mathcal{N}_{\mathcal{U}, \mathbb{Q}}$ has the following properties:
\begin{enumerate}
\item $\mathcal{M} \prec_{\mathrm{cf}} \mathcal{N}_{\mathcal{U}, \mathbb{Q}}$,
\item $|N_{\mathcal{U}, \mathbb{Q}}|= \aleph_0$,
\item there is an embedding $j \mapsto \tilde{j}$ of $\mathrm{Aut}(\mathbb{Q})$ into $\mathrm{Aut}(\mathcal{N}_{\mathcal{U}, \mathbb{Q}})$ such that if $j \in \mathrm{Aut}(\mathbb{Q})$ has no fixed points then $\mathcal{I}_\mathrm{fix}(\tilde{j})= \mathcal{I}$.
\end{enumerate}
We will show below (Theorem \ref{Th:IfixIsI2}) that (3) holds for all non-trivial $j \in \mathrm{Aut}(\mathbb{Q})$, not just the $j \in \mathrm{Aut}(\mathbb{Q})$ with no fixed points. First, however, we show that $\mathcal{N}_{\mathcal{U}, \mathbb{Q}}$ can be identified with $\mathcal{M}$ by an isomorphism that fixes $I$. The existence of this isomorphism will follow from
Theorem~\ref{theorem-6.8}. The fact that $\mathcal{U}$ is $I$-conservative ensures that $\mathrm{SSy}_I(\mathcal{N}_{\mathcal{U}, \mathbb{Q}})= \mathrm{SSy}_I(\mathcal{M})$.

\begin{Lemma1}
$\mathrm{SSy}_I(\mathcal{N}_{\mathcal{U}, \mathbb{Q}})= \mathrm{SSy}_I(\mathcal{M})$.
\end{Lemma1}

\begin{proof}
Since $\mathcal{I} \subseteq_{\mathrm{end}} \mathcal{N}_{\mathcal{U}, \mathbb{Q}}$, it follows that
$$\mathrm{SSy}_I(\mathcal{N}_{\mathcal{U}, \mathbb{Q}}) \supseteq \mathrm{SSy}_I(\mathcal{M}).$$
We need to show the reverse inclusion. Let $[\hat{f}(c_{i_0}, \ldots, c_{i_n})] \in N_{\mathcal{U}, \mathbb{Q}}$ where $f \in \mathcal{F}$ and $i_0 < \cdots < i_n \in \mathbb{Q}$. Let
$$A= \{x \in I \mid \mathcal{N}_{\mathcal{U}, \mathbb{Q}} \models (x \in [\hat{f}(c_{i_0}, \ldots, c_{i_n})])\}.$$
So $A \in \mathrm{SSy}_I(\mathcal{N}_{\mathcal{U}, \mathbb{Q}})$ and we need to show that $A \in \mathrm{SSy}_I(\mathcal{M})$. Working inside $\mathcal{M}_\mathcal{F}$ define $g: H_{\bar{\kappa}} \longrightarrow \mathcal{P}([\bar{\kappa}]^{n+1})$ such that for all $x \in H_{\bar{\kappa}}$,
$$g(x)= \{\{\alpha_0 < \cdots < \alpha_n\} \in [\bar{\kappa}]^{n+1} \mid x \in \hat{f}(\alpha_0, \ldots, \alpha_n)\}.$$
So, $g \in M$. From Lemma \ref{Th:LosLemma} we have for all $x \in H_{\bar{\kappa}}^*$,
$$\mathcal{N}_{\mathcal{U}, \mathbb{Q}} \models (x \in [\hat{f}(c_{i_0}, \ldots, c_{i_n})]) \textrm{ if and only if there exists } Y \in \mathcal{U} \textrm{ s.t. } \mathcal{M}_\mathcal{F} \models ([Y]^{n+1} \subseteq g(x)).$$
Let $X \in \mathcal{U}$ and let $\lambda \in \bar{\kappa}^*$ be an $I$-large $\mathcal{M}$-cardinal such that for all $x \in (H_{\lambda^+}^\mathcal{M})^*$,
$$\mathcal{M}_\mathcal{F} \models ([X]^{n+1} \subseteq g(x)) \textrm{ or } \mathcal{M}_\mathcal{F} \models ([X]^{n+1} \subseteq [\bar{\kappa}]^{n+1} \backslash g(x)).$$
Working inside $\mathcal{M}_\mathcal{F}$, let 
$$C=\{x \in H_{\bar{\kappa}} \mid [X]^{n+1} \subseteq g(x)\}.$$
So, $C \in M$. Since $I \subseteq (H_{\lambda^+}^\mathcal{M})^*$, it follows that for all $x \in I$,
$$x \in C^* \textrm{ if and only if } \mathcal{N}_{\mathcal{U}, \mathbb{Q}} \models (x \in [\hat{f}(c_{i_0}, \ldots, c_{i_n})]).$$
Therefore $A= C^* \cap I$ and $A \in \mathrm{SSy}_I(\mathcal{M})$.  
\Square
\end{proof}

The next result shows that for any set in $\mathcal{N}_{\mathcal{U}, \mathbb{Q}}$ there is a set in $\mathcal{M}$ with the same $\mathcal{M}$-members on some initial segment of $\mathcal{M}$ that contains $I$.

\begin{Lemma1} \label{Th:GoodApproximations}
Let $\lambda \in M$ be an $I$-large $\mathcal{M}$-cardinal. If $u \in (H_\lambda^{\mathcal{N}_{\mathcal{U}, \mathbb{Q}}})^*$ then there is an $I$-large $\mathcal{M}$-cardinal $\mu \in M$ and $w \in M$ such that:
\begin{equation} \label{eq:GoodApproximationsCondition}
\{x \in M \mid (\mathcal{M} \models (x \in H_\mu)) \land (\mathcal{N}_{\mathcal{U}, \mathbb{Q}} \models ([\hat{h}_x(c_0)] \in u))\}= \{x \in M \mid \mathcal{M} \models (x \in w)\}.
\end{equation}
\end{Lemma1}

\begin{proof}
Let $u \in (H_\lambda^{\mathcal{N}_{\mathcal{U}, \mathbb{Q}}})^*$. Working inside $\mathcal{N}_{\mathcal{U}, \mathbb{Q}}$, define $f: \mathrm{Card} \cap \lambda \longrightarrow H_\lambda$ such that for all $\gamma \in \mathrm{Card} \cap \lambda$,
$$f(\gamma)= u \cap H_\gamma.$$
Since $\mathrm{SSy}_I(\mathcal{N}_{\mathcal{U}, \mathbb{Q}})= \mathrm{SSy}_I(\mathcal{M})$, there is $g \in M$ such that $g$ agrees with $f$ on $I$. Work inside $\mathcal{N}_{\mathcal{U}, \mathbb{Q}}$. Define
\begin{displaymath}
\begin{split}
C=\{ \delta \in \lambda \mid & (\delta \in \mathrm{Card}) \land \\
& (\forall \gamma \in \delta^+) ( (\gamma \in \mathrm{Card}) \longrightarrow ( (g \cap (2^\gamma \times H_{2^\gamma}) \textrm{ is a function}) \land (f(\gamma)=g(\gamma))))\}.
\end{split}
\end{displaymath}
Since $I$ is a $H$-cut of $\mathcal{N}_{\mathcal{U}, \mathbb{Q}}$, there is an $I$-large cardinal $\eta \in C$. Using Lemma \ref{Th:MCardinalsDownwardCofinalInN} we can find $\mu \in \mathrm{Card}^\mathcal{M} \backslash \mathrm{Card}^\mathcal{I}$ such that $\mathcal{N}_{\mathcal{U}, \mathbb{Q}} \models ([\hat{h}_\mu] \leq \eta)$. Since $\mathcal{M} \prec \mathcal{N}_{\mathcal{U}, \mathbb{Q}}$,
$$\mathcal{M} \models (g \cap (2^\mu \times H_{2^\mu}) \textrm{ is a function}).$$
Let $w \in M$ be such that $\mathcal{M} \models (g(\mu)=w)$. Therefore
$$\mathcal{N}_{\mathcal{U}, \mathbb{Q}} \models (g([\hat{h}_\mu(c_0)])= [\hat{h}_w(c_0)]) \land ([\hat{h}_w(c_0)] = u \cap H_{[\hat{h}_\mu(c_0)]}).$$
And so $w, \mu \in M$ satisfy (\ref{eq:GoodApproximationsCondition}).  
\Square
\end{proof}

We are now in a position to show that $\mathcal{M}$ and $\mathcal{N}_{\mathcal{U}, \mathbb{Q}}$ can be identified by an isomorphism that fixes $I$.

\begin{Theorems1}\label{theorem-6.8}
There exists an isomorphism $\Theta: \mathcal{M} \longrightarrow \mathcal{N}_{\mathcal{U}, \mathbb{Q}}$ such that for all $x \in I$, $\Theta(x)=x$.
\end{Theorems1}

\begin{proof}
We will construct $\Theta$ using a back-and-forth construction. We begin by endowing $\mathcal{M}$ and $\mathcal{N}_{\mathcal{U}, \mathbb{Q}}$ with satisfaction classes. By Theorem \ref{Th:SatisfactionClassesForRecSatModels} there is an $X^\mathcal{M} \subseteq M$ such that 
\begin{itemize}
\item[(I)] $\langle M, \in^\mathcal{M}, X^\mathcal{M} \rangle \models \mathrm{ZFC}(X)$, 
\item[(II)] $X^\mathcal{M}$ is a satisfaction class for $\mathcal{M}$.
\end{itemize}
Since $\mathcal{M}$ is recursively saturated, there is a non-standard $s \in (\omega^\mathcal{M})^* (=(\omega^\mathcal{I})^*= (\omega^{\mathcal{N}_{\mathcal{U}, \mathbb{Q}}})^*)$ such that $X^{\mathcal{M}}$ is $s$-correct for $\mathcal{M}$. Throughout this proof we will identify formulae from the point of view of $\mathcal{M}$ with their G\"{o}del codes in $\omega^\mathcal{M}$. We will also abbrieviate our notation by identifying elements of $M$ with the equivalence classes of their corresponding constant functions in $N_{\mathcal{U}, \mathbb{Q}}$. By shortening $X^\mathcal{M}$ if necessary, we can assume without loss of generality that $\langle M, \in^\mathcal{M}, X^\mathcal{M} \rangle$ satisfies:
\begin{equation}\label{eq:XisShort}
\langle \phi, a \rangle \in X \Rightarrow \ulcorner\phi\urcorner < s. 
\end{equation}
Using Theorem \ref{Th:CofinalExpansionLemma} we can expand $\mathcal{N}_{\mathcal{U}, \mathbb{Q}}$ to an $\mathcal{L}_X$-structure $\langle N_{\mathcal{U}, \mathbb{Q}}, \in^{\mathcal{N}_{\mathcal{U}, \mathbb{Q}}}, X^{\mathcal{N}_{\mathcal{U}, \mathbb{Q}}}\rangle$ such that 
$$\langle M, \in^\mathcal{M}, X^\mathcal{M} \rangle \prec \langle N, \in^{\mathcal{N}_{\mathcal{U}, \mathbb{Q}}}, X^{\mathcal{N}_{\mathcal{U}, \mathbb{Q}}}\rangle.$$
It immediately follows that $X^{\mathcal{N}_{\mathcal{U}, \mathbb{Q}}}$ is a satisfaction class that is $s$-correct for $\mathcal{N}_{\mathcal{U}, \mathbb{Q}}$ and $\mathcal{N}_{\mathcal{U}, \mathbb{Q}}$ is recursively saturated. Moreover, $\langle N_{\mathcal{U}, \mathbb{Q}}, \in^{\mathcal{N}_{\mathcal{U}, \mathbb{Q}}}, X^{\mathcal{N}_{\mathcal{U}, \mathbb{Q}}}\rangle$ satisfies (\ref{eq:XisShort}). Let $\langle p_i \mid i \in \omega \rangle$ be an enumeration of $M$ and let $\langle q_i \mid i \in \omega \rangle$ be an enumeration of $N_{\mathcal{U}, \mathbb{Q}}$. We will construct $\Theta: \mathcal{M} \longrightarrow \mathcal{N}_{\mathcal{U}, \mathbb{Q}}$ by constructing sequences $\langle u_i \mid i \in \omega \rangle$ and $\langle v_i \mid i \in \omega \rangle$, together with decreasing sequences $\langle \gamma_i \mid i \in \omega \rangle$ and $\langle r_i \mid i \in \omega \rangle$, such that
\begin{itemize}
\item[(I)] $\langle u_i \mid i \in \omega \rangle$ enumerates $M$,
\item[(II)] $\langle v_i \mid i \in \omega \rangle$ enumerates $N$,
\item[(III)] for all $i \in \omega$, $\gamma_i$ is an $I$-large $\mathcal{M}$-cardinal,
\item[(IV)] for all $i \in \omega$, $r_i \in (\omega^\mathcal{M})^*$ is non-standard. 
\end{itemize}
We then define 
$$\Theta(u_i)= v_i \textrm{ for all } i \in \omega.$$
We begin by letting $\gamma_0 \in \bar{\kappa}^*$ be an $I$-large $\mathcal{M}$-cardinal and $r_0=s$. The fact that $\langle M, \in^\mathcal{M}, X^\mathcal{M} \rangle \prec \langle N, \in^{\mathcal{N}_{\mathcal{U}, \mathbb{Q}}}, X^{\mathcal{N}_{\mathcal{U}, \mathbb{Q}}}\rangle$ ensures that if $\phi(y) \in (\mathrm{Form}^\mathcal{M})^*$ with $\phi(y) < r_0$ then for all $a \in (H_{\gamma_0}^\mathcal{M})^*$,
\begin{equation} \label{eq:Dagger0}
\langle \phi(y), \langle a \rangle \rangle \in X^\mathcal{M} \textrm{ if and only if } \langle \phi(y), \langle a \rangle \rangle \in X^{\mathcal{N}_{\mathcal{U}, \mathbb{Q}}}.
\end{equation}
At stage $j > 0$, after having defined $u_0, \ldots, u_{j-1} \in M$, $v_0, \ldots, v_{j-1} \in N_{\mathcal{U}, \mathbb{Q}}$, $r_{j} \in (\omega^\mathcal{M})^*$ and $\gamma_{j} \in M$ we will ensure that the following condition is maintained:
\begin{itemize}
\item[]($\dagger_j$) if $\phi(x_0, \ldots, x_{j-1}, y) \in (\mathrm{Form}^\mathcal{M})^*$ with $\phi(x_0, \ldots, x_{j-1}, y) < r_j$ then for all $a \in (H_{\gamma_j}^\mathcal{M})^*$,
\begin{displaymath}
\begin{split}
& \langle \phi(x_0, \ldots, x_{j-1}, y), \langle u_0, \ldots, u_{j-1}, a \rangle \rangle \in X^\mathcal{M} \textrm{ if and only if} \\
& \langle \phi(x_0, \ldots, x_{j-1}, y), \langle v_0, \ldots, v_{j-1}, a \rangle \rangle \in X^{\mathcal{N}_{\mathcal{U}, \mathbb{Q}}}.
\end{split}
\end{displaymath}
\end{itemize}
We construct the sequences $\langle u_i \mid i \in \omega \rangle$, $\langle v_i \mid i \in \omega \rangle$, $\langle \gamma_i \mid i \in \omega \rangle$ and $\langle r_i \mid i \in \omega \rangle$ by induction. Note that (\ref{eq:Dagger0}) corresponds to the condition $\dagger_0$ that forms the base case of our induction before any of the $u$s or $v$s have been defined. Suppose that we have defined $u_0, \ldots, u_{2k-1} \in M$ and $v_0, \ldots, v_{2k-1} \in N_{\mathcal{U}, \mathbb{Q}}$, $r_{2k} \in (\omega^\mathcal{M})^*$ and $\gamma_{2k}$ such that $\dagger_{2k}$ holds.\\
\textbf{STAGE} $2k+1$: We need to choose $u_{2k} \in M$, $v_{2k} \in N$, $r_{2k+1} \in (\omega^\mathcal{M})^*$ and $\gamma_{2k+1}$ such that $\dagger_{2k+1}$ is maintained. Let $v_{2k}= q_k$. This choice will eventually ensure that $\Theta$ is onto $N_{\mathcal{U}, \mathbb{Q}}$. Using Lemma \ref{Th:BethIterate}, let $\lambda$ be an $I$-large $\mathcal{M}$-cardinal such that $2^{\lambda} < \gamma_{2k}$. Working inside $\langle N_{\mathcal{U}, \mathbb{Q}}, \in^{\mathcal{N}_{\mathcal{U}, \mathbb{Q}}}, X^{\mathcal{N}_{\mathcal{U}, \mathbb{Q}}}\rangle$, define
$$u= \{\langle \phi(x_0, \ldots, x_{2k}, y), a \rangle \mid (\phi < r_{2k}) \land (a \in H_\lambda) \land (\langle \phi, \langle v_0, \ldots, v_{2k}, a \rangle\rangle \in X) \}.$$ 
Note that $u \in (H_{2^\lambda}^{\mathcal{N}_{\mathcal{U}, \mathbb{Q}}})^*$. The fact that $X^{\mathcal{N}_{\mathcal{U}, \mathbb{Q}}}$ is $s$-correct ensures that for all $l \in \omega$,
$$\mathcal{N}_{\mathcal{U}, \mathbb{Q}} \models \exists x \left( (\forall y \in H_\lambda) \left( \bigwedge_{\phi < l} (\langle \phi(x_0, \ldots, x_{2k}, y), y \rangle \in u) \iff \phi(v_0, \ldots, v_{2k-1}, x, y) \right) \right).$$
Using Lemma \ref{Th:GoodApproximations} we can find an $I$-large $\mathcal{M}$-cardinal $\gamma_{2k+1} \in M$ with $\gamma_{2k+1} \leq \lambda$ and $w \in M$ such that
$$\{x \in M \mid (\mathcal{M} \models (x \in H_{\gamma_{2k+1}})) \land (\mathcal{N}_{\mathcal{U}, \mathbb{Q}} \models ([\hat{h}_x(c_0)] \in u))\}= \{x \in M \mid \mathcal{M} \models (x \in w)\}.$$
Therefore, for all $l \in \omega$,
$$\mathcal{N}_{\mathcal{U}, \mathbb{Q}} \models \exists x \left( (\forall y \in H_{\gamma_{2k+1}}) \left( \bigwedge_{\phi < l} (\langle \phi(x_0, \ldots, x_{2k}, y), y \rangle \in w \iff \phi(v_0, \ldots, v_{2k-1}, x, y)) \right) \right).$$
Since $2^{\gamma_{2k+1}} < \gamma_{2k}$ and $w, \gamma_{2k+1} \in (H_{2^{\gamma_{2k}}}^\mathcal{M})^*$, $\dagger_{2k}$ implies that for all $l \in \omega$,
$$\mathcal{M} \models \exists x \left( (\forall y \in H_{\gamma_{2k+1}}) \left( \bigwedge_{\phi < l} (\langle \phi(x_0, \ldots, x_{2k}, y), y \rangle \in w \iff \phi(u_0, \ldots, u_{2k-1}, x, y)) \right) \right).$$ 
For all $l \in \omega$, let $\psi_l(x, y, x_0, \ldots, x_{2k-1}, z)$ be the $\mathcal{L}$-formula:
$$\bigwedge_{\phi < l} (\langle \phi(x_0, \ldots, x_{2k}, y), y \rangle \in z \iff \phi(x_0, \ldots, x_{2k-1}, x, y)).$$
And let $\Psi(l, u_0, \ldots, u_{2k-1}, w, \gamma_{2k+1})$ be the $\mathcal{L}_X$-formula:
$$\exists x((\forall y \in H_{\gamma_{2k+1}})(\langle \psi_l(x, y, x_0, \ldots, x_{2k-1}, z), \langle x, y, u_0, \ldots, u_{2k-1}, w \rangle \rangle \in X)).$$
Therefore, for all $l \in \omega$,
$$\langle M, \in^\mathcal{M}, X^\mathcal{M} \rangle \models \Psi(l, u_0, \ldots, u_{2k-1}, w, \gamma_{2k+1}).$$
So, by overspill in $\langle M, \in^\mathcal{M}, X^\mathcal{M} \rangle$, we can find a non-standard $r_{2k+1} \in (\omega^\mathcal{M})^*$ with $r_{2k+1} < r_{2k}$ such that
$$\langle M, \in^\mathcal{M}, X^\mathcal{M} \rangle \models \exists x((\forall y \in H_{\gamma_{2k+1}}) (\langle \psi_{r_{2k+1}}(x, y, x_0, \ldots, x_{2k-1}, z), \langle x, y, u_0, \ldots, u_{2k-1}, w \rangle \rangle \in X)).$$
Let $u_{2k} \in M$ be such that
$$\langle M, \in^\mathcal{M}, X^\mathcal{M} \rangle \models (\forall y \in H_{\gamma_{2k+1}}) (\langle \psi_{r_{2k+1}}(x, y, x_0, \ldots, x_{2k-1}, z), \langle u_{2k}, y, u_0, \ldots, u_{2k-1}, w \rangle \rangle \in X).$$
Therefore the choices of $v_{2k}$, $u_{2k}$, $r_{2k+1}$ and $\gamma_{2k+1}$ made in this stage maintain $\dagger_{2k+1}$.\\   
\textbf{STAGE} $2k+2$: We need to choose $u_{2k+1} \in M$, $v_{2k+1} \in N_{\mathcal{U}, \mathbb{Q}}$, $r_{2k+2} \in (\omega^\mathcal{M})^*$ and $\gamma_{2k+2}$ such that $\dagger_{2k+2}$ is maintained. Let $u_{2k+1}= p_k$. This choice will eventually ensure that $\Theta$ is defined on all of $M$. Using Lemma \ref{Th:BethIterate}, let $\gamma_{2k+2} \in M$ be an $I$-large $\mathcal{M}$-cardinal such that $2^{\gamma_{2k+2}} < \gamma_{2k+1}$. Working inside $\langle M, \in^\mathcal{M}, X^\mathcal{M} \rangle$, define
$$u= \{\langle \phi(x_0, \ldots, x_{2k+1}, y), a \rangle \mid (\phi < r_{2k+1}) \land (a \in H_{\gamma_{2k+2}}) \land (\langle \phi, \langle u_0, \ldots, u_{2k+1}, a \rangle\rangle \in X) \}.$$ 
Note that $u \in (H_{2^{\gamma_{2k+2}}}^{\mathcal{M}})^*$. The fact that $X^{\mathcal{M}}$ is $s$-correct ensures that for all $l \in \omega$,
$$\mathcal{M} \models \exists x \left( (\forall y \in H_{\gamma_{2k+2}}) \left( \bigwedge_{\phi < l} (\langle \phi(x_0, \ldots, x_{2k+1}, y), y \rangle \in u) \iff \phi(u_0, \ldots, u_{2k}, x, y) \right) \right).$$
Since $u, \gamma_{2k+2} \in (H_{\gamma_{2k+1}}^\mathcal{M})^*$, $\dagger_{2k+1}$ implies that
$$\mathcal{N}_{\mathcal{U}, \mathbb{Q}} \models \exists x \left( (\forall y \in H_{\gamma_{2k+2}}) \left( \bigwedge_{\phi < l} (\langle \phi(x_0, \ldots, x_{2k+1}, y), y \rangle \in u) \iff \phi(v_0, \ldots, v_{2k}, x, y) \right) \right).$$
For all $l \in \omega$, let $\psi_l(x, y, x_0, \ldots, x_{2k}, z)$ be the $\mathcal{L}$-formula:
$$\bigwedge_{\phi < l} (\langle \phi(x_0, \ldots, x_{2k+1}, y), y \rangle \in z \iff \phi(x_0, \ldots, x_{2k}, x, y)).$$
And let $\Psi(l, v_0, \ldots, v_{2k}, u, \gamma_{2k+2})$ be the $\mathcal{L}_X$-formula:
$$\exists x((\forall y \in H_{\gamma_{2k+2}})(\langle \psi_l(x, y, x_0, \ldots, x_{2k}, z), \langle x, y, v_0, \ldots, v_{2k}, u \rangle \rangle \in X)).$$
Therefore, for all $l \in \omega$,
$$\langle N_{\mathcal{U}, \mathbb{Q}}, \in^{\mathcal{N}_{\mathcal{U}, \mathbb{Q}}}, X^{\mathcal{N}_{\mathcal{U}, \mathbb{Q}}}\rangle \models \Psi(l, v_0, \ldots, v_{2k}, u, \gamma_{2k+2}).$$
So, by overspill in $\langle N_{\mathcal{U}, \mathbb{Q}}, \in^{\mathcal{N}_{\mathcal{U}, \mathbb{Q}}}, X^{\mathcal{N}_{\mathcal{U}, \mathbb{Q}}}\rangle$, we can find a non-standard $r_{2k+2} \in (\omega^\mathcal{M})^*$ with $r_{2k+2} < r_{2k+1}$ such that
\begin{displaymath}
\begin{split}
\langle N_{\mathcal{U}, \mathbb{Q}}, \in^{\mathcal{N}_{\mathcal{U}, \mathbb{Q}}}, X^{\mathcal{N}_{\mathcal{U}, \mathbb{Q}}}\rangle \models \exists x((\forall y \in H_{\gamma_{2k+2}}) (\langle \psi_{r_{2k+2}}(x, y, x_0, \ldots, x_{2k}, z),\\
\langle x, y, v_0, \ldots, v_{2k}, u \rangle \rangle \in X)).
\end{split}
\end{displaymath}
Let $v_{2k+1} \in N_{\mathcal{U}, \mathbb{Q}}$ be such that
\begin{displaymath}
\begin{split}
\langle N_{\mathcal{U}, \mathbb{Q}}, \in^{\mathcal{N}_{\mathcal{U}, \mathbb{Q}}}, X^{\mathcal{N}_{\mathcal{U}, \mathbb{Q}}}\rangle \models (\forall y \in H_{\gamma_{2k+2}}) (\langle \psi_{r_{2k+2}}(x, y, x_0, \ldots, x_{2k-1}, z),\\
\langle v_{2k+1}, y, v_0, \ldots, v_{2k}, u \rangle \rangle \in X).
\end{split}
\end{displaymath}
Therefore the choices of $v_{2k+1}$, $u_{2k+1}$, $r_{2k+2}$ and $\gamma_{2k+2}$ made in this stage maintain $\dagger_{2k+2}$.\\ 
This completes the proof of the theorem.        
\Square
\end{proof}

Before concluding the proof of Theorem \ref{Th:ExtensionOfToghasTheorem} we first show that the fact that $\mathcal{U}$ is $I$-conservative can be used to demonstrate that Theorem \ref{Th:IFixIsI} can be strengthened to show that \emph{for every non-trivial } $j\in \mathrm{Aut}(\mathbb{Q})$, $\mathcal{I}_{\mathrm{fix}}(\tilde{j})= \mathcal{I}$.

\begin{Lemma1} \label{Th:LargeSetsOfNAreInhabitted}
If $[\tau] \in N_{\mathcal{U}, \mathbb{Q}} \backslash I$ then there exists $x \in M \backslash I$ with
$$\mathcal{N}_{\mathcal{U}, \mathbb{Q}} \models (\hat{h}_x(c_0) \in \mathrm{TC}(\{[\tau]\})).$$
\end{Lemma1}

\begin{proof}
Let $[\tau] \in N_{\mathcal{U}, \mathbb{Q}} \backslash I$. Let $f \in \mathcal{F}$ and $i_1 < \cdots < i_n \in \mathbb{Q}$ be such that 
$$\mathcal{N}_{\mathcal{U}, \mathbb{Q}} \models ([\hat{f}(c_{i_1}, \ldots, c_{i_n})]= \mathrm{TC}(\{[\tau]\})).$$
Let $\lambda \in M$ be an $I$-large $\mathcal{M}$-cardinal such that
$$\mathcal{N}_{\mathcal{U}, \mathbb{Q}} \models ([\hat{f}(c_{i_1}, \ldots, c_{i_n})] \in H_{[\hat{h}_\lambda(c_0)]}).$$
Using Lemma \ref{Th:GoodApproximations} we can find an $I$-large $\mathcal{M}$-cardinal $\mu \in M$ and $w \in M$ such that
$$\{x \in M \mid (\mathcal{M} \models (x \in H_\mu)) \land(\mathcal{N}_{\mathcal{U}, \mathbb{Q}} \models ([\hat{h}_x(c_0)] \in [\hat{f}(c_{i_1}, \ldots, c_{i_n})]))\}$$
$$=\{x \in M \mid \mathcal{M} \models (x \in w)\}.$$
Working inside $\mathcal{M}$, define
$$C= \{ \gamma \in \mu \mid (\gamma \in \mathrm{Card}) \land (w \cap (H_{\gamma^+} \backslash H_{\gamma}) \neq \emptyset)\}.$$
We will prove that $C \notin I$. Suppose that $C \in I$. Work inside $\mathcal{M}$. Let $\eta^\prime = \sup C$ and let $\eta = (\eta^\prime)^+$. Therefore $\eta \in I$ is a cardinal and for all cardinals $\eta \leq \gamma < \mu$, $w \cap (H_{\gamma^+} \backslash H_\gamma) = \emptyset$. Let $a \in w \backslash H_\eta$ be a set with minimal rank. Therefore, for all $y \in a$, $y \in H_\eta$. But $\mathcal{P}(H_\eta) \subseteq H_{2^{2^\eta}}$, and so there exists a cardinal $\eta \leq \gamma < 2^{2^\eta}$ such that $a \in H_{\gamma^+} \backslash H_{\gamma}$. This contradicts the fact that $w \cap (H_{\gamma^+} \backslash H_{\gamma}) = \emptyset$ and shows that $C \notin I$. Therefore, since $I \subseteq M$ is a $H$-cut, there exists an $I$-large $\mathcal{M}$-cardinal $\xi \in M$ with $\xi \in \mu^*$ such that $w \cap (H_{\xi^+} \backslash H_\xi) \neq \emptyset$. Let $x \in (w \backslash (H_{\xi^+} \backslash H_\xi))$. Therefore
$$\mathcal{N}_{\mathcal{U}, \mathbb{Q}} \models (\hat{h}_x(c_0) \in [\hat{f}(c_{i_1}, \ldots, c_{i_n})]).$$        
\Square
\end{proof}

\begin{Theorems1} \label{Th:IfixIsI2}
If $j \in \mathrm{Aut}(\mathbb{Q})$ is non-trivial then $\mathcal{I}_{\mathrm{fix}}(\tilde{j})= \mathcal{I}$.
\end{Theorems1} 

\begin{proof}
Let $j \in \mathrm{Aut}(\mathbb{Q})$ be non-trivial. The embedding of $\mathcal{I}$ into an initial segment of $\mathcal{N}_{\mathcal{U}, \mathbb{Q}}$ shows that $I \subseteq I_{\mathrm{fix}}(\tilde{j})$. Therefore, we need to show that for every $[\tau] \in N_{\mathcal{U}, \mathbb{Q}} \backslash I$, there exists $y \in (\mathrm{TC}(\{[\tau]\})^{\mathcal{N}_{\mathcal{U}, \mathbb{Q}}})^*$ such that $\tilde{j}(y) \neq y$. Let $l, m \in \mathbb{Q}$ such that $j(l) = m \neq l$. Let $[\tau] \in N_{\mathcal{U}, \mathbb{Q}} \backslash I$. Using Lemma \ref{Th:LargeSetsOfNAreInhabitted}, let $x \in M \backslash I$ be such that
$$\mathcal{N}_{\mathcal{U}, \mathbb{Q}} \models (\hat{h}_x(c_0) \in \mathrm{TC}(\{[\tau]\})).$$
Therefore, by Lemma \ref{Th:TCsArePreserved},
$$\mathcal{N}_{\mathcal{U}, \mathbb{Q}} \models (\hat{h}_{\mathrm{TC}(x)}(c_0) \subseteq \mathrm{TC}(\{[\tau]\})).$$
And, since $x \in M \backslash I$, $|\mathrm{TC}(x)|^\mathcal{M} \notin I$. Therefore, there is $X \in \mathcal{U}$ with $\mathcal{M} \models (|X| \leq  |\mathrm{TC}(x)|)$. Let $g \in M$ be such that
$$\mathcal{M} \models (g: X \longrightarrow \mathrm{TC}(x)) \land (g \textrm{ is injective}).$$
Therefore, there is $f \in M$ such that $\mathcal{M} \models (f: \bar{\kappa} \longrightarrow \mathrm{TC}(x))$ and for all $z_1 < z_2 \in X^*$,
$$\mathcal{M}_\mathcal{F} \models \hat{f}(z_1) \neq \hat{f}(z_2).$$
So, by Lemma \ref{Th:LosLemma}, $\tilde{j}([\hat{f}(c_l)]) \neq [\hat{f}(c_l)]$. Moreover, for all $z \in X^*$,
$$\mathcal{M}_\mathcal{F} \models (\hat{f}(z) \in \mathrm{TC}(x)).$$
Therefore, by Lemma \ref{Th:LosLemma},
$$\mathcal{N}_{\mathcal{U}, \mathbb{Q}} \models ([\hat{f}(c_l)] \in [\hat{h}_{\mathrm{TC}(x)}(c_0)]).$$
And so,
$$\mathcal{N}_{\mathcal{U}, \mathbb{Q}} \models ([\hat{f}(c_l)] \in \mathrm{TC}(\{[\tau]\})),$$
which proves the theorem.        
\Square
\end{proof} 

Let $\Theta: \mathcal{M} \longrightarrow \mathcal{N}_{\mathcal{U}, \mathbb{Q}}$ be the isomorphism obtained from Theorem \ref{theorem-6.8} such that for all $x \in I$, $\Theta(x)= x$. To complete the proof of Theorem \ref{Th:ExtensionOfToghasTheorem} we use $\Theta$ to turn the embedding $j \mapsto \tilde{j}$ of $\mathrm{Aut}(\mathbb{Q})$ into $\mathrm{Aut}(\mathcal{N}_{\mathcal{U}, \mathbb{Q}})$, into an embedding $j \mapsto \check{j}$ of $\mathrm{Aut}(\mathbb{Q})$ into $\mathrm{Aut}(\mathcal{M})$ such that if $j \in \mathrm{Aut}(\mathbb{Q})$ is non-trivial then $\mathcal{I}_\mathrm{fix}(\check{j})= \mathcal{I}$. For all $j \in \mathrm{Aut}(\mathbb{Q})$, define $\check{j}: \mathcal{M} \longrightarrow \mathcal{M}$ such that
$$\check{j}(x)= \Theta^{-1}(\tilde{j}(\Theta(x))) \textrm{ for all } x \in M.$$
It follows immediately from the fact that $\Theta$ is isomorphism that that the map $j \mapsto \check{j}$ is an injective group homomorphism of $\mathrm{Aut}(\mathbb{Q})$ into $\mathrm{Aut}(\mathcal{M})$.

\begin{Theorems1}
If $j \in \mathrm{Aut}(\mathbb{Q})$ is non-trivial then $\mathcal{I}_{\mathrm{fix}}(\check{j})= \mathcal{I}$. 
\end{Theorems1}

\begin{proof}
Let $j \in \mathrm{Aut}(\mathbb{Q})$ be non-trivial. The fact that $\Theta$ fixes $I$ immediately implies that $I \subseteq I_{\mathrm{fix}}(\check{j})$. We need to show that if $x \in M \backslash I$ then there exists $y \in (\mathrm{TC}(\{x\})^\mathcal{M})^*$ such that $\check{j}(y) \neq y$. Let $x \in M \backslash I$. Therefore $\Theta(x) \in N_{\mathcal{U}, \mathbb{Q}} \backslash I$. By Theorem \ref{Th:IfixIsI2} there is $y \in (\mathrm{TC}(\{\Theta(x)\})^{\mathcal{N}_{\mathcal{U}, \mathbb{Q}}})^*$ such that $\tilde{j}(y) \neq y$. Now, $\Theta^{-1}(y) \in (\mathrm{TC}(\{x\})^\mathcal{M})^*$ and
$$\check{j}(\Theta^{-1}(y))= \Theta^{-1}(\tilde{j}(\Theta(\Theta^{-1}(y))))= \Theta^{-1}(\tilde{j}(y)) \neq \Theta^{-1}(y).$$
\Square
\end{proof}

This completes the proof of Theorem \ref{Th:ExtensionOfToghasTheorem}.  

\bibliographystyle{alpha}
\bibliography{.}          

\end{document}